\documentclass[11pt]{article}
\usepackage{amsmath,amssymb,amsfonts} % Typical maths resource packages
\usepackage{graphics}                 % Packages to allow inclusion of graphics
\usepackage{makeidx}
\usepackage{graphicx}
\usepackage{sidecap}
\usepackage[caption=false,font=normalsize,labelfont=sf,textfont=sf]{subfig}
\usepackage{algorithmic}
\usepackage[algoruled,vlined]{algorithm2e}
\usepackage{fullpage}
\usepackage{float}

\numberwithin{equation}{section}
\numberwithin{figure}{section}
\numberwithin{table}{section}

\newtheorem{theorem}{Theorem}[section]
\newtheorem{proposition}[theorem]{Proposition}
\newtheorem{corollary}[theorem]{Corollary}
\newtheorem{lemma}[theorem]{Lemma}

\newtheorem{definition}[theorem]{Definition}

\let \e=\varepsilon
\let \f=\varphi

\let \th=\theta

\let \vt=\vartheta
\let \rr=\rho
\let \k=\kappa
\let \n=\nu

\let \w=\omega

\let \la=\lambda

\let \La=\Lambda
\let \g=\gamma
\let \za=\zeta

\let \o=\sigma
\def \a{\alpha}
\def \b{\beta}
\def \da{\delta}

\def \az{analytic}

\def \as{approximat}
\def \nn{\nonumber}
\def \pf{p-filter}

\def \fr{filter}

\def \fb{filter bank}
\def \pb{p-filter bank}
\def \sq{sequence}
\def \fg{frequency response}

\def \sy{symmetr}
\def \df{differen}

\def \ds{discrete spline}
\def \dss{discrete-spline}

\def \mg{magnitude response}
\def \ir{impulse response}

\def \srr{\stackrel{\Delta}{=}}

\def \pa{\paragraph}

\def \v{vector}

\def \ww{wavelet}

\def \btt{\begin{theorem}}
\def \bll{\begin{lemma}}
\def \bdd{\begin{definition} }
\def \bpp{\begin{proposition} }
\def \bcc{\begin{corollary} }

\def \ett{\end{theorem}}
\def \ell{\end{lemma}}
\def \edd{\end{definition}}
\def \epp{\end{proposition}}
\def \ecc{\end{corollary} }

\def \vm{vanishing moment}

\def \we{waveform}

\def \pr{processing}

\def \ss{signal}

\def \on{orthonormal}
\def \oo{orthogonal}

\def \ry{represent}

\def \cc{convolution}

\def \c{coefficient}
\def \fd{function}

\def \ff{frequency}

\def \aa{analysis}
\def \sa{synthesis}
\def \s{spline}
\def \ds{discrete spline}
\def \dss{discrete-spline}
\def \sp{spectr}
\def \de{direction}

\def \p{periodic}
\def \ir{impulse response}

\def \wq{wavelet packet}

\def \mv{modulation  matri}

\def \ehh{Equation \rf}
\def \eh{Eq. \rf}

\def \d{decomposition}
\def \r{reconstruction}

\def \dd{discrete}
\def \dt{discrete-time}

\def \t{transform}
\def \F{Fourier}

\newcommand{\rf}[1]{(\ref{#1})}

\newtheorem{thm}{Theorem}
\newtheorem{rmk}{Remark}[section]

\newtheorem{exm}{Example}

\def \srr{\stackrel{\mathrm{def}}{=}}

\def \bq{\begin{quote}}
\def \eq{\end{quote}}

\def \bs{B-spline}
\def \bt{\begin{thm}}
\def \et{\end{thm}}

\def \br{\begin{rmk}}
\def \er{\end{rmk}}

\def\proof{\goodbreak\medskip\noindent{\small\bf Proof: }}
\def\eop{{\vrule height7pt width7pt depth0pt}\par\bigskip}

\def \bex{\begin{exm}}
\def \eex{\end{exm}}

\def\eop{{\vrule height7pt width7pt depth0pt}\par\bigskip}

\begin{document}
\title{Analytic and directional  \wq s in the space of periodic signals}

\author{Amir Averbuch$^1$~~Pekka Neittaanm\"aki$^2$~~Valery  Zheludev$^1$  \\
$^1$School of Computer Science\\
Tel Aviv University, Tel Aviv 69978, Israel\\~\\
$^2$Department of Mathematical Information Technology\\
 University of Jyv\"askyl\"a, Finland}
 \date{ }
\maketitle
\begin{abstract}
The paper presents a versatile library of \az\ and quasi-\az\ complex-valued \wq s (WPs) which originate from \ds s of arbitrary orders. The real parts of the quasi-\az\ WPs are the regular \s-based \on\ WPs designed in \cite{ANZ_book3}. The imaginary parts are the so-called complementary \on\ WPs, which, unlike the \sy ic regular  WPs, they are  anti\sy ic. Tensor products of 1D quasi-\az\ WPs  provide a diversity of 2D WPs oriented in multiple directions. For example, a set of the fourth-level WPs comprises 62 \df t directions. The designed computational scheme in the paper enables us to get fast and easy implementation of the WP \t s. The presented WPs proved to be efficient in signal/image processing applications
%on \aa\ of vibration \ss s recorded from the rotating bearings and
such as restoration of images degraded by either additive noise or missing of up to 90\% of their pixels.
\end{abstract}

\section{Introduction}\label{sec:s1}
%\subsection{Related works}\label{sec:ss11}
Since the introduction  in  Kingsbury \cite{king1,king2} of  complex \ww\ \t s implemented by the dual-tree scheme, the complex \ww s (DT\_CW), \ww\ frames and \wq s (WPs) have become a field of active research that appears in multiple applications (\cite{jalob1,jalob2,jalob3,xie_wang_zhao_chen,barakin,bay_sele,bhan_zhao,bhan_zhao_zhu}, to name a few). The advantages of the DT\_CWs over the standard real \ww\ \t s stem from  their approximate shift invariance and some  directionality inherent to tensor-products of the DT\_CWs.

However, the \de ality of the DT\_CWs is very limited (only 6 directions) and this is a drawback for image \pr\ applications.  The tight tensor-product complex  \ww\ frames TP\_CTF$_{n}$ with \df t number of directions, are designed in \cite{bhan_zhao,bhan_zhao_zhu} and some of them, in particular TP\_CTF$_{6}$ and TP\_CTF$^{\downarrow}_{n}$, demonstrate excellent performance (in terms of PSNR) for image denoising and inpainting.  The number of directions in both 2D  TP\_CTF$_{6}$ and TP\_CTF$^{\downarrow}_{n}$ frames is 14 and remains the same for all \d\ levels.

Some of the disadvantages of the above 2D  TP\_CTF$_{6}$ and TP\_CTF$^{\downarrow}_{n}$ frames are mentioned in
\cite{che_zhuang}. For example, ``limited and fixed number of  directions is undesirable in practice especially
when the resolution of an image is very high that requires large number of directional filters in order to capture as many
features with different orientations as possible" (\cite{che_zhuang}).
 In addition, ``due to the fixed number
of 1D filters, the number of free parameters
is limited which prevents the search of optimal filter bank
systems for image processing" (\cite{che_zhuang}).

  According to \cite{che_zhuang}, the remedy for the above shortcomings is in the incorporation of the two-layer structure that is inherent in the TP\_CTF$_{6}$ and TP\_CTF$^{\downarrow}_{n}$ frames into directional \fb s introduced in \cite{bhan_zhuang,zhuang}.

The complex \wq s (Co\_WPs) is an alternative way to overcome the above disadvantages. The first version of complex WPs appears in \cite{jalob1}  after the introduction of the complex \ww s by Kingsbury. The complex \ww\ \t s in \cite{jalob1,jalob2,jalob3} are extended to the Co\_WP \t s by the application of the same \fr s as used in the DT\_CW \t s
to the high-\ff\ bands. Although  the low- and high-\ff\ bands in DT\_CW are \as ely \az, this is not the case for the Co\_WPs designed in \cite{jalob1,jalob2,jalob3}.
In addition, as shown in \cite{bay_sele} (Fig. 1), much energy passes into the negative half-bands of the \sp a.  Another approach to the design of  Co\_WPs is described in \cite{bay_sele}. It is suggested  in \cite{bay_sele} that in order to retain an \as e \az ity of the dual-tree WP \t s, the \fb s for the second decomposition  of the \t s should be the same for both stems of the tree.

%\subsection{Our motivation and approach}\label{sec:ss12}
Although the potential advantages of the Co\_WP \t s are apparent, so far, to the best of our knowledge, none of the Co\_WP schemes  in the literature have the desired properties such as perfect \ff\ separation,
   Hilbert \t\ relation between real and imaginary parts of the Co\_WPs,
  orthonormality of shifts of real and imaginary parts of the Co\_WPs,
unlimited number of \de s in the multidimensional case,
a variety of free parameters,
and fast and easy implementation.

Our motivation in this paper is to fill this gap. We design a family of Co\_WP \t s which possess all the above properties. As a base for the design, we use the family of  \dt\ WPs originated from \p\ \ds s of different orders that are described in \cite{ANZ_book3} (Chapter 4). The \wq s $\psi_{[m],l}^{2r}$, where $m$ is the \d\ level, $l=0,...,2^{m}-1$ is the index of the related \ff\ band and $2r$ is the order of the generating \ds, are \sy ic, well localized in time domain (although are not compactly supported), their DFTs \sp a are flat, and provide a refined split of the \ff\ domain. The WP \t s are executed in the \ff\ domain using the Fast \F\ \t\ (FFT). By varying the order $2r$, we can supply the WPs $\psi_{[m],l}^{2r}$ with any number of local \vm s without increase of the computational cost. Different combinations of the shifts in these WPs provide a variety  of \on\   bases of the space of $N$-\p\ \ss s.

To derive the Co\_WPs, we expand the WPs $\psi_{[m],l}^{2r}$ to \p\ \az\ \dt\ \ss s $\bar{\psi}_{\pm[m],l}^{2r}=\psi_{[m],l}^{2r}\pm i\,\th_{[m],l}^{2r}$, where $\th_{[m],l}^{2r}$ is the \dd\ \p\ Hilbert \t\ (HT) of the WP $\psi_{[m],l}^{2r}$. The \we s $\th_{[m],l}^{2r}$ are anti\sy ic, and for all   $l\neq0, \,2^{m}-1,$  orthonormal properties similar to the properties of the WPs $\psi_{[m],l}^{2r}$ take place. To achieve  \on ity, the \we s $\th_{[m],l}^{2r}, \;l=0, \,2^{m}-1$ are slightly corrected at the expense of minor deviation from anti\sy y and we get a new \on\ complimentary  WP (cWP)  system $\left\{\f_{[m],l}^{2r}\right\}, \;m=1,...,M,\;l=0, ...,2^{m}-1$, where for $l\neq0, \,2^{m}-1$, the WPs satisfy $\f_{[m],l}^{2r}=\th_{[m],l}^{2r}$. The magnitude \sp a of the cWPs $\f_{[m],l}^{2r}$ coincide with the \sp a of the respective WPs $\psi_{[m],l}^{2r}$ and, similarly to the WPs $\psi_{[m],l}^{2r}$, the cWPs $\f_{[m],l}^{2r}$ provide a variety  of \on\   bases for the space of $N$-\p\ \ss s.
%%%%%???

Correspondingly, we define the quasi-\az\ WP systems (qWP) as $$  \Psi^{2r}_{\pm[m],l}=\psi^{2r}_{[m],l}   \pm i\f^{2r}_{[m],l}, \quad m=1,...,M,\;l=0,...,2^{m}-1,$$ where all the WPs with indices other than $l=0,2^{m}-1$ are \az.
For the implementation of the \t s with the complex qWPs we  do not use the dual-tree scheme with \df t \fb s for real and imaginary \ww s but use the scheme with a single complex \fb\ in  the first step of the \t, and a real \fb\ in  the additional steps.

A dual-tree structure type appears in the 2D case when two sets of qWPs are defined as the tensor products of 1D qWPs

\begin{equation}\label{psipsi}
\Psi_{++[m],j,l}^{2r}[k,n] \srr \Psi_{+[m],\k}^{2r}[k]\,\Psi_{+[m],l}^{2r}[n], \quad
  \Psi_{+-[m],j,l}^{2r}[k,n] \srr \Psi_{+[m],\k}^{2r}[k]\,\Psi_{-[m],l}^{2r}[n]
\end{equation}
and \pr\ with the qWPs $\Psi_{+\pm[m],j,l}^{2r}$ is executed separately.

The real and imaginary parts of the qWPs $\Psi_{+\pm[m],j,l}^{2r}$ are the 2D \we s oriented in multiple \de s, specifically the $2(2^{m+1}-1)$ directions at the $m$-th \d\ level. Such an abundant \de ality proved to be efficient in the examples on image denoising and inpainting. It is worth mentioning that the WPs of  one- and two-dimensions have a localized oscillating  structure, which is useful for detection of transient oscillating events in 1D \ss s and oscillating patterns in the images (for example, ``Barbara" in Fig. \ref{barb30}).

Both one- and two-dimensional \t s are implemented in a very  fast ways by using FFT.

The paper is organized as follows:
Section \ref{sec:s2} outlines briefly \p\ \dt\ WPs originated from \ds s and corresponding \t s that form a basis for the design of Co\_WPs. The \aa\ $\tilde{\mathbf{F}}$ and \sa\ $\mathbf{F}$ \fb s for the WP \t s are described.
  Section \ref{sec:s3} outlines the construction of \dt\ \p\ \az\ \ss s. This section also introduces complimentary sets of WPs (cWPs), \az\ and quasi-\az\ WPs (qWPs).
   Section \ref{sec:s4} describes the implementation of the cWP and qWP \t s. The \fb s for   one step  of  \aa\ and \sa\ \t s are introduced. It is interesting to note that subsequent application of the direct and inverse qWP \fb s to a \ss\  $\mathbf{x}$ produces the \az\ \ss\ $\bar{\mathbf{x}}=\mathbf{x}+i\,H(\mathbf{x})$. All the subsequent steps of cWP and qWP \t s are implemented with the same \fb s $\tilde{\mathbf{F}}$ and $\mathbf{F}$ as used in the above WP \t s (section \ref{sec:s2}).
      Section \ref{sec:s5} extends the design of 1D qWPs to the 2D case. The 2D qWPs are defined via tensor products as shown in \eh{psipsi}. Directionality of the 2D qWPs is discussed.
 Section \ref{sec:s6} describes the implementation of the 2D qWP \t s by a dual-tree.
 Section \ref{sec:s7} presents a few experimental results for images restoration degraded by either strong additive noise or by  missing many of the pixels. In one example, both the degradation sources are present.
 Section \ref{sec:s8} discusses the results.
The Appendix contains proof of a proposition.

\paragraph{Notations and abbreviations}
$N=2^{j}$, $M=2^{m},\,m<j$, $N_{m}=2^{j-m}$ and $\Pi[N]$ is a space of real-valued  $N$-\p\ \ss s.
$\Pi[N,N]$ is the space of two-dimensional   $N$-\p\ arrays in both vertical and horizontal directions.
$\w\srr e^{2\pi\,i/N}$.

A \ss\ $\mathbf{x}=\left\{ x[k]\right\}\in\Pi[N]$ is \ry ed by its  inverse \dd\ \F\ \t\ (DFT)
\begin{equation}\label{fsx}
\begin{array}{lll}%\label{fsx}
  x[k]&=&\frac{1}{N}\sum_{n=0}^{N-1}\hat{x}[n]\,\w^{kn}=\frac{1}{N}\sum_{n=-N/2}^{N/2-1}\hat{x}[n]\,\w^{kn},\\%\nn
 \hat{x}[n] &=&\sum_{k=0}^{N-1}{x}[k]\,\w^{-kn},\quad \hat{x}[-n]=\hat{x}[N-n]= \hat{x}[n]^{\ast},
\end{array}
\end{equation}
where $\cdot^{\ast}$ means  complex conjugate. In particular,
\(
  \hat{x}[0]=\sum_{k=0}^{N-1}{x}[k]\)  and \(  \hat{x}[N/2]=\sum_{k=0}^{N-1}(-1)^{k}\,{x}[k]\) are real numbers.

The DFT of an $N_{m}$-\p\ \ss\ is $\hat{x}[n]_{m} =\sum_{k=0}^{N_{m}-1}{x}[k]\,\w^{-kn2^{m}}$.
The abbreviation PR means perfect \r.
The abbreviations 1D and 2D mean one-dimensional and two-dimensional, respectively. FFT is the fast \F\ \t, HT is the Hilbert \t,  $H(\mathbf{x})$ is the \dd\ \p\ HT of a \ss\ $\mathbf{x}$.

The abbreviations WP, cWP and qWP mean \wq s (typically \s-based \wq s $\psi^{2r}_{[m],l}$), complimentary \wq s $\f^{2r}_{[m],l}$ and quasi-\az\ \wq s $\Psi^{2r}_{\pm[m],l}$, respectively, in 1D case, and  \wq s $\psi^{2r}_{[m],j,l}$, complimentary \wq s $\f^{2r}_{[m],j,l}$ and quasi-\az\ \wq s $\Psi^{2r}_{+\pm[m],l,j}$, respectively, in 2D case.

\begin{equation}\label{u2r}
  U^{4r}[n]\srr \frac{1}{2}\left(\cos^{4r}\frac{\pi\,n}{N} +\sin^{4r}\frac{\pi\,n}{N}\right).
\end{equation}

Peak Signal-to-Noise ratio (PSNR) in decibels (dB) is
\[
  PSNR\srr10\log_{10}\left(\frac{N\,255^2}{\sum_{k=1}^N(x_{k}-\tilde
  x_{k})^2}\right)\; dB.
\]
SBI stands for split Bregman iterations and \pf\ means \p\ \fr.

\section{Outline of orthonormal WPs originated from \ds s: preliminaries}\label{sec:s2}
This section provides a brief  outline of  \p\ \dt\ \wq s originated from \ds s and corresponding \t s. For details see Chapter 4 in \cite{ANZ_book3}.
\subsection{ Periodic  \ds s }\label{sec:ss21}

The centered span-two $N$-\p\ \dd\ \bs\ of order $2r$ is defined as the IDFT of the \sq\
\begin{eqnarray*}\label{bs_2}
 \hat{b}^{2r}[n]&=&\cos^{2r}\frac{\pi\,n}{N},\quad {b}^{2r}[k]=\frac{1}{N}\sum_{n=-N/2}^{N/2-1}\w^{kn}\, \cos^{2r}\frac{\pi\,n}{N}.
\end{eqnarray*}
The \bs s are non-negative \sy ic  finite-length \ss s (up to periodization). Only the samples ${b}^{2r}[k], \;k=-r,...r,$ are non-zero.

The \ss s $
  {s}^{2r}[k]\srr \sum_{l=0}^{N/2-1}q[l]\,{b}^{2r}[k-2l],
$
which are referred to as \ds s,  form an $N/2$-dimensional subspace $^{2r}{\mathcal{S}}_{[1]}^{0}$ of the space  $\Pi[N]$ whose basis consists of two-sample shifts of the \bs\ $\mathbf{b}^{2r}$.
Here $\mathbf{q}=\left\{ q[l]\right\}, \;l=0, ...  ,N/2-1,$ is a real-valued \sq. The DFT of the \ds\  $\mathbf{{s}}^{2r}$ is
\begin{eqnarray*}
\label{adsf1}
 \hat{ {s}}^{2r}[n] =\hat{q}[n]_{1}\, \hat{ {b}}^{2r}[n]  =\hat{q}[n]_{1}\,\cos^{2r}\frac{\pi\,n}{N}.
\end{eqnarray*}
A \ds\  $\psi_{[1],0}^{2r}\in{} ^{2r}{\mathcal{S}}_{[1]}^{0}$ is defined through its inverse DFT (iDFT):
\begin{equation*}\label{psi0}
  \psi_{[1],0}^{2r}[k]\srr\frac{1}{N}\sum_{n=-N/2}^{N/2-1}\w^{kn}\,\frac{\cos^{2r}\frac{\pi\,n}{N}}{\sqrt{ U^{4r}[n]}},
\end{equation*} where  $ U^{4r}[n]$ is defined in \eh{u2r}.
\bpp[\cite{ANZ_book3}, Proposition 3.6]\label{pro:psi0}
Two-sample shifts  $\left\{ \psi_{[1],0}^{2r}[\cdot-2l]\}\right\}, \;l=0, ...  ,N/2-1,$ of the \ds s $ \psi_{[1],0}^{2r}$ form an \on\ basis of the subspace $^{2r}{\mathcal{S}}_{[1]}^{0}\subset\Pi[N]$.

The \oo\ projection of a \ss\  $\mathbf{x}\in \Pi[N]$  onto the space   ${} ^{2r}{\mathcal{S}}_{[1]}^{0}$ is the \ss\ $\mathbf{x}_{[1]}^{0}\in\Pi[N]$ such that
\begin{eqnarray}\nn
          x_{[1]}^{0}[k]&=&\sum_{l=0}^{N/2-1}y_{[1]}^{0}[l]\, \psi_{[1],0}^{2r}[k-2l] ,\quad
             y_{[1]}^{0}[l]  =\left\langle \mathbf{x},\, \psi_{[1],0}^{2r}[\cdot-2l]   \right\rangle
       = \sum_{k=0}^{N-1}h_{[1]}^{0}[k-2l] \,x[k], \\\label{y0_del}  h_{[1]}^{0}[k] &=&\psi_{[1],0}^{2r}[k],\quad
         \hat{{h}}_{[1]}^{0}[n] = \hat{ \psi}_{[1],0}^{2r}[n]=\frac{\cos^{2r}\frac{\pi\,n}{N}}{\sqrt{U^{4r}[n]}}.
         \end{eqnarray}
\epp
\subsection{ Orthogonal complement for   subspace $ ^{2r}{\mathcal{S}}_{[1]}^{0}$  }\label{sec:ss22}
The  \oo\ complement for   ${} ^{2r}{\mathcal{S}}_{[1]}^{0}$ in the \ss\ space $\Pi[N]$ is denoted by
   ${} ^{2r}{\mathcal{S}}_{[1]}^{1}$. Thus, $\Pi[N]={\mathcal{S}}_{[0]}={} ^{2r}{\mathcal{S}}_{[1]}^{0}\bigoplus^{2r}{\mathcal{S}}_{[1]}^{1}$. The \on\ basis in
the subspace        is formed by the two-sample shifts  $\left\{ \psi_{[1],1}^{2r}[\cdot-2l]\}\right\}, \;l=0, ...  ,N/2-1,$ of the \ds\ $ \psi_{[1],1}^{2r}$, which is  defined through its inverse DFT (iDFT):
\begin{equation*}\label{psi1}
  \psi_{[1],1}^{2r}[k]\srr\frac{1}{N}\sum_{n=-N/2}^{N/2-1}\w^{kn}\,\frac{\w^{n}\,\sin^{2r}\frac{\pi\,n}{N}}{\sqrt{ U^{4r}[n]}}.
\end{equation*}

\bpp[\cite{ANZ_book3}, Proposition 4.1]\label{pro:psi1}
The \oo\ projection of a \ss\  $\mathbf{x}\in \Pi[N]$  onto the space   ${} ^{2r}{\mathcal{S}}_{[1]}^{1}$ is the \ss\ $\mathbf{x}_{[1]}^{1}\in\Pi[N]$ such that
\begin{eqnarray}\nn
          x_{[1]}^{1}[k]&=&\sum_{l=0}^{N/2-1}y_{[1]}^{1}[l]\, \psi_{[1],1}^{2r}[k-2l] ,\quad
          y_{[1]}^{1}[l]  =\left\langle \mathbf{x},\, \psi_{[1],1}^{2r}[\cdot-2l]   \right\rangle
       = \sum_{k=0}^{N-1}h_{[1]}^{1}[k-2l] \,x[k], \\\label{y1_del}  h_{[1]}^{1}[k] &=&\psi_{[1],1}^{2r}[k],\quad
         \hat{{h}}_{[1]}^{1}[n] = \hat{ \psi}_{[1],1}^{2r}[n]=\frac{\w^{n}\,\sin^{2r}\frac{\pi\,n}{N}}{\sqrt{ U^{4r}[n]}}.
         \end{eqnarray}
\epp
The \ss s $\psi_{[1],0}^{2r}$ and $\psi_{[1],1}^{2r}$ are referred to as the \dss\ \wq s of order $2r$ from the first \d\ level . They are the \ir s of the low- and high-pass \p\ \fr s (\pf s) $\mathbf{h}_{[1]}^{0}$ and $\mathbf{h}_{[1]}^{1}$, respectively.
\br\label{0_N2rem}We emphasise that the DFTs $\hat{\psi}_{[1],0}^{2r}[N/2]=0$ and $\hat{\psi}_{[1],1}^{2r}[0]=0$.\er

\begin{figure}[H]
\begin{center}
\resizebox{12cm}{5cm}{
\includegraphics{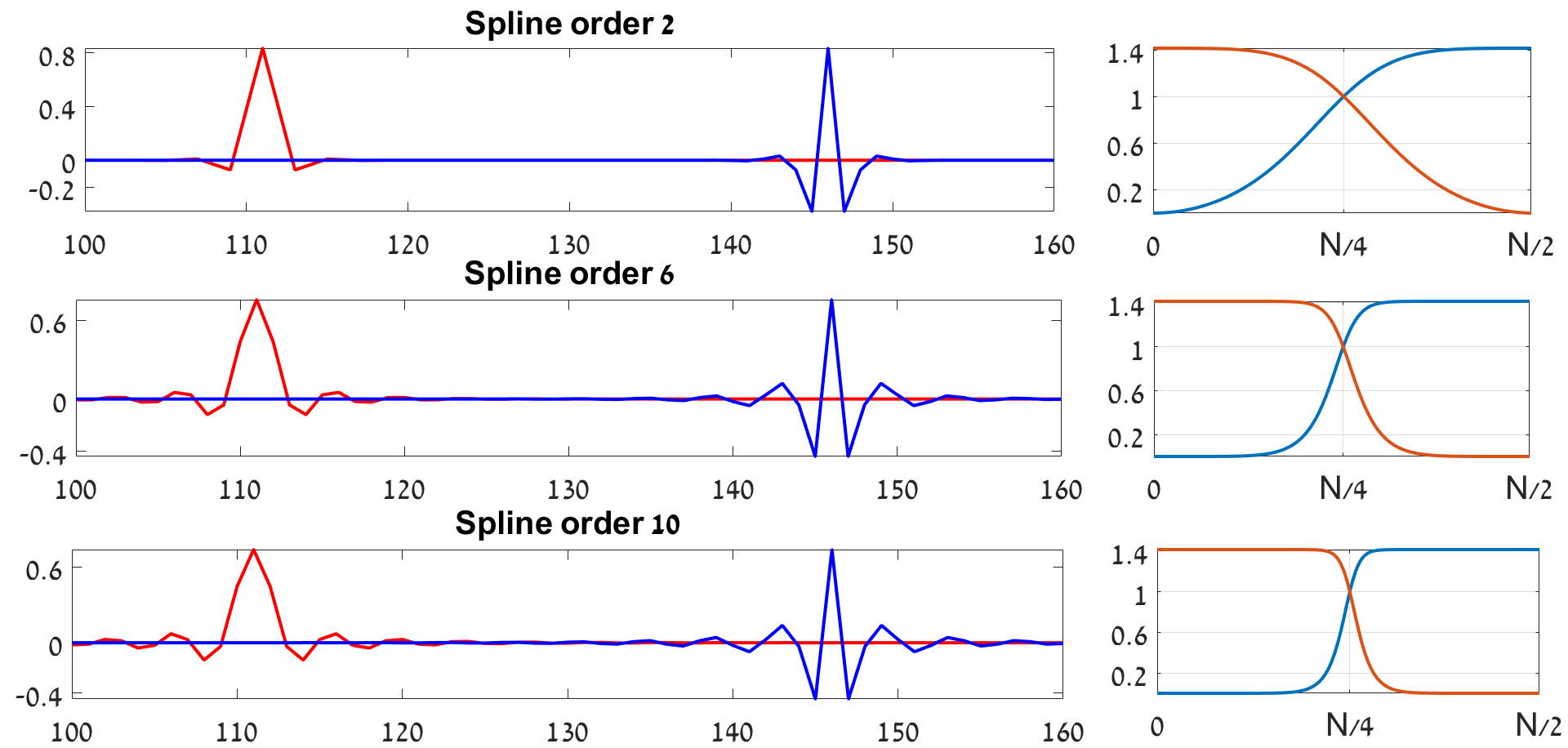}
}
\end{center}
\caption{ Left:  \wq s  $  \psi_{[1],0}^{2r}$ (red lines) and $  \psi_{[1],1}^{2r}$  (blue lines), $r=1,3,5$.
Right: magnitude \sp a of  $  \psi_{[1],0}^{2r}$ (red  lines) and $  \psi_{[1],1}^{2r}$  (blue lines)}
     \label{ds_wav1LA}
 \end{figure}

Figure \ref{ds_wav1LA} displays the \dss\ \wq s  $  \psi_{[1],0}^{2r}$ and $  \psi_{[1],1}^{2r}$  of different orders
and magnitudes of their DFT \sp a (which are the \pf s $\mathbf{h}_{[1],0}$ and $\mathbf{h}_{[1],1}$ \mg s).
It is seen that the \ww s are well localized in time domain. The \sp a are flat and their shapes tend to rectangular as their orders increase.

\subsection{One-level \wq\ \t\ of a \ss}\label{sec:ss23}

   The \t\ of a \ss\ $\mathbf{x}\in\Pi[N]$ into the pair $\left\{\mathbf{ y}_{[1]}^{0},\mathbf{ y}_{[1]}^{1}\right\}$ of \ss s from $\Pi[N/2]$ is
   referred to as the one-level \wq\ \t\ (WPT) of the \ss\ $\mathbf{x}$. According to Propositions \ref{pro:psi0} and \ref{pro:psi1},
   the \t\ is implemented by \fr ing $\mathbf{x}$ with time-reversed half-band low- and high-pass \pf s ${\mathbf{h}}_{[1]}^{0}$
    and ${\mathbf{h}}_{[1]}^{1}$, respectively, which is followed by downsampling. Thus the \pf s ${\mathbf{h}}_{[1]}^{0}$
    and ${\mathbf{h}}_{[1]}^{1}$ form a critically sampled \aa\ \pb\ $\tilde{\mathbf{H}}_{[1]}^{1}$.  Eqs. \rf{y0_del} and \rf{y1_del} imply that  its
  modulation matrix  is
  \begin{equation}\label{aa_modma10}%\label{alphabet0}
  \begin{array}{lll}
 \tilde{\mathbf{M}}[n]&=& \left(
                         \begin{array}{cc}
                            \hat{h}_{[1]}^{0}[n] &   \hat{h}_{[1]}^{0}\left[n+\frac{N}{2}\right]\\
                            \ \hat{h}^{1}_{[1]}[n] &   \hat{h}^{1}_{[1]}\left[n+\frac{N}{2}\right]\\
                         \end{array}
                       \right)=\left(
                                 \begin{array}{cc}
                                   \b[n] & \b\left[n+\frac{N}{2}\right]  \\
                                 \a[n] &\a\left[n+\frac{N}{2}\right]  \\
                                 \end{array}
                               \right)=\left(
                                 \begin{array}{cc}
                                   \b[n] & \w^{-n}\a[n]  \\
                                 \a[n] &-\w^{n}\b[n] \\
                                 \end{array}
                               \right),\\   %\label{alphabet1}
    \b[n]      &=&      \frac{\cos^{2r}\frac{\pi\,n}{N}}{\sqrt{U^{4r}[n]}}, \quad
       \a[n]=\w^{n}\,\b\left[n+\frac{N}{2}\right]=\w^{n}\, \frac{\sin^{2r}\frac{\pi\,n}{N}}{\sqrt{U^{4r}[n]}}.
  \end{array}
  \end{equation}

The \aa\ \mv x $\tilde{\mathbf{M}}[n]/\sqrt{2}$, as well as the matrix $\tilde{\mathbf{M}}[-n]/\sqrt{2}$ are unitary matrices. Therefore, the \sa\ \mv x
  is
  \begin{eqnarray}\label{sa_modma10}
    \mathbf{M}[n]=\left(
                                 \begin{array}{cc}
                                   \b[n] &\a[n]   \\
                             \w^{-n}\a[n]     &-\w^{n}\b[n]  \\
                                 \end{array}
                               \right)= \tilde{\mathbf{M}}[n]^{T}.
  \end{eqnarray}
  Consequently, the \sa\ \pb\ coincides with the \aa\ \pb\ and, together, they  form a perfect  reconstruction (PR) \pb .

  The one-level WP \t\ of a \ss\ $\mathbf{x}$ and its inverse are \ry ed in a matrix  form:
   \begin{equation}\label{mod_repAS1}
    \left(
     \begin{array}{c}
       \hat{y}_{[1]}^{0}[n]_{1} \\
         \hat{y}_{[1]}^{1}[n]_{1}\\
     \end{array}
   \right)=\frac{1}{2}\tilde{\mathbf{M}}[-n]\cdot \left(
     \begin{array}{l}
      \hat{x}[n] \\
       \hat{x}[\vec{n}]
     \end{array}
   \right),\quad  \left(
     \begin{array}{l}
      \hat{x}[n] \\
       \hat{x}[\vec{n}]
     \end{array}
   \right)={\mathbf{M}}[n]\cdot \left(
     \begin{array}{c}
        \hat{y}_{[1]}^{0}[n]_{1} \\
         \hat{y}_{[1]}^{1}[n]_{1}\\
     \end{array}
   \right)
    \end{equation}
    where $\vec{n}=n+{N}/{2}$.

\subsection{Extension of \t s to deeper \d\ levels}\label{sec:ss24}
\subsubsection{Second-level \pb s}\label{sec:sss241}
The \ss s $\mathbf{y}_{[1]}^{\la},\;\la=0,1,$ belong to the space $\Pi[N/2]\subset\Pi[N]$.  The space $\Pi[N/2]$ can be split into mutually \oo\ subspaces, which we denote by $\Pi^{0}[N/2]$  and $\Pi^{1}[N/2]$,
  in a way that is similar to the split of the space  $\Pi[N]$. Projection of a \ss\  $\mathbf{Y}\in \Pi[N/2]$ onto these subspaces and the inverse operation are done using  the \aa\ and \sa\ \pb s
 $\tilde{\mathbf{H}}_{[2]}=\left\{\mathbf{h}^{0}_{[2]},{\mathbf{h}}^{1}_{[2]}\right\}=\mathbf{H}_{[2]}$ (Eq. \ref{pfs2}), which operate in the space $\Pi[N/2]$. The \fg s of the \pf s are
 \begin{eqnarray}\label{pfs2}
  \hat{{h}}_{[2]}^{0}[n]_{1}= \b[2n] \quad
   \hat{{h}}_{[2]}^{1}[n]_{1} = \a[2n],
\end{eqnarray}
where $\b[n]$ and $\a[n]$ are defined in \eh{aa_modma10}. The \ir s of the \pf s $\mathbf{h}^{0}_{[2]}$ and $\mathbf{h}^{1}_{[2]}$ are \oo\ to each other in the space $\Pi[N/2]$ and their 2-sample shifts are mutually \oo\
\begin{equation*}\label{or_pf2}
 \sum_{k=0}^{N/2-1}{h}_{[2]}^{\la}[k-2l] \,{h}_{[2]}^{\mu}[k-2p] =\da[\la-\mu]\,\da[l-p],\quad \la,\mu=0,1.
\end{equation*}
The \oo\ projections of  a \ss\  $\mathbf{Y}\in \Pi[N/2]$ onto the subspaces $\Pi^{0}[N/2]$  and $\Pi^{1}[N/2]$ are
\begin{equation*}\label{ort_proj2}
         Y^{\mu}[k]=\sum_{l=0}^{N/4-1}y_{[2]}^{\mu}[l]\, {h}_{[2]}^{\mu}[k-2l]  , \quad
         y_{[2]}^{\mu}[l] = \sum_{k=0}^{N/2-1}{h}_{[2]}^{\mu}[k-2l] Y[k],
         \end{equation*}
         where $\mu=0,1.$
The \mv ces of the     \pb\
${\mathbf{H}}_{[2]}$ are
 \begin{eqnarray}
 \tilde{\mathbf{M}}_{[2]}[n]=\tilde{\mathbf{M}}[2n],\quad
    {\mathbf{M}}_{[2]}[n]={\mathbf{M}}[2n]\label{sa_modma20},
  \end{eqnarray}
  where the \mv ces $\tilde{\mathbf{M}}[n]$ and ${\mathbf{M}}[n]$ are defined in Eqs. \rf{1a_modma10} and \rf{sa_modma10}, respectively.
\subsubsection{Second-level WPTs}\label{sec:sss242}
By the application of the \aa\ \pb\ $\tilde{\mathbf{H}}_{[2]}$ (section \ref{sec:sss241} and Eq. \ref{pfs2}) to the \ss s $
 y_{[1]}^{\la}[k] = \sum_{n=0}^{N-1}h_{[1]}^{\la}[n-2k] x[n],\;\mu, ~\la=0,1,$ that belong to $\Pi[N/2]$, we get their  \oo\ projections $\mathbf{y}_{[1]}^{\la,0}$ and $\mathbf{y}_{[1]}^{\la,1}\in \Pi[N/2]$  onto the subspaces $\Pi^{0}[N/2]$  and $\Pi^{1}[N/2]$:
\begin{eqnarray*}\label{ort_proj2y00}
         y_{[1]}^{\la,\mu}[k]&=&\sum_{l=0}^{N/4-1}y_{[2]}^{\rr}[l]\, {h}_{[2]}^{\mu}[k-2l]  , \quad
         y_{[2]}^{\rr}[l] = \sum_{k=0}^{N/2-1}h_{[2]}^{\mu}[k-2l] \, y_{[1]}^{\la}[k]
\\\nn
     &=& \sum_{k=0}^{N/2-1}{h}_{[2]}^{\mu}[k-2l] \, \sum_{n=0}^{N-1}{h}_{[1]}^{\la}[n-2k] x[n]=   \sum_{n=0}^{N-1}x[n]\,{\psi}_{[2],\rr}^{2r}[n-4l], \\\label{psi20}
   {\psi}_{[2],\rr}^{2r}[n]  &\srr& \sum_{k=0}^{N/2-1}{h}_{[2]}^{\mu}[k] \,{h}_{[1]}^{\la}[n-2k]=\sum_{k=0}^{N/2-1}{h}_{[2]}^{\mu}[k] \, {\psi}_{[1],\la}^{2r}[n-2k].
%\\\label{ort_proj2y01}   y_{[1]}^{0,1}[k]&=&\sum_{l=0}^{N/4-1}y_{[2]}^{1}[l]\, {h}_{[2]}^{1}[k-2l]  , \quad
%         y_{[2]}^{1}= \sum_{k=0}^{N/2-1}h_{[2]}^{1}[k-2l] \, y_{[1]}^{0}[k]
         \end{eqnarray*}
where % $\la,\mu=0,1,\;
$\rr=\left\{
     \begin{array}{ll}
     \mu, & \hbox{if $\la=0$;} \\
    3-\mu, & \hbox{if $\la=1$.}
   \end{array}
  \right.$

 The \ss\  ${\psi}_{[2],\rr}^{2r}$ is a linear combination of 2-sample shifts of the \dss\ WP ${\psi}_{[1],\la}^{2r}$, therefore ${\psi}_{[2],\rr}^{2r}\in{} ^{2r}{\mathcal{S}}_{[1]}^{\la}\subset \Pi[N]$. Its DFT is
 \begin{equation}\label{spec_psi20}
 \hat{ {\psi}}_{[2],\rr}^{2r}[n]=\hat{{\psi}}_{[1],\la}^{2r}[n]\,\hat{h}_{[2]}^{\mu}[n]_{1}.
 \end{equation}

 \bpp[\cite{ANZ_book3}]\label{psi20_pro}  The norms of the \ss s ${\psi}_{[2],\rr}^{2r}\in \Pi[N]$ are equal to one. The 4-sample shifts $\left\{\psi_{[2],\rr}^{2r}[\cdot-4l] \right\},\;l=0,...,N/4-1,$ of this \ss\ are mutually \oo\ and \ss s with \df t indices $\rr$ are \oo\ to each other.\epp
 Thus, the \ss\ space $\Pi[N]$  splits into four mutually \oo\ subspaces $\Pi[N]=\bigoplus_{\rr=0}^{3} {}^{2r}\mathcal{S}_{[1]}^{\rr}$ whose \on\ bases are formed by 4-sample shifts  $\left\{\psi_{[2],\rr}^{2r}[\cdot-4l] \right\},\;l=0,...,N/4-1,$ of the \ss s $\psi_{[2],\rr}^{2r}$, which are referred to as the second-level \dss\ \wq s of order $2r$.

 The \oo\  projection of a \ss\ $\mathbf{x}\in \Pi[N]$ onto the subspace  ${}^{2r}\mathcal{S}_{[2]}^{\rr}$ is the \ss\
 \begin{equation*}\label{s20_rep}
 x_{[2]}^{\rr}[k]= \sum_{l=0}^{N/4-1}\left\langle \mathbf{x}, \,\psi_{[2],\rr}^{2r}[\cdot-4l] \right\rangle\,{\psi}_{[2],\rr}^{2r}[k-4l]=\sum_{l=0}^{N/4-1}y_{[2]}^{\rr}[l]\,{\psi}_{[2],\rr}^{2r}[k-4l], ~~
 k=0, \ldots , N -1.
\end{equation*}

   Practically, derivation of the \wq\ \t\  \c s $\mathbf{y}_{[1]}^{\la},\;\la=0,1,$ from $\mathbf{x}$ and the inverse operation are implemented using \eh{mod_repAS1}, while the \t\
 $\mathbf{y}_{[1]}^{\la}\longleftrightarrow\mathbf{y}_{[2]}^{\rr}$ are implemented similarly using the  \mv ces of the     \pb\
${\mathbf{H}}_{[2]}$ defined in \eh{sa_modma20}.
  The second-level  \wq s $\psi_{[2],\rr}^{2r}$ are derived from the first-level  \wq s $\psi_{[1],\la}^{2r}$ by \fr ing the latter with the \pf s  $\mathbf{h}_{[2]}^{\mu},\;\la,\mu=0,1,\;\rr=\left\{
                                                   \begin{array}{ll}
                                                     \mu, & \hbox{if $\la=0$;} \\
                                                     3-\mu, & \hbox{if $\la=1$.}
                                                   \end{array}
                                                 \right.
$.

\begin{figure}[H]
\begin{center}
\resizebox{14cm}{5cm}{
\includegraphics{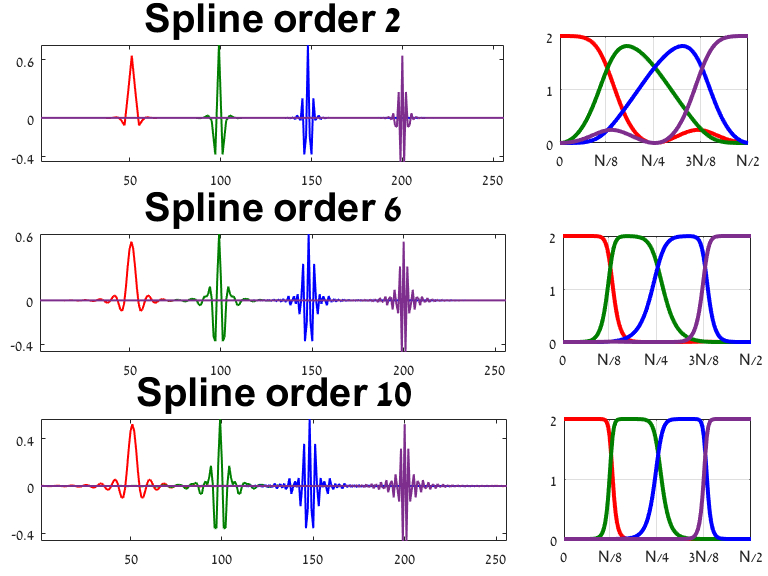}
}
\end{center}
\caption{Left: second-level \dss\ \wq s of different orders; left to right: $\psi_{[2],0}^{2r}\to\psi_{[2],1}^{2r}\to\psi_{[2],2}^{2r}\to\psi_{[2],3}^{2r}$. Right: magnitude DFT \sp a of these \wq s}
     \label{dss_wq_s2}
  \end{figure}

 Figure \ref{dss_wq_s2} displays the second-level \wq s originating from \dd\ \s s  of orders 2, 6 and  10   and their  DFTs. One can observe that the \wq s are \sy ic and well localized  in time domain.
    Their \sp a are flat and their shapes tend to rectangular as their orders increase. They split the \ff\ domain into four quarter-bands.
\subsubsection{Transforms to deeper levels}\label{sec:sss243}
The WPTs to deeper \d\ levels are implemented iteratively,  while the \t\ \c s $\left\{\mathbf{y}_{[m+1]}^{\rr}\right\}$ are derived by \fr ing the \c s $\left\{\mathbf{y}_{[m]}^{\la}\right\}$ with the \pf s $\mathbf{h}^{\mu}_{[m+1]},$ where $\la=0,...,2^m-1,\;\mu=0,1$ and $\rr=\left\{
                                                                                                                               \begin{array}{ll}
                                                                                                                                 2\la +\mu, & \hbox{if $\la$ is even;} \\
                                                                                                                                 2\la +(1-\mu), & \hbox{if $\la$ is odd.}
                                                                                                                               \end{array}
                                                                                                                             \right.
$
The \t\ \c s are ${y}_{[m]}^{\la}[l]=\left\langle \mathbf{x},\psi^{2r}_{[m],\la}[\cdot, -2^{m}l] \right\rangle$, where the \ss s $\psi^{2r}_{[m],\la}$ are normalized, \oo\ to each other in the space $\Pi[N]$, and their $2^{m}l-$sample shifts are mutually \oo. They are referred to as level-$m$ \dss\ \wq s of order $2r$. The set $\left\{\psi^{2r}_{[m],\la}[\cdot, -2^{m}l] \right\},\;\la=0,...,2^m-1,\;l=0,...N/2^m-1,$ constitutes an \on\ basis of the  space $\Pi[N]$ and generates its split into $2^m$ \oo\ subspaces. The next-level \wq s $\psi^{2r}_{[m+1],\rr}$ are derived by \fr ing the \wq s  $\psi^{2r}_{[m],\la}$ with the \pf s $\mathbf{h}^{\mu}_{[m+1]}$ such that
\begin{equation}\label{mlev_wq}
  {\psi}_{[m+1],\rr}^{2r}[n]  =\sum_{k=0}^{N/2^{m}-1}{h}_{[m+1]}^{\mu}[k] \, {\psi}_{[m],\la}^{2r}[n-2^{m}k].
\end{equation}
Note that the \fg\ of an $m-$level \pf\ is $ \hat{h}^{\mu}_{[m]}[n]=\hat{h}^{\mu}_{[1]}[2^{m-1}n].$

A scheme of fast implementation of the \dss-based WPT  is described in \cite{ANZ_book3}.   The \t s are executed in the spectral domain using the  Fast Fourier transform  (FFT) by the application of critically sampled two-channel \fb s to the half-band \sp al  components of a \ss. For example, the Matlab execution of the 8-level 12-th-order WPT of a signal  comprising  245760 samples, takes  0.2324  seconds.

\subsection{ 2D WPTs}\label{sec:ss25}
A standard way to extend the one-dimensional (1D) WPTs to multiple dimensions is the tensor-product extension.
The 2D one-level WPT of a  \ss\ $\mathbf{x}=\left\{x[k,n]\right\},\;k,n=0,...,N-1,$  which  belongs to $\Pi[N,N]$, consists of the application of  1D WPT to columns of  $\mathbf{x}$, which is followed by the application of the \t\ to  rows of the  \c s array. As a result of the 2D WPT of \ss s from $\Pi[N,N]$, the space becomes split
 into four mutually \oo\ subspaces
$ \Pi[N,N]=\bigoplus_{j,l=0}^{1}\,^{2r} {\mathcal{S}}^{j,l}_{[1]}.$

The subspace ${}^{2r} {\mathcal{S}}^{j,l}_{[1]}$ is a linear hull of two-sample shifts of the 2D \wq s\\
$\left\{\psi_{[1],j,l}^{2r}[k-2p,n-2t]\right\} ,\;p,t,=0,...,N/2-1,$ that form an \on\ basis of ${}^{2r} {\mathcal{S}}^{j,l}_{[1]}$. The \oo\ projection of the  \ss\ $\mathbf{x}\in\Pi[N,N]$ onto the subspace ${}^{2r} {\mathcal{S}}^{j,l}_{[1]}$ is the  \ss\ $\mathbf{x}_{[1]}^{j,l}\in\Pi[N,N]$ such that
\begin{equation*}\label{op_2d}
 {x}_{[1]}^{j,l}[k,n]=\sum_{p,t=0}^{N/2-1} y_{[1]}^{j,l}[p,t] \,\psi_{[1],j,l}^{2r}[k-2p,n-2t], \quad j,l =0,1,
\end{equation*}

The 2D \wq s are $\psi_{[1],j,l}^{2r}[n,m]\srr \psi_{[1],j }^{2r}[n]\, \psi_{[1],l}^{2r}[m],   \quad j,l=0,1.$ They are normalized and  \oo\ to each other in the space $\Pi[N,N]$. It means that \\ $\sum_{n,m=0}^{N-1}\psi_{[1],j1 ,l1}^{2r}[n,m]\,\psi_{[1],j2 ,l2}^{2r}[n,m]=\da[j1-j2]\,\da[l1-l2]$. Their two-sample shifts in either direction are mutually \oo. The \t\ \c s are $$y_{[1]}^{j,l}[p,t] =\left\langle \mathbf{x},\psi_{[1],j,l}^{2r}[\cdot-2p,\cdot-2t] \right\rangle=\sum_{n,m=0}^{N-1} \psi_{[1],j,l}^{2r}[n-2p,m-2t]\: x[n,m].$$

By the application of the above \t s iteratively to blocks of the \t\ \c s down to $m$-th level, we get that the space $ \Pi[N,N]$ is decomposed into $4^{m}$ mutually \oo\ subspaces
$ \Pi[N,N]=\bigoplus_{j,l=0}^{2^{m}-1}\,^{2r} {\mathcal{S}}^{j,l}_{[m]}.$  The \oo\ projection of the  \ss\ $\mathbf{x}\in\Pi[N,N]$ onto the subspace ${}^{2r} {\mathcal{S}}^{j,l}_{[m]}$ is the  \ss\ $\mathbf{x}_{[m]}^{j,l}\in\Pi[N,N]$ such that
\begin{eqnarray*}\label{op_2d}
 {x}_{[m]}^{j,l}[k,l]&=&\sum_{p,t=0}^{N/2^{m}-1} y_{[m]}^{j,l}[p,t] \,\psi_{[m],j ,l}^{2r}[k-2^{m}p,l-2^{m}t], \quad j,l =0,...,2^{m}-1,\\
\psi_{[m],j ,l}^{2r}[k,n]&=&\psi_{[m],j}^{2r}[k]\,\psi_{[m],l}^{2r}[n],\quad  y_{[m]}^{j,l}[p,t]
=\left\langle \mathbf{x},\psi_{[m],j ,l}^{2r}[\cdot-2^{m}p,\cdot-2^{m}t] \right\rangle.
\end{eqnarray*}

The 2D tensor-product  \wq s $\psi_{[m],j ,l}^{2r}$ are well localized in the spatial domain, their 2D DFT \sp a are flat and provide a refined split of the \ff\ domain of \ss s from $ \Pi[N,N].$\footnote{Especially it is true for WPs derived from  higher-order \ds s.} The drawback is that the \wq s are  oriented  in ether horizontal or vertical directions or are not oriented  at all.

\begin{figure}[H]
\centering
\includegraphics[width=3.2in]{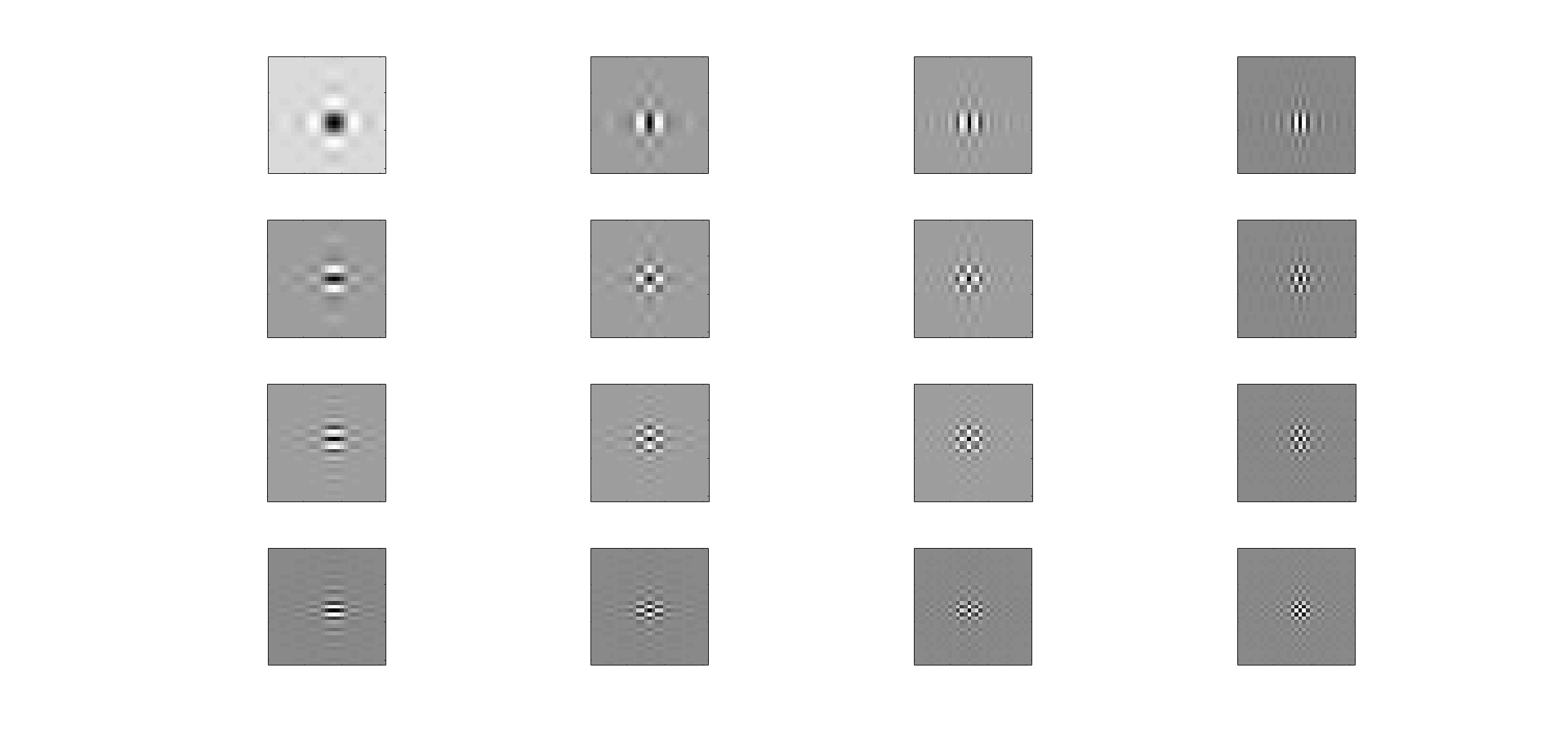}
\hfil
\includegraphics[width=3.2in]{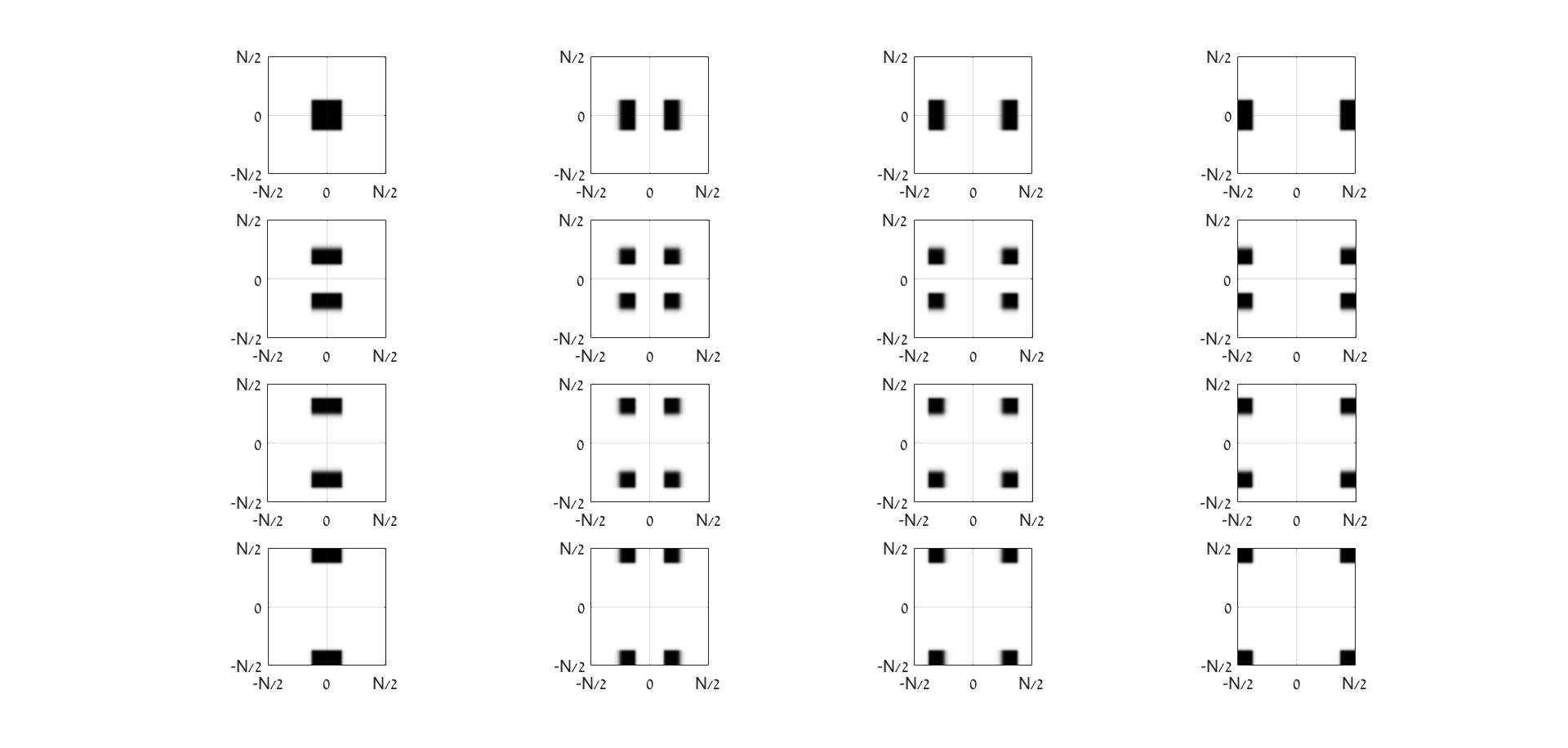}%
\caption{WPs from the second \d\ level (left) and their magnitude \sp a (right)}
\label{psifpsi2_2}
\end{figure}

Figure \ref{psifpsi2_2} displays the tenth-order 2D \wq s  from the second \d\ level and their magnitude \sp a.

\section{(Quasi-)analytic and complementary WPs}\label{sec:s3}
In this section we define \az\ and the so-called quasi-\az\  WPs related to the \dss-based WPs discussed in Section \ref{sec:s2} and introduce an \on\ set of waveforms which are complementary to the above WPs.
\subsection{ Analytic \p\ \ss s}\label{sec:ss31}

A \ss\ $\mathbf{x}\in\Pi[N]$ is \ry ed by its  inverse DFT. Then, \eh{fsx} can be written as follows:
\begin{eqnarray*}\label{fsx+}
  x[k]&=&\frac{\hat{x}[0]+(-1)^{k}\hat{x}[N/2]}{N}+\frac{2}{N}\sum_{n=1}^{N/2-1}\frac{\hat{x}[n]\,\w^{kn}+(\hat{x}[n]\,\w^{kn})^{\ast}}{2}.
\end{eqnarray*}
Define the real-valued  \ss\ $\mathbf{h}\in\Pi[N]$ and two complex-valued \ss s $\mathbf{\bar{x}}_{+}$ and  $\mathbf{\bar{x}}_{-}$ such that
\begin{equation}
\label{yy}
\begin{array}{lll}
 h[k]&\srr&\frac{2}{N}\sum_{n=1}^{N/2-1}\frac{\hat{x}[n]\,\w^{kn}-\hat{x}[n]^{\ast}\,\w^{-kn}}{2i},\\
 \bar{x}_{\pm}[k]&\srr&x[k]\pm ih[k]=\frac{\hat{x}[0]+(-1)^{k}\hat{x}[N/2]}{N}\\&+&\frac{2}{N}\sum_{n=1}^{N/2-1}
 \left\{
   \begin{array}{ll}
     \hat{x}[ n]\,\w^{ kn}, & \hbox{for $\bar{x}_{+}$;} \\
      \hat{x}[ n]^{\ast}\,\w^{- kn}=\hat{x}[ N-n]\,\w^{- k(N-n)}, & \hbox{for $\bar{x}_{-}$.}
   \end{array}
 \right.
\end{array}
\end{equation}

The \ss s' $\mathbf{\bar{x}}_{\pm}$ DFT \sp a are
\begin{equation}
\label{xpm_sp}
\begin{array}{lll}
  \hat{\bar{x}}_{+}[n] &=& \left\{
               \begin{array}{ll}
                 \hat{x}[n], & \hbox{if $n=0$, or  $n=N/2$ ;} \\
                 2\hat{x}[n], & \hbox{if $0<n<N/2$;} \\
                 0, & \hbox{if $-N/2<n<0\Longleftrightarrow N/2<n<N$,}
               \end{array}
             \right. \\
  \hat{\bar{x}}_{-}[n] &=& \left\{
               \begin{array}{ll}
                 \hat{x}[n], & \hbox{if $n=0$, or  $n=N/2$ ;} \\
                 2\hat{x}[n], & \hbox{if $-N/2<n<0\Longleftrightarrow N/2<n<N$;} \\
                 0, & \hbox{if $n=0$, or  $n=N/2$.}
               \end{array}
             \right.
\end{array}
\end{equation}
 %The \ss\ $\mathbf{y}$ is a \dd\ \p\ version of the Hilbert \t\ (HT) of the \ss\ $\mathbf{x}\in\Pi[N]$.
The \sp um of $\mathbf{\bar{x}}_{+}$ comprises only non-negative frequencies and vice versa for $\mathbf{\bar{x}}_{-}$.
$\mathbf{x}=\mathfrak{Re}(\mathbf{\bar{x}}_{\pm})$ and $\mathfrak{Im}(\mathbf{\bar{x}}\pm)=\pm\mathbf{h}$. The  \ss s  $\mathbf{\bar{x}}_{\pm}$ are referred to as  \p\ analytic \ss s.

The \ss's $\mathbf{h}$ DFT \sp um is
\begin{equation*}\label{yft}
\hat{h}[n]=\left\{
               \begin{array}{ll}
                  -i\,\hat{x}[n], & \hbox{if $0<n<N/2$;} \\
                i\, \hat{x}[n], & \hbox{if $-N/2<n<0\Longleftrightarrow N/2<n<N$;} \\
                 0, & \hbox{if $n=0$, or  $n=N/2$.}
               \end{array}
             \right.
\end{equation*}
Thus, the \ss\ $\mathbf{h}$ where  $\mathbf{h}=H(\mathbf{x})$ can be regarded as the  Hilbert \t\ (HT) of a \dt\ \p\ \ss\  $\mathbf{x}$,  (see \cite{opp}, for example).
\bpp\label{pro:ysym}
\par\noindent
\begin{enumerate}
  \item The HT  $\mathbf{h}=H(\mathbf{x})$  is invariant with respect to circular shift in $\Pi[N]$. That means that $\mathbf{\tilde{h}}=\mathbf{h}[\cdot +m]$ is the HT of  $\mathbf{\tilde{x}}=\mathbf{x}[\cdot +m]$.
      \item  If   the \ss\ $\mathbf{x}\in\Pi[N]$  is \sy ic about a grid point $k=K$ than   $\mathbf{h}=H(\mathbf{x})$ is anti\sy ic about \emph{K}  and $h[K]=0$.
 \item Assume that a   \ss\ $\mathbf{x}\in\Pi[N]$ and $\hat{x}[0]=\hat{x}[N/2]=0$. Then,
\begin{enumerate}
  \item The norm of  its HT  is $\|H(\mathbf{x}) \|=\|\mathbf{x} \|$.
  \item Magnitude \sp a of the \ss s $\mathbf{x}$ and $\mathbf{h}=H(\mathbf{x})$ coincide.
  \end{enumerate}
\end{enumerate}\epp
\proof
\begin{enumerate}
 \item The DFT of the signal $\mathbf{\tilde{x}}$ is  $\hat{\tilde{x}}[n]=\w^{mn}\,\hat{x}[n]$. Denote by $\mathbf{\bar{\tilde{x}}}_{+}$ the \az\ \ss\ related to $\mathbf{\tilde{x}}$.  \ehh{xpm_sp} implies that $\hat{\bar{\tilde{x}}}_{+}[n]=\w^{mn}\,\hat{\bar{x}}[n]$. Consequently, $\bar{\tilde{x}}_{+}[k]=\bar{x}_{+}[k+m]$. The same relation holds for $\mathbf{\tilde{h}}=\mathfrak{Im}(\mathbf{\bar{\tilde{x}}})$.
     \item Assume that $\mathbf{x}\in\Pi[N]$  is \sy ic about  $K=0$. Then, its DFT is $$\hat{x}[n]=x[0]+(-1)^{n}x[N/2] +2\sum_{k=1}^{N/2-1}x[k] \,\cos(2\pi kn/N)=\hat{x}[-n].$$
Then, due to \eh{yy}, $h[k]=2/N\sum_{n=1}^{N/2-1}x[n]\sin(2\pi kn/N)=h[-k]$ and $h[0]=0$. Extension of the proof to $K\neq0$ is straightforward.
\item
\begin{enumerate}
\item The squared norm is
\(
  \|\mathbf{h} \|^{2}=\frac{1}{N}\sum_{n=-N/2}^{N/2-1}|\hat{h}[n]|^{2}=\frac{1}{N}\sum_{n=-N/2}^{N/2-1}|\hat{x}[n]|^{2}.
\)%end{equation*}
\item The coincidence of  the magnitude \sp a is straightforwared.

\end{enumerate}\end{enumerate}\eop

\subsection{Analytic  WPs}\label{sec:ss32}
The \az\ \s-based WPs and their DFT \sp a are  derived from the corresponding WPs $\left\{\psi^{2r}_{[m],l}\right\},\;m=1,...,M,\;l=0,...,2^{m}-1,$ in line with the scheme in Section \ref{sec:ss31}. Recall that for all $l\neq0$, the DFT $\hat{\psi}^{2r}_{[m],l}[0]=0$ and  for all $l\neq2^{m}-1$, the DFT $\hat{\psi}^{2r}_{[m],l}[N/2]=0$.

Denote by $\th^{2r}_{[m],l}=H(\psi^{2r}_{[m],l})$ the \dd\ \p\ HT  of the \wq\ $\psi^{2r}_{[m],l}$, such that the DFT is
\begin{equation*}\label{th_df}
  \hat{\th}^{2r}_{[m],l}[n]=\left\{
               \begin{array}{ll}
                  -i\,\hat{\psi}^{2r}_{[m],l}[n], & \hbox{if $0<n<N/2$;} \\
                i\, \hat{\psi}^{2r}_{[m],l}, & \hbox{if $-N/2<n<0$;} \\
                 0, & \hbox{if $n=0$, or  $n=N/2$ .}
              \end{array}
             \right.
\end{equation*}
 Then, the corresponding \az\ WPs are
\begin{equation*}\label{awq}
\bar{\psi}^{2r}_{\pm[m],l}=\psi^{2r}_{[m],l}   \pm i\th^{2r}_{[m],l}.
\end{equation*}
\paragraph{Properties of the \az\ WPs}
\begin{enumerate}
  \item The DFT \sp a of the \az\ WPs $\bar{\psi}^{2r}_{+[m],l}$ and $\bar{\psi}^{2r}_{-[m],l}$ are located within the bands $[0,N/2]$ and $[N/2,N]\Longleftrightarrow[-N/2,0]$, respectively.
  \item The real component ${\psi}^{2r}_{[m],l}$ is the same for both WPs $\bar{\psi}^{2r}_{\pm[m],l}$. {It} is a \sy ic oscillating waveform.
  \item\label{prop3}
The HT WPs $\th^{2r}_{[m],l}=H({\psi}^{2r}_{[m],l})$ are anti\sy ic oscillating waveforms.
\item For all $l\neq0, \,2^{m}-1$, the  norms  $\left\| \th^{2r}_{[m],l}\right\|=1$.  Their magnitude \sp a $\left|\hat{\th}^{2r}_{[m],l}[n]\right|$  coincide with the magnitude \sp a of the respective WPs  $\psi^{2r}_{[m],l}$.
    \item\label{prop5}  When $l=0$ or  $l=2^{m}-1$,  the magnitude \sp a of $\th^{2r}_{[m],l}$ coincide with that of ${\psi}^{2r}_{[m],l}$ everywhere except for the points $n=0$ or $N/2,$ respectively,  and the waveforms'  norms are no longer equal to 1.
\end{enumerate}
Properties in items \ref{prop3}--\ref{prop5} follow directly from Proposition \ref{pro:ysym}.

\begin{figure}[H]
\resizebox{16cm}{10cm}{
\centering
\includegraphics{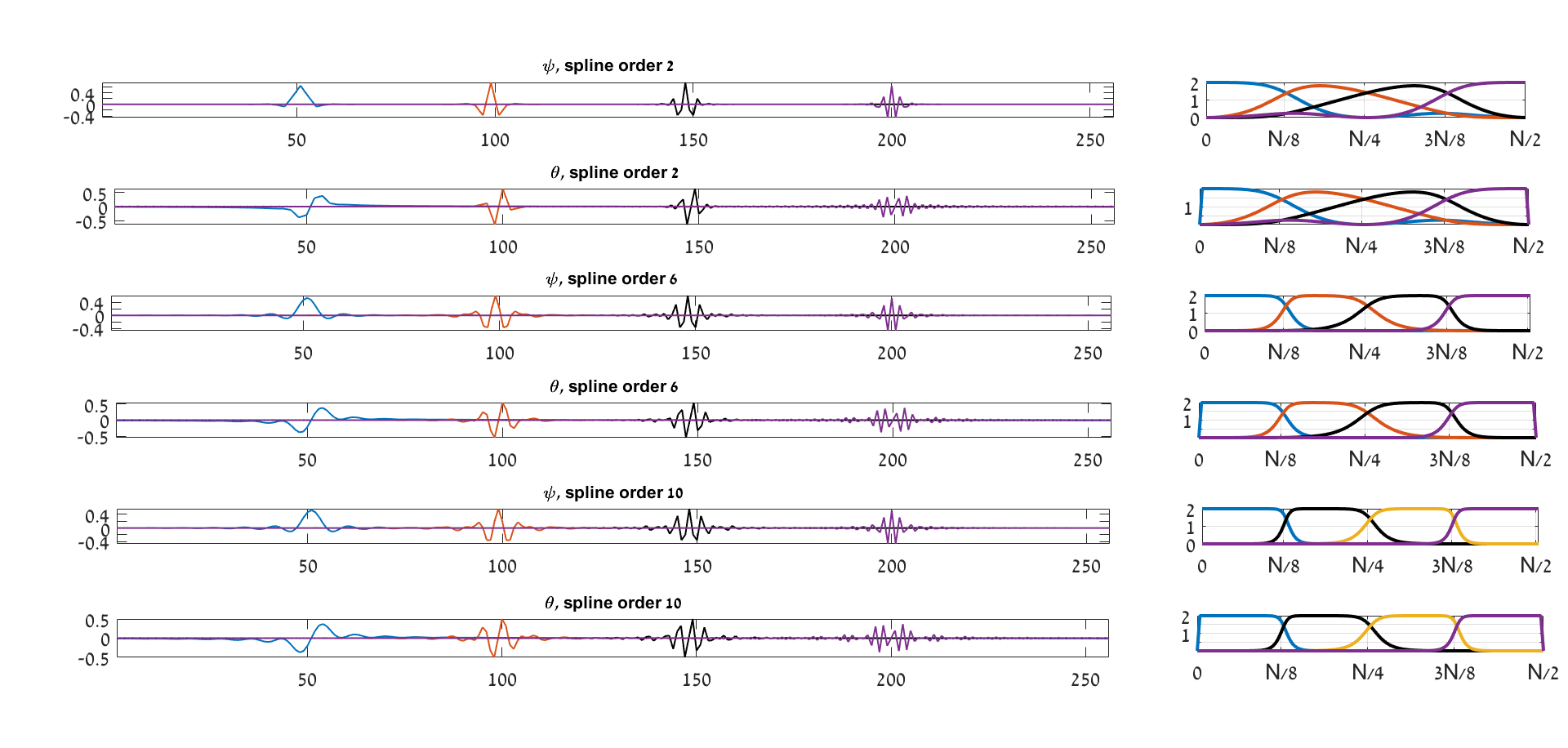}
}
\caption{WPs ${\psi}^{2r}_{[2],l}$ (first, third and fifth  from the top left frames) and ${\th}^{2r}_{[2],l}$ (second, fourth and sixth  from the top left frames) from the second \d\ level and their magnitude \sp a, respectively (right frames)}
\label{psi_theta2}
\end{figure}

\bpp\label{pro:teta_oo} { For all $l\neq0, \,2^{m}-1$,
the shifts of  the HT WPs $\left\{\th^{2r}_{[m],l}[\cdot-2^{m}l]\right\}$ are \oo\ to each other in the space $\Pi[N]$. The \oo ity does not take place for  for $\th^{2r}_{[m],0}$ and $\th^{2r}_{[m],2^{m}-1}$.}
\epp
\proof The inner product is
\begin{eqnarray*}
   && \left\langle \th^{2r}_{[m],l},\th^{2r}_{[m],l}[\cdot-2^{m}l]  \right\rangle =\frac{1}{N}\sum_{n=-N/2}^{N/2-1}\w^{2^{m}ln}\left|\hat{\th}^{2r}_{[m],l}[n] \right|^{2}\\
&& =\frac{1}{N}\sum_{n=-N/2}^{N/2-1}\w^{2^{m}ln}\left|\hat{\psi}^{2r}_{[m],l}[n] \right|^{2}=\left\langle \psi^{2r}_{[m],l},\psi^{2r}_{[m],l}[\cdot-2^{m}l]  \right\rangle=0.
\end{eqnarray*}
\eop
Figure \ref{psi_theta2} displays the \wq s ${\psi}^{2r}_{[2],l}$ and ${\th}^{2r}_{[2],l},\;r=1,3,5,\; l=0,1,2,3,$ and their magnitude  \sp a.

\subsection{Complementary set of wavelet packets and quasi-\az\ WPs}\label{sec:ss33}
\subsubsection{Complementary \on\  WPs}\label{sec:sss331}
The \dss-based WPs $\left\{\psi^{2r}_{[m],l}\right\}$ are normalized and their $2^{m}$-sample shifts are mutually \oo. Combinations of shifts of several \wq s can form \on\ bases for the \ss\ space $\Pi[N]$. On the other hand, it is not true for the set $\left\{\th^{2r}_{[m],l}\right\},\;l=0,...2^{m}-1,$ of the anty\sy ic \we s, which are the HTs of the WPs  $\left\{\psi^{2r}_{[m],l}\right\}$.
  At  the \d\ level \emph{m},  the \we s $\left\{\th^{2r}_{[m],l}\right\},\;l=1,...2^{m}-2,$ are normalized and their $2^{m}$-sample shifts are mutually \oo, but the norms of the
\we s $\th^{2r}_{[m],0}$ and $\th^{2r}_{[m],2^{m}-1}$ are close  but not equal to 1 and their shifts are not mutually \oo. It happens because the values $\hat{\th}^{2r}_{[m],j}[0]$ and $\hat{\th}^{2r}_{[m],j}[N/2]$ are missing in their DFT \sp a\footnote{Recall that these values are real}. This keeping in mind, we upgrade the set $\left\{\th^{2r}_{[m],l}\right\},\;l=0,...2^{m}-1$ in the following way.

Define
a set $\left\{\f^{2r}_{[m],l}\right\},\;m=1,...,M, \;l=0,...,2^{m}-1,$  of \ss s from the space $\Pi[N]$ via their DFTs:
\begin{equation}\label{phi_df}
  \hat{\f}^{2r}_{[m],l}[n]=\left\{
               \begin{array}{ll}
                  -i\,\hat{\psi}^{2r}_{[m],l}[n], & \hbox{if $0<n<N/2$;} \\
                i\, \hat{\psi}^{2r}_{[m],l}[n], & \hbox{if $-N/2<n<0$ ;}\\
                 \hat{\psi}^{2r}_{[m],l}[n], & \hbox{if $n=0$, or  $n=N/2 .$}
               \end{array}
             \right.
\end{equation}
For all $l\neq0, 2^{m}-1,$ the \ss s $\f^{2r}_{[m],l}$ coincide with $\th^{2r}_{[m],l}=H(\psi^{2r}_{[m],l})$.
\bpp\label{pro:phi_oo}
\par\noindent
\begin{description}
  \item[-] The magnitude \sp a  $\left|\hat{\f}^{2r}_{[m],l}[n]\right|$  coincide with the magnitude \sp a of the respective WPs  $\psi^{2r}_{[m],l}$.
  \item[-] For any   $m=1,...,M,$ and $l=1,...,2^{m}-2,$ the \ss s \ $\f^{2r}_{[m],l}$ are anti\sy ic oscillating waveforms.   For $l=0$ and $l=2^{m}-1$, the shapes of the \ss s are near anti\sy ic.
      \item[-] The \on ity properties that are similar to the properties of WPs  $\psi^{2r}_{[m],l}$ hold for the \ss s  $\f^{2r}_{[m],l}$ such that
       \begin{eqnarray*}
      \label{onp}
        \left\langle\f^{2r}_{[m],l}[\cdot -p\,2^{m}],\f^{2r}_{[m],\la}[\cdot -s\,2^{m}] \right\rangle= \da[\la,l]\,\da[p,s].
      \end{eqnarray*}
\end{description}
\epp
The proof of Proposition \ref{pro:phi_oo} is similar to the proof of Proposition \ref{pro:teta_oo}. %Note that due to the addition of values $\hat{\psi}^{2r}_{[m],l}[0]$ and $\hat{\psi}^{2r}_{[m],l}[N/2]$ into the DFT of \ss s $\f^{2r}_{[m],l}$, the \ss s $\f^{2r}_{[m],0}$ and $\f^{2r}_{[m],2^{m}-1}$ are no longer anti\sy ic.

\begin{figure}[H]
\resizebox{17cm}{4cm}{
\centering
\includegraphics{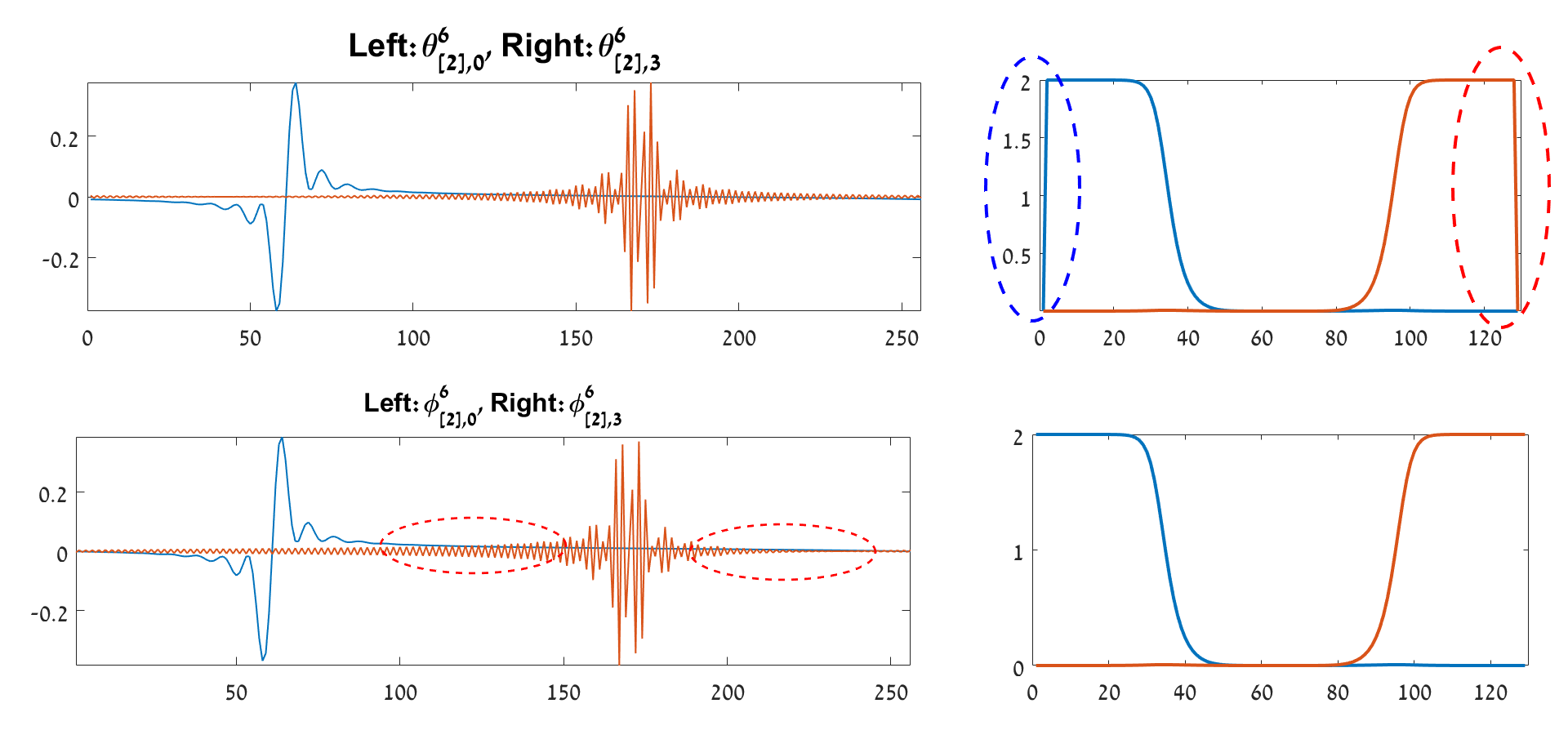}
}
\caption{Left: \ss s ${\th}^{6}_{[2],l},\;l=0,3$ (top), and ${\f}^{6}_{[2],l},\;l=0,3$  (bottom). Right:  their magnitude DFT \sp a, respectively}
     \label{teta_phi}
\end{figure}

  Figure \ref{teta_phi} displays  the \ss s ${\th}^{6}_{[2],l},\;l=0,3$ and ${\f}^{6}_{[2],l},\;l=0,3$, from the second \d\ level and their magnitude \sp a. Lack of the values $\hat{\th}^{2r}_{[m],j}[0]$ and $\hat{\th}^{2r}_{[m],j}[N/2]$ in the DFTs of ${\th}^{6}_{[2],l},\;l=0,3$, are seen. Addition of $\hat{\psi}^{2r}_{[m],j}[0]$ and $\hat{\psi}^{2r}_{[m],j}[N/2]$ to the above \sp a results in the anti\sy y distortion.
%\paragraph{Properties}

We call the \ss s     $\left\{\f^{2r}_{[m],l}\right\},\;m=1,...,M, \;l=0,...,2^{m}-1$, the \emph{complementary  wavelet packets} (cWPs). Similarly to the WPs $\left\{\psi^{2r}_{[m],l}\right\},$ \df ent combinations of the cWPs can provide \df ent  \on\ bases for the space $\Pi[N]$.  These  can be, for example, the \ww\ bases $$\left\{\bigoplus_{r=0}^{N/2^{M}}\f^{2r}_{[M],0}[\cdot-r\,2^{M}]\right\}\bigoplus_{m=1}^{M}
\left\{\bigoplus_{r=0}^{N/2^{m}}\f^{2r}_{[m],1}[\cdot-r\,2^{m}]\right\}.$$
or a  Best Basis \cite{coiw1} type.
\subsubsection{Quasi-\az\  WPs}\label{sec:sss332}
The sets of complex-valued WPs, which we refer to as the quasi-\az\  wavelet packets (qWP),  are defined as
\begin{equation*}\label{qaz}
  \Psi^{2r}_{\pm[m],l}=\psi^{2r}_{[m],l}   \pm i\f^{2r}_{[m],l}, \quad m=1,...,M,\;l=0,...,2^{m}-1,
\end{equation*}
where $\f^{2r}_{[m],l}$ are the cWPs defined in \eh{phi_df}. The qWPs $\Psi^{2r}_{\pm[m],l}$ differ from the \az\ WPs $\bar{\psi}^{2r}_{\pm[m],l}$ by the addition of the two values $\pm i\,\hat{\psi}^{2r}_{[m],l}[0]$ and $\pm i\,\hat{\psi}^{2r}_{[m],l}[N/2]$ into their DFT \sp a, respectively. For a given \d\ level $m$, these values are zero for all $l$ except for $l_{0}=0$ and $l_{m}=2^{m}-1$. It means that for all $l$ except for $l_{0}$ and $l_{m}$, the qWPs $\Psi^{2r}_{\pm[m],l}$ are \az.
The DFTs of qWPs are
\begin{eqnarray}\label{qa_df}
  \hat{\Psi}^{2r}_{+[m],l}[n]&=&\left\{
               \begin{array}{ll}
  (1+i)\hat{\psi}^{2r}_{[m],l}[n], & \hbox{if $n=0$ or $n= N/2$;} \\
                  2\hat{\psi}^{2r}_{[m],l}[n], & \hbox{if $0< n< N/2$;} \\
               0 & \hbox{if $ N/2<n<N$,}
               \end{array}
             \right.\\\nn
 \hat{\Psi}^{2r}_{-[m],l}[n]&=&\left\{
               \begin{array}{ll}
(1-i)\hat{\psi}^{2r}_{[m],l}[n], & \hbox{if $n=0$ or $n= N/2$;} \\
                  0 & \hbox{if $0< n< N/2$;} \\
               2\hat{\psi}^{2r}_{[m],l}[n], & \hbox{if $ N/2< n< N$.}
               \end{array}
             \right.
\end{eqnarray}
\subsubsection{Design of cWPs and qWPs}\label{sec:sss333}
The DFTs of the first-level WPs are
\begin{equation*}\label{df_wq1}
  \hat{\psi}^{2r}_{[1],0}[n]=\frac{\cos^{2r}\frac{\pi\,n}{N}}{\sqrt{ U^{4r}[n]}}=\b[n],\quad  \hat{\psi}^{2r}_{[1],1}[n]=\frac{\w^{n}\,\sin^{2r}\frac{\pi\,n}{N}}{\sqrt{ U^{4r}[n]}}=\a[n],
\end{equation*}
where the \sq\ $ U^{4r}[n]$ is defined in \eh{u2r}. Consequently, the DFTs of  the  first-level cWPs are

\begin{eqnarray}\label{df_cwq1}
  \hat{\f}^{2r}_{[1],0}[n]=\left\{
               \begin{array}{ll}
                  -i\,\b[n], & \hbox{if $0<n<N/2$;} \\
                i\, \b[n], & \hbox{if $N/2<n<N$;} \\
                 \sqrt{2}, & \hbox{if $n=0$; }
             \\
                 0, & \hbox{ if $n=N/2,$}
               \end{array}
\right.\quad \hat{\f}^{2r}_{[1],1}[n]=\left\{
               \begin{array}{ll}
                  -i\,\a[n], & \hbox{if $0<n<N/2$;} \\
                i\, \a[n], & \hbox{if $N/2<n<N$;} \\
                 0, & \hbox{if $n=0$; }
             \\
                 -\sqrt{2}, & \hbox{ if $n=N/2.$}
               \end{array}
\right.
\end{eqnarray}
Due to \eh{spec_psi20}, the  DFT of the second-level WPs are
\begin{eqnarray}\nn
 \hat{ {\psi}}_{[2],\rr}^{2r}[n]&=&\hat{{\psi}}_{[1],\la}^{2r}[n]\,\hat{h}_{[2]}^{\mu}[n]_{1},\quad \la,\mu=0,1,\;\rr=2\la+\left\{
                                                   \begin{array}{ll}
                                                     \mu, & \hbox{if $\la=0$;} \\
                                                     1-\mu, & \hbox{if $\la=1$.}
                                                   \end{array}
                                                 \right.,\\\label{df_wq2}
\hat{h}_{[2]}^{0}[n]_{1}&=&   \b[2n],\quad \hat{h}_{[2]}^{1}[n]_{1}=\a[2n].
 \end{eqnarray}
For example, assume $\la=\mu=0$ then we have
\[\hat{ \psi}_{[2],0}^{2r}[n]=\hat{{\psi}}_{[1],0}^{2r}[n]\,\hat{h}_{[2]}^{0}[n]_{1}=\frac{\cos^{2r}\frac{\pi\,n}{N}}{\sqrt{ U^{4r}[n]}}\,\frac{\cos^{2r}\frac{2\pi\,n}{N}}{\sqrt{U^{4r}[2n]}}.\]
Keeping in mind that the \sq\ $\b[2n]=\cos^{2r}(2\pi\,n/{N})/\sqrt{U^{4r}[2n]}$  is $N/2-$\p,  we have that the DFT of the corresponding cWP is
\[\hat{\f}^{2r}_{[2],0}[n]=\widehat{H({ \psi}_{[2],0}^{2r})}[n]=\b[2n]\,\left\{
               \begin{array}{ll}
                  -i\,\b[n], & \hbox{if $0<n<N/2$;} \\
                i\, \b[n], & \hbox{if $N/2<n<N$;} \\
                 2, & \hbox{if $n=0$; }
             \\
                 0, & \hbox{ if $n=N/2,$}
               \end{array}
\right.=\hat{{\f}}_{[1],0}^{2r}[n]\,\hat{h}_{[2]}^{0}[n]_{1}=\hat{{\f}}_{[1],0}^{2r}[n]\,\hat{h}_{[1]}^{0}[2n]_{1}.
\]
Similar relations hold for all the second-level cWPs and a general statement is true.
\bpp\label{pro:cwq_des} Assume that for a WP $\psi_{[m+1],\rr}^{2r}$ the relation in \eh{mlev_wq} holds. Then, for the cWP $\f_{[m+1],\rr}^{2r}$ we have
\begin{eqnarray*}\label{mlev_cwq}
  {\f}_{[m+1],\rr}^{2r}[n]  &=&\sum_{k=0}^{N/2^{m}-1}{h}_{[m+1]}^{\mu}[k] \, {\f}_{[m],\la}^{2r}[n-2^{m}k]\Longleftrightarrow\hat{\f}_{[m+1],\rr}^{2r}[\n]=
 \hat{h}_{[1]}^{\mu}[2^{m}\n]_{m}\,\hat{\f}_{[m],\la}^{2r}[\n],\\\nn
\hat{h}_{[1]}^{0}[\n] &=&\hat{\psi}^{2r}_{[1],0}[\n]=\b[\n],\quad  \hat{h}_{[1]}^{1}[\n] = \hat{\psi}^{2r}_{[1],1}[\n]=\a[\n].
\end{eqnarray*}
\epp
\bcc\label{cor:cwq_des} Assume that for a WP $\psi_{[m+1],\rr}^{2r}$, the relation in \eh{mlev_wq} holds. Then, for the qWP $\Psi_{\pm[m+1],\rr}^{2r}$ we have
\begin{eqnarray}\label{mlev_qwq}
  \Psi_{\pm[m+1],\rr}^{2r}[n]  =\sum_{k=0}^{N/2^{m}-1}{h}_{[m+1]}^{\mu}[k] \, \Psi_{\pm[m],\la}^{2r}[n-2^{m}k]\Longleftrightarrow\hat{\Psi}_{\pm[m+1],\rr}^{2r}[\n]=
 \hat{h}_{[1]}^{\mu}[2^{m}\n]_{m}\,\hat{\Psi}_{\pm[m],\la}^{2r}[\n].
\end{eqnarray}
\ecc
\section{Implementation of cWP and qWP \t s}\label{sec:s4}
Implementation of \t s with WPs ${\psi}_{[m],\la}^{2r}$ was discussed in Section \ref{sec:s2}. In this section, we extent that the \t\ scheme to the \t s with cWPs ${\f}_{[m],\la}^{2r}$ and qWPs $\Psi_{[m],\la}^{2r}$.
\subsection{One-level \t s}\label{sec:ss41}
Denote by ${} ^{2r}{\mathcal{C}}_{[1]}^{0}$ the subspace of the \ss\ space $\Pi[N]$, which is the linear hull of the set $\mathbf{W}_{[1]}^{0}=\left\{{\f}_{[1],0}^{2r}[\cdot-2k]\right\},\;k=0,...,N/2-1$. The \ss s from the set $\mathbf{W}_{[1]}^{0}$ form an \on\ basis of the subspace  ${} ^{2r}{\mathcal{C}}_{[1]}^{0}$. Denote by ${} ^{2r}{\mathcal{C}}_{[1]}^{1}$ the \oo\ complement of the  subspace ${} ^{2r}{\mathcal{C}}_{[1]}^{0}$ in the  space $\Pi[N]$. The \ss s from   the set $\mathbf{W}_{[1]}^{1}=\left\{{\f}_{[1],1}^{2r}[\cdot-2k]\right\},\;k=0,...,N/2-1$ form an \on\ basis of the subspace  ${} ^{2r}{\mathcal{C}}_{[1]}^{1}$.
\bpp\label{pro:phi1}
The \oo\ projections of a \ss\  $\mathbf{x}\in \Pi[N]$  onto the spaces   ${} ^{2r}{\mathcal{C}}_{[1]}^{\mu},\;\mu=0,1$ are the \ss s $\mathbf{x}_{[1]}^{\mu}\in\Pi[N]$ such that
\begin{eqnarray*}\label{c1_del1}
          x_{[1]}^{\mu}[k]&=&\sum_{l=0}^{N/2-1}c_{[1]}^{\mu}[l]\, \f_{[1],\mu}^{2r}[k-2l] ,\quad
          %\\label{c1_del}
           c_{[1]}^{\mu}[l]  =\left\langle \mathbf{x},\, \f_{[1],\mu}^{2r}[\cdot-2l]   \right\rangle
       = \sum_{k=0}^{N-1}g_{[1]}^{\mu}[k-2l] \,x[k], \\\nn  g_{[1]}^{\mu}[k] &=&\f_{[1],\mu}^{2r}[k],\quad
         \hat{{g}}_{[1]}^{\mu}[n] = \hat{ \f}_{[1],\mu}^{2r}[n],\quad \mu=0,1.
         \end{eqnarray*}
The DFTs $ \hat{ \f}_{[1],\mu}^{2r}[n]$ of the first-level cWPs are given in \eh{df_cwq1}.
\epp
The  \t s $\mathbf{x}\rightarrow \mathbf{c}_{[1]}^{0}\bigcup\mathbf{c}_{[1]}^{1}$ and back are implemented using the \aa\ $ \tilde{\mathbf{M}}^{c}[n] $ and the \sa\
$ \mathbf{M}^{c}[n] $ \mv ces:

 \begin{equation}
 \label{aa_modma11}
   \begin{array}{lll}
 \tilde{\mathbf{M}}^{c}[n]&\srr& \left(
                         \begin{array}{cc}
                            \hat{g}_{[1]}^{0}[n] &   \hat{g}_{[1]}^{0}\left[n+\frac{N}{2}\right]\\
                            \ \hat{g}^{1}_{[1]}[n] &   \hat{g}^{1}_{[1]}\left[n+\frac{N}{2}\right]\\
                         \end{array}
                       \right)=\left(\begin{array}{cc}
                                   \check{\b}[n] & -\check{\b}\left[n+\frac{N}{2}\right]  \\
                                  \check{\a}[n] &- \check{\a}\left[n+\frac{N}{2}\right]  \\
                                 \end{array}
                               \right),\\
    {\mathbf{M}}^{c}[n]&\srr& \left(\begin{array}{cc}
                                   \check{\b}[n] &            \check{\a}[n] \\
                 -\check{\b}\left[n+\frac{N}{2}\right]     &- \check{\a}\left[n+\frac{N}{2}\right]  \\
                                 \end{array}
                               \right),
\\  % \label{alphabetC}
    \check{\b}[n]      &=&      \left\{
                                                                                                   \begin{array}{ll}
                                                                                                     {\b}[0], & \hbox{if $n=0$;} \\
                                                                                                     -i{\b}[n], & \hbox{otherwise,}
                                                                                                   \end{array}
                                                                                                 \right.
     \quad
       \check{\a}[n]= \left\{
                                                                                                   \begin{array}{ll}
                                                                                                     {\a}[N/2], & \hbox{if $n=N/2$;} \\
                                                                                                     -i\a[n], & \hbox{otherwise.}
                                                                                                   \end{array}
                                                                                                 \right.
  \end{array}
  \end{equation}
The \sq s $\b[n]$ and $\a[n]$ are given in \eh{aa_modma10}.

Similarly to \eh{mod_repAS1}, the one-level cWP \t\ of a \ss\ $\mathbf{x}$ and its inverse are:
   \begin{equation*}\label{mod_repAC1}
    \left(
     \begin{array}{c}
       \hat{c}_{[1]}^{0}[n]_{1} \\
         \hat{c}_{[1]}^{1}[n]_{1}\\
     \end{array}
   \right)=\frac{1}{2} \tilde{\mathbf{M}}^{c}[-n]\cdot \left(
     \begin{array}{l}
      \hat{x}[n] \\
       \hat{x}[\vec{n}]
     \end{array}
   \right),\quad  \left(
     \begin{array}{l}
      \hat{x}[n] \\
       \hat{x}[\vec{n}]
     \end{array}
   \right)={\mathbf{M}}^{c}[n]\cdot \left(
     \begin{array}{c}
        \hat{c}_{[1]}^{0}[n]_{1} \\
         \hat{c}_{[1]}^{1}[n]_{1}\\
     \end{array}
   \right),
    \end{equation*}
    where $\vec{n}=n+{N}/{2}$.

Define the \pf s $$\mathbf{q}^{l}_{\pm[1]}\srr \mathbf{h}^{j}_{[1]}\pm i\,\mathbf{g}^{j}_{[1]}=\psi^{2r}_{[1],l}\pm i\,\f^{2r}_{[1],l}={\Psi}^{2r}_{\pm [1],l}, \quad l=0,1.$$
\ehh{qa_df} implies that their \fg s are
\begin{eqnarray*}\label{Psi1_dfp}
  \hat{q}^{0}_{+[1]}[n]=\left\{
               \begin{array}{ll}
  (1+i)\sqrt{2}, & \hbox{if $n=0$;} \\
                  2\b[n], & \hbox{if $0< n< N/2$;} \\
               0 & \hbox{if $ N/2\leq n<N$,}
               \end{array}
             \right.\quad
 \hat{q}^{1}_{+[1]}[n]=\left\{
               \begin{array}{ll}
-(1+i)\sqrt{2}, & \hbox{if  $n= N/2$;} \\
                   2\a[n], & \hbox{if $0< n< N/2$;} \\
               0, & \hbox{if $ N/2< n\leq N$.}
               \end{array}
             \right.\\\label{Psi1_dfm}
  \hat{q}^{0}_{-[1]}[n]=\left\{
               \begin{array}{ll}
  (1-i)\sqrt{2}, & \hbox{if $n=0$;} \\
                  2\b[n], &  \hbox{if $ N/2<n<N$,} \\
               0 & \hbox{if $0< n\leq  N/2$;}
               \end{array}
             \right.\quad
 \hat{q}^{1}_{-[1]}[n]=\left\{
               \begin{array}{ll}
-(1-i)\sqrt{2}, & \hbox{if  $n= N/2$;} \\
                   2\a[n] & \hbox{if $ N/2< n\leq  N$;} \\
               0, & \hbox{if $0\leq n< N/2$.}
               \end{array}
             \right.
\end{eqnarray*}
Thus, the \aa\ \mv ces for the \pf s $\mathbf{q}^{l}_{\pm[1]}$ are
\begin{eqnarray}\label{aa_modma10p}
 \tilde{\mathbf{M}}_{+}^{q}[n]&=& \left(
                         \begin{array}{cc}
                            \hat{q}_{+[1]}^{0}[n] &   0\\
                             \hat{q}^{1}_{+[1]}[n] &   -\sqrt{2}(1+i)\,\da[n-N/2]\\
                         \end{array}
                       \right)= \tilde{\mathbf{M}}[n]+i\, \tilde{\mathbf{M}}^{c}[n],\\\label{aa_modma10m}
 \tilde{\mathbf{M}}_{-}^{q}[n]&=& \left(
                         \begin{array}{cc}
                     (1-i)\sqrt{2}\da[n]      &   \hat{q}_{-[1]}^{0}[n] \\
                             0 &  \hat{q}^{1}_{-[1]}[n]\\
                         \end{array}
                       \right)= \tilde{\mathbf{M}}[n]-i\, \tilde{\mathbf{M}}^{c}[n],
  \end{eqnarray}
   where the \mv x $\tilde{\mathbf{M}}[n]$ is  defined in Eq. \rf{aa_modma10} and $\tilde{\mathbf{M}}^{c}[n]$ is  defined in \eh{aa_modma11}.
Application of the matrices $\tilde{\mathbf{M}}_{\pm}^{q}[n]$ to the \v\ $( \hat{x}[n] ,
       \hat{x}[\vec{n}])^{T}$ produces the \v s
\begin{equation*}\label{mod_decAQ1}
    \left(
     \begin{array}{c}
       \hat{z}_{\pm[1]}^{0}[n]_{1} \\
         \hat{z}_{\pm[1]}^{1}[n]_{1}\\
     \end{array}
   \right)=\frac{1}{2}( \tilde{\mathbf{M}}_{\pm}^{q}[n])^{*}\cdot \left(
     \begin{array}{l}
      \hat{x}[n] \\
       \hat{x}[\vec{n}]
     \end{array}
   \right)= \left(
     \begin{array}{c}
       \hat{y}_{[1]}^{0}[n]_{1} \\
         \hat{y}_{[1]}^{1}[n]_{1}\\
     \end{array}
   \right)\mp i\, \left(
     \begin{array}{c}
       \hat{c}_{[1]}^{0}[n]_{1} \\
         \hat{c}_{[1]}^{1}[n]_{1}\\
     \end{array}
   \right).
 \end{equation*}

Define the matrices ${\mathbf{M}}_{\pm}^{q}[n]\srr\tilde{\mathbf{M}}_{\pm}^{q}[n]={\mathbf{M}}[n]\pm i\,{\mathbf{M}}^{c}[n]$ and apply these matrices to the \v s\\ $(\hat{z}_{\pm[1]}^{0}[n]_{1} ,
         \hat{z}_{\pm[1]}^{1}[n]_{1})^{T}$. Here   the \mv x  ${\mathbf{M}}[n]$ is defined in Eq.  \rf{sa_modma10} and ${\mathbf{M}}^{c}[n]$ is  defined in \eh{aa_modma11}.
\bpp\label{pro:Mq_z}The following relations hold
\begin{eqnarray*}\label{Mq_z1}
 && {\mathbf{M}}_{\pm}^{q}[n]\cdot \left(
     \begin{array}{c}
        \hat{z}_{\pm[1]}^{0}[n]_{1} \\
         \hat{z}_{\pm[1]}^{1}[n]_{1}\\
     \end{array}
   \right)=\mathbf{M}[n]\cdot \left(
     \begin{array}{c}
        \hat{y}_{[1]}^{0}[n]_{1} \\
         \hat{y}_{[1]}^{1}[n]_{1}\\
     \end{array}
   \right) +{\mathbf{M}}^{c}[n]\cdot \left(
     \begin{array}{c}
        \hat{c}_{[1]}^{0}[n]_{1} \\
         \hat{c}_{[1]}^{1}[n]_{1}\\
     \end{array}
   \right) \\\nn&&\pm i \left( \mathbf{M}^{c}[n]\cdot \left(
     \begin{array}{c}
        \hat{y}_{[1]}^{0}[n]_{1} \\
         \hat{y}_{[1]}^{1}[n]_{1}\\
     \end{array}
   \right)   -\mathbf{M}[n]\cdot \left(
     \begin{array}{c}
        \hat{c}_{[1]}^{0}[n]_{1} \\
         \hat{c}_{[1]}^{1}[n]_{1}\\
     \end{array}
   \right)     \right)
    \\\label{Mq_z2}&&=2\left(\left(
     \begin{array}{l}
      \hat{x}[n] \\
       \hat{x}[n+N/2]
     \end{array}
   \right)
{\pm}i \,\left(
    \begin{array}{l}
      \hat{h}[n] \\
       \hat{h}[n+N/2]
     \end{array}
   \right)\right)
=2\left(
     \begin{array}{l}
      \hat{\bar{x}}_{\pm}[n] \\
       \hat{\bar{x}}_{\pm}[n+N/2]
     \end{array}
   \right)
,
\end{eqnarray*}
where $\mathbf{h}$ is the HT of the \ss\ $\mathbf{x}\in \Pi[N]$ and  $\mathbf{\bar{x}}_{\pm}$ are the \az\ \ss s associated with  $\mathbf{x}$.
\epp
\proof In Appendix section \ref{sec:ap2}
\bdd\label{def:a_mvs} The matrices  $\tilde{\mathbf{M}}_{\pm}^{q}[n]$ and  ${\mathbf{M}}_{\pm}^{q}[n]$ are called the \aa\ and \sa\ \mv ces for the qWP \t, respectively.\edd
\bcc\label{cor:a_mvs}Successive application of \fb s defined by the \aa\ and \sa\ \mv ces  $\tilde{\mathbf{M}}_{\pm}^{q}[n]$ and  ${\mathbf{M}}_{\pm}^{q}[n]$ to a \ss\ $\mathbf{x}\in \Pi[N]$ produces the \az\ \ss s  $\mathbf{\bar{x}}_{\pm}$ associated with  $\mathbf{x}$.\ecc
\bcc\label{cor:on_sys}A \ss\ $\mathbf{x}\in \Pi[N]$  is \ry ed by the redundant  \on\ system
\begin{eqnarray*}\label{on_sys}
         x[k]&=&\frac{1}{2}\sum_{\mu=0}^{1}\sum_{l=0}^{N/2-1}\left(y_{[1]}^{\mu} [l]\psi_{[1],\mu}^{2r}[k-2l] +
          c_{[1]}^{\mu} [l] \f_{[1],\mu}^{2r}[k-2l]\right),\\\label{on_sys_yc}
          y_{[1]}^{\mu} [l]&=& \left\langle \mathbf{x},\, \psi_{[1],\mu}^{2r}[\cdot-2l]   \right\rangle,\quad
         c_{[1]}^{\mu} [l]= \left\langle \mathbf{x},\, \f_{[1],\mu}^{2r}[\cdot-2l]   \right\rangle.
         \end{eqnarray*}
         Thus, the system
        \begin{equation*}\label{on_sysTF}
\mathbf{F} \srr  \left\{\psi_{[1],0}^{2r}[\cdot-2l] \right\}\bigoplus \left\{\psi_{[1],1}^{2r}[\cdot-2l] \right\}\bigcup\left\{\f_{[1],0}^{2r}[\cdot-2l] \right\}\bigoplus \left\{\f_{[1],1}^{2r}[\cdot-2l] \right\}
\end{equation*}
form a tight frame of the space $ \Pi[N]$.
\ecc
\subsection{Multi-level \t s}\label{sec:ss42}

%\subsubsection{Second-level \t s}\label{sec:sss421}
It was explained in Section \ref{sec:sss242} that the second-level \t\ \c s $\mathbf{y}_{[2]}^{\rr}$ are
\begin{eqnarray*}
 {y}_{[2]}^{\rr} [l]&=&   \sum_{n=0}^{N-1}x[n]\,{\psi}_{[2],\rr}^{2r}[n-4l], \quad {\psi}_{[2],\rr}^{2r}[n]  =\sum_{k=0}^{N/2-1}{h}_{[2]}^{\mu}[k] \, {\psi}_{[1],\la}^{2r}[n-2k]\Longrightarrow\\
   {y}_{[2]}^{\rr} [l]&=&  \sum_{k=0}^{N/2-1}h_{[2]}^{\mu}[k-2l] \, y_{[1]}^{\la}[k], \quad \la,\mu=0,1,\;\rr=\left\{
                                                                                                                               \begin{array}{ll}
                                                                                                                                 \mu, & \hbox{if $\la=0$ ;} \\
                                                                                                                                 3-\mu, & \hbox{if $\la=1$.}
                                                                                                                               \end{array}
                                                                                                                             \right.
\end{eqnarray*}
The \fg s of the \pf s $\mathbf{h}_{[2]}^{\mu}$ are given in \eh{pfs2} and \eh{df_wq2}. Recall that $\hat{h}_{[2]}^{\mu}[n]=\hat{h}_{[1]}^{\mu}[2n]$. The direct and inverse \t s
$\mathbf{y}_{[1]}^{\la}\longleftrightarrow\mathbf{y}_{[2]}^{2\la}\bigcup\mathbf{y}_{[2]}^{2\la+1}$ are  implemented using the \aa\ and \sa\ \mv ces $\tilde{\mathbf{M}}[2n]$ and $\mathbf{M}[2n]$,  that are defined in Eqs. \rf{aa_modma10} and \rf{sa_modma10} respectively:
\begin{equation*}\label{mod_decA20}
    \left(
     \begin{array}{c}
       \hat{y}_{[2]}^{\rr0}[n]_{2} \\
         \hat{y}_{[2]}^{\rr1}[n]_{2}\\
     \end{array}
   \right)=\frac{1}{2} \tilde{\mathbf{M}}[-2n]\cdot \left(
     \begin{array}{l}
      \hat{y}_{[1]}^{\la}[n]_{1} \\
        \hat{y}_{[1]}^{\la}[\vec{n}]_{1}
     \end{array}
   \right),\quad\left(
     \begin{array}{l}
      \hat{y}_{[1]}^{\la}[n]_{1} \\
        \hat{y}_{[1]}^{\la}[\vec{n}]_{1}
     \end{array}
   \right)={\mathbf{M}}[2n]\cdot\left(
     \begin{array}{c}
       \hat{y}_{[2]}^{\rr0}[n]_{2} \\
         \hat{y}_{[2]}^{\rr1}[n]_{2}\\
     \end{array}
   \right),
 \end{equation*}
where
\begin{equation*}\label{rr_01}
\rr0=\left\{
              \begin{array}{ll}
                0, & \hbox{if $\la=0$;} \\
                3, & \hbox{if $\la=1$,}
              \end{array}
            \right. \quad \rr1=\left\{
              \begin{array}{ll}
                1, & \hbox{if $\la=0$;} \\
                2, & \hbox{if $\la=1$,}
              \end{array}
            \right.\quad \vec{n}=n+N/4.
\end{equation*}

The second-level \t\ \c s $\mathbf{c}_{[2]}^{\rr}$ are
\begin{eqnarray*}
 {c}_{[2]}^{\rr} [l]&=&   \sum_{n=0}^{N-1}x[n]\,{\f}_{[2],\rr}^{2r}[n-4l], \quad {\f}_{[2],\rr}^{2r}[n]  =\sum_{k=0}^{N/2-1}{h}_{[2]}^{\mu}[k] \, {\f}_{[1],\la}^{2r}[n-2k]\Longrightarrow\\
   {c}_{[2]}^{\rr} [l]&=&  \sum_{k=0}^{N/2-1}h_{[2]}^{\mu}[k-2l] \, c_{[1]}^{\la}[k], \quad \la,\mu=0,1,\;\rr=\left\{
                                                                                                                               \begin{array}{ll}
                                                                                                                                 \mu, & \hbox{if $\la=0$ ;} \\
                                                                                                                                 3-\mu, & \hbox{if $\la=1$.}
                                                                                                                               \end{array}
                                                                                                                             \right.
\end{eqnarray*}
We emphasize that the \pf s $\mathbf{h}_{[2]}^{\mu}$ for the \t\ $\mathbf{c}_{[1]}^{\la}\longleftrightarrow\mathbf{c}_{[2]}^{2\la}\bigcup\mathbf{c}_{[2]}^{2\la+1}$ are the same that the \pf s for the \t\ $\mathbf{y}_{[1]}^{\la}\longleftrightarrow\mathbf{y}_{[2]}^{2\la}\bigcup\mathbf{y}_{[2]}^{2\la+1}$. Therefore, the direct and inverse \t s
$\mathbf{c}_{[1]}^{\la}\longleftrightarrow\mathbf{c}_{[2]}^{2\la}\bigcup\mathbf{c}_{[2]}^{2\la+1}$ are  implemented using the same \aa\ and \sa\ \mv ces $\tilde{\mathbf{M}}[2n]$ and $\mathbf{M}[2n]$. Apparently, it is the case also for the \t s $\mathbf{z}_{\pm[1]}^{\la}\longleftrightarrow\mathbf{z}_{\pm[2]}^{2\la}\bigcup\mathbf{z}_{\pm[2]}^{2\la+1}$. The \t s to subsequent \d\  levels are implemented in an iterative way:

\begin{eqnarray*}\label{mod_decA2m}
    \left(
     \begin{array}{c}
       \hat{z}_{\pm[m+1]}^{\rr0}[n]_{m+1} \\
         \hat{z}_{\pm[m+1]}^{\rr1}[n]_{m+1}\\
     \end{array}
   \right)&=&\frac{1}{2} \tilde{\mathbf{M}}[-2^{m}n]\cdot \left(
     \begin{array}{l}
      \hat{z}_{\pm[m]}^{\la}[n]_{m} \\
        \hat{z}_{\pm[m]}^{\la}[\vec{n}]_{m}
     \end{array}
   \right),\\\nn\left(
     \begin{array}{l}
      \hat{z}_{\pm[m]}^{\la}[n]_{m} \\
        \hat{z}_{\pm[m]}^{\la}[\vec{n}]_{m}
     \end{array}
   \right)&=&{\mathbf{M}}[2^{m}n]\cdot\left(
     \begin{array}{c}
       \hat{z}_{\pm[m+1]}^{\rr0}[n]_{m+1} \\
         \hat{z}_{\pm[m+1]}^{\rr1}[n]_{m+1}\\
     \end{array}
   \right),
 \end{eqnarray*}
where
\(
\rr0=\left\{
              \begin{array}{ll}
                2\la, & \hbox{if $\la$ is even;} \\
                2\la+1, & \hbox{if $\la$ is odd,}
              \end{array}
            \right. \) and vice versa for $\rr1$, $\vec{n}=n+N/2^{m+1}$ and $m=1,...,M$.
By the application of the inverse DFT to the arrays $\left\{ \hat{z}_{\pm[m+1]}^{\rr}[n]_{m+1}\right\}$, we get the arrays\\ $\left\{ z_{\pm[m+1]}^{\rr}[k]=y_{[m+1]}^{\rr}[k]\pm i\,c_{[m+1]}^{\rr}[k]\right\}$ of the \t\ \c s with the qWPs $\Psi^{2r}_{\pm[m+1],\rr}$.
\br\label{rem:zyc_cyz}We stress that by operating on  the \t\  \c s $\left\{ z_{\pm[m]}^{\rr}[k]\right\}$, we simultaneously operate on  the  arrays $\left\{ y_{[m]}^{\rr}[k]\right\}$ and $\left\{ c_{[m]}^{\rr}[k]\right\}$, which are the \c s for the \t s with the WPs $\psi^{2r}_{[m],\rr}$ and cWPs $\f^{2r}_{[m],\rr}$, respectively. The execution speed of the \t\ with the qWPs  $\left\{ \Psi_{\pm[m]}^{2r}\right\}= \psi_{[m]}^{2r}\pm i \f_{[m]}^{2r}$ is the same as the speed of the \t s with either WPs $\left\{ \psi_{[m]}^{2r}\right\}$ or cWPs $\left\{ \f_{[m]}^{2r}\right\}$.\er
 The \t s are executed in the spectral domain using the  FFT by the application of critically sampled two-channel \fb s to the half-band \sp al components $(\hat{x}[n],\hat{x}[n+N/2])^{T}$  of a \ss.

\begin{figure}[H]
\begin{center}
\resizebox{15cm}{8cm}{
\includegraphics{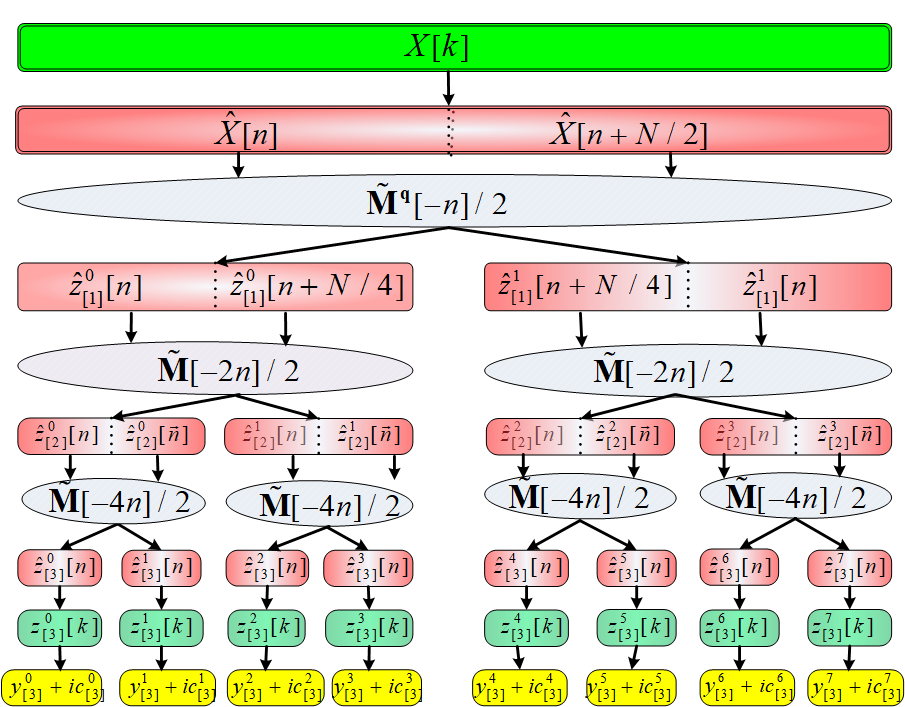}
%\hfill
%\vline
%\hfill
%\includegraphics{png/dia_wq_SAh.png}
}
\end{center}
\caption{ Forward qWTP of a \ss\ $\mathbf{X}$ down to the third \d\ level with quasi-\az\ \wq s. Here $\vec{n}$ means ${n}+N/8$ }
     \label{dia_wq_A}
 \end{figure}

\begin{figure}[H]
\begin{center}
\resizebox{15cm}{8cm}{
%\includegraphics{png/dia_wq_aA.png}
%\hfill
%\vline
%\hfill
\includegraphics{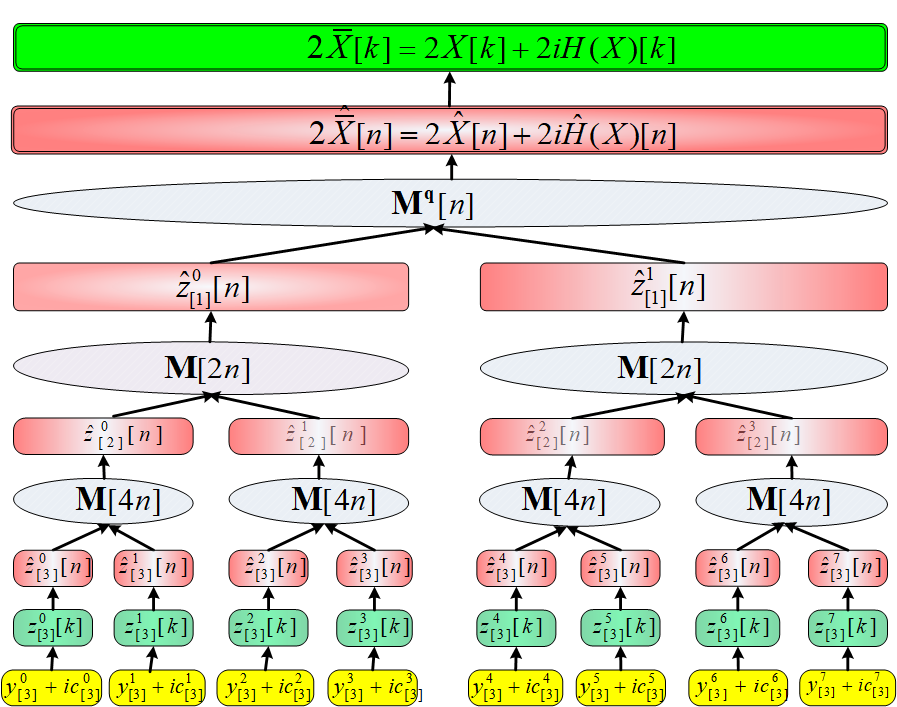}
}
\end{center}
\caption{  Inverse  qWTP from the  \t\  \c s from the third \d\ level that results in restoration of the \ss\  $\mathbf{X}$ and its HT   $H(\mathbf{X})$ }
     \label{dia_wq_S}
 \end{figure}

The diagrams in Figs. \ref{dia_wq_A} and \ref{dia_wq_S} illustrate the three-level forward and inverse qWPTs of a \ss\ with quasi-\az\ \wq s, which use the \aa\ $\tilde{\mathbf{M}}^{q}[n]$ and the \sa\ ${\mathbf{M}}^{q}[n]$ \mv ces, respectively,  for the \t s to and from the first \d\ level, respectively, and the \mv ces  $\tilde{\mathbf{M}}[2^{m}n]$ and  ${\mathbf{M}}[2^{m}n]$ for the subsequent levels.

%\paragraph{Brief discussion}
\br\label{rem:recon}The \d\ of a \ss\  $\mathbf{x}\in\Pi[N]$ down to the $M$-th level produces  $2MN$ \t\ \c s $\left\{ y_{[m]}^{\rr}[k]\right\}\bigcup\left\{ c_{[m]}^{\rr}[k]\right\}$. Such a redundancy provides many options for the \ss\ \r. Some of them are listed below.
\begin{itemize}
  \item A basis compiled from either WPs $\left\{ \psi_{[m]}^{2r}\right\}$ or $\left\{ \f_{[m]}^{2r}\right\}$.
\begin{itemize}
  \item Wavelet basis.
  \item Best bases \cite{coiw1}, Local discriminant bases \cite{sai,sai2}.
  \item WPs from a single \d\ level.
\end{itemize}
  \item Combination of bases compiled from both $\left\{ \psi_{[m]}^{2r}\right\}$ and $\left\{ \f_{[m]}^{2r}\right\}$ WPs generates a tight frame of the space $\Pi[N]$  with redundancy rate 2. The bases for  $\left\{ \psi_{[m]}^{2r}\right\}$ and $\left\{ \f_{[m]}^{2r}\right\}$ can have a different structure.
  \item Frames with increased redundancy rate. For example, a combined \r\ from several \d\ levels.
\end{itemize}
 \er
The collection of WPs $\left\{ \psi_{[m]}^{2r}\right\}$ and cWPs $\left\{ \f_{[m]}^{2r}\right\}$, which originate from \ds  s of \df t orders $2r$, provides a variety of \we s that are (anti)\sy ic, well localized in time domain. Their DFT \sp a are flat and the \sp a shapes tend to rectangles when  the order $2r$ increases. Therefore, they can be utilized as a collection of band-pass \fr s which produce a refined split of the \ff\ domain into bands of \df t widths. The (c)WPs can be used as testing \we s for the \ss\ \aa, such as a dictionary for the Matching Pursuit procedures \cite{mal,azk_MP}.

\section{Two-dimensional complex \wq s}\label{sec:s5}
A standard design scheme for 2D \wq s is outlined in  Section \ref{sec:ss25}. The 2D \wq s are defined as the tensor products of 1D WPs such that
\begin{equation*}\label{psipsi0}
  \psi_{[m],j ,l}^{2r}[k,n]=\psi_{[m],j}^{2r}[k]\,\psi_{[m], l}^{2r}[n].
\end{equation*}
The $2^{m}$-sample shifts of the WPs $\left\{\psi_{[m],j ,l}^{2r}\right\},\;j , l=0,...,2^{m}-1,$ in both directions form an \on\ basis for the space $\Pi[N,N]$ of arrays that are $N$-\p\ in both directions. The DFT \sp um  of such a WP is  concentrated in four \sy ic spots in the \ff\ domain as it is seen in Fig. \ref{psifpsi2_2}.

Similar properties are inherent to the 2D cWPs such that
\begin{equation*}\label{phiphi}
  \f_{[m],j ,l}^{2r}[k,n]=\f_{[m],j}^{2r}[k]\,\f_{[m], l}^{2r}[n].
\end{equation*}

\subsection{Design of  2D directional WPs}\label{sec:ss51}
\subsubsection{2D complex WPs and their \sp a }\label{sec:sss511}
The WPs  $\left\{\psi_{[m],j ,l}^{2r}\right\}$ as well as the cWPs  $\left\{\f_{[m],j ,l}^{2r}\right\}$  lack  the directionality property  which is needed in many applications that process 2D data.  However, real-valued 2D \wq s oriented in multiple directions  can be
derived from  tensor  products of complex quasi-\az\ qWPs $\Psi_{\pm[m],\rr}^{2r}$.

\begin{figure}[H]
\centering
\includegraphics[width=6.2in]{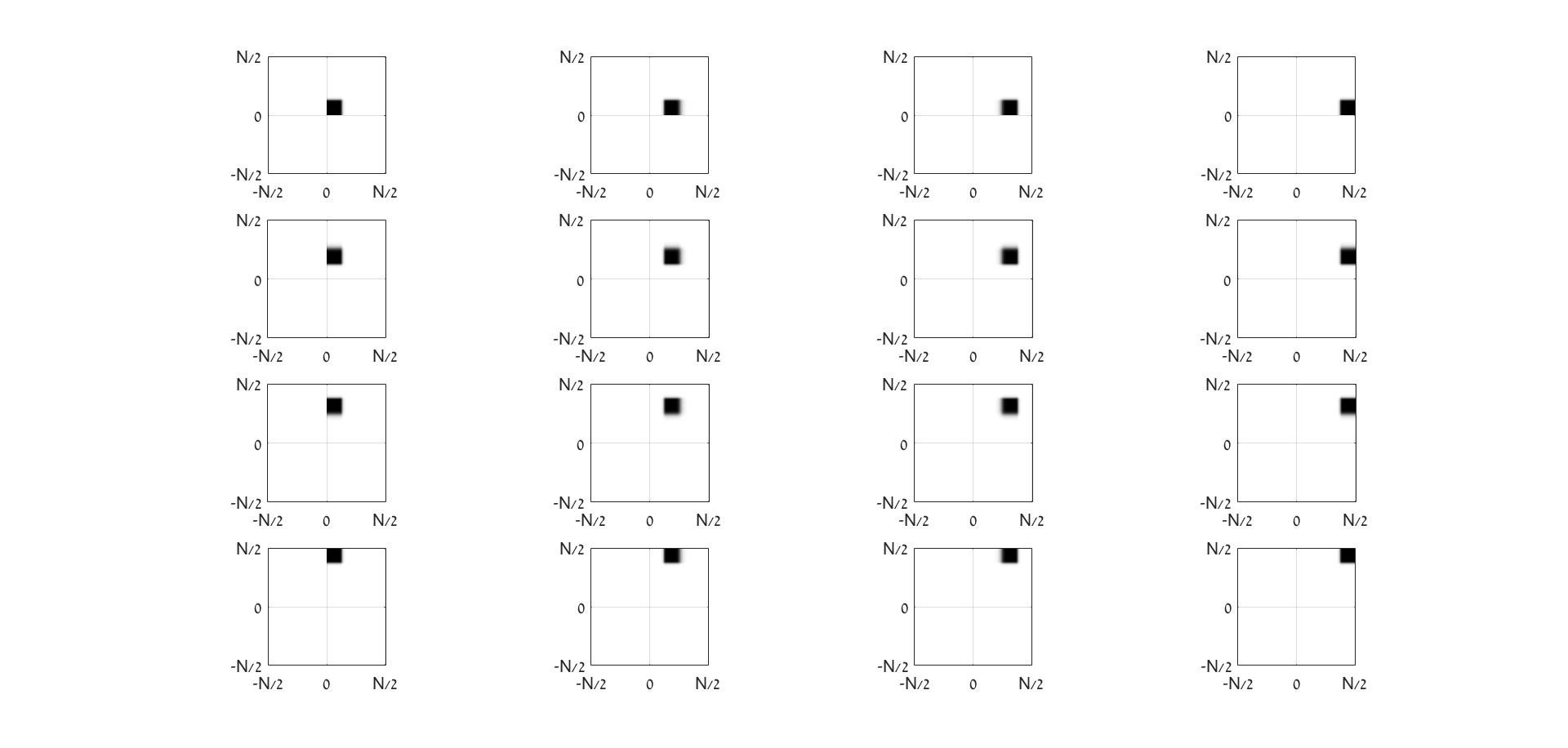}
%\hfil
%\includegraphics[width=3.2in]{png/fpm_2.png}%
\caption{Magnitude \sp a of 2D qWPs $\Psi_{++[2],j ,l}^{10}$  from the second \d\ level}
\label{fpp_2}
\end{figure}

\begin{figure}[H]
\centering
%\includegraphics[width=3.2in]{png/fpp_2.png}
%\hfil
\includegraphics[width=6.2in]{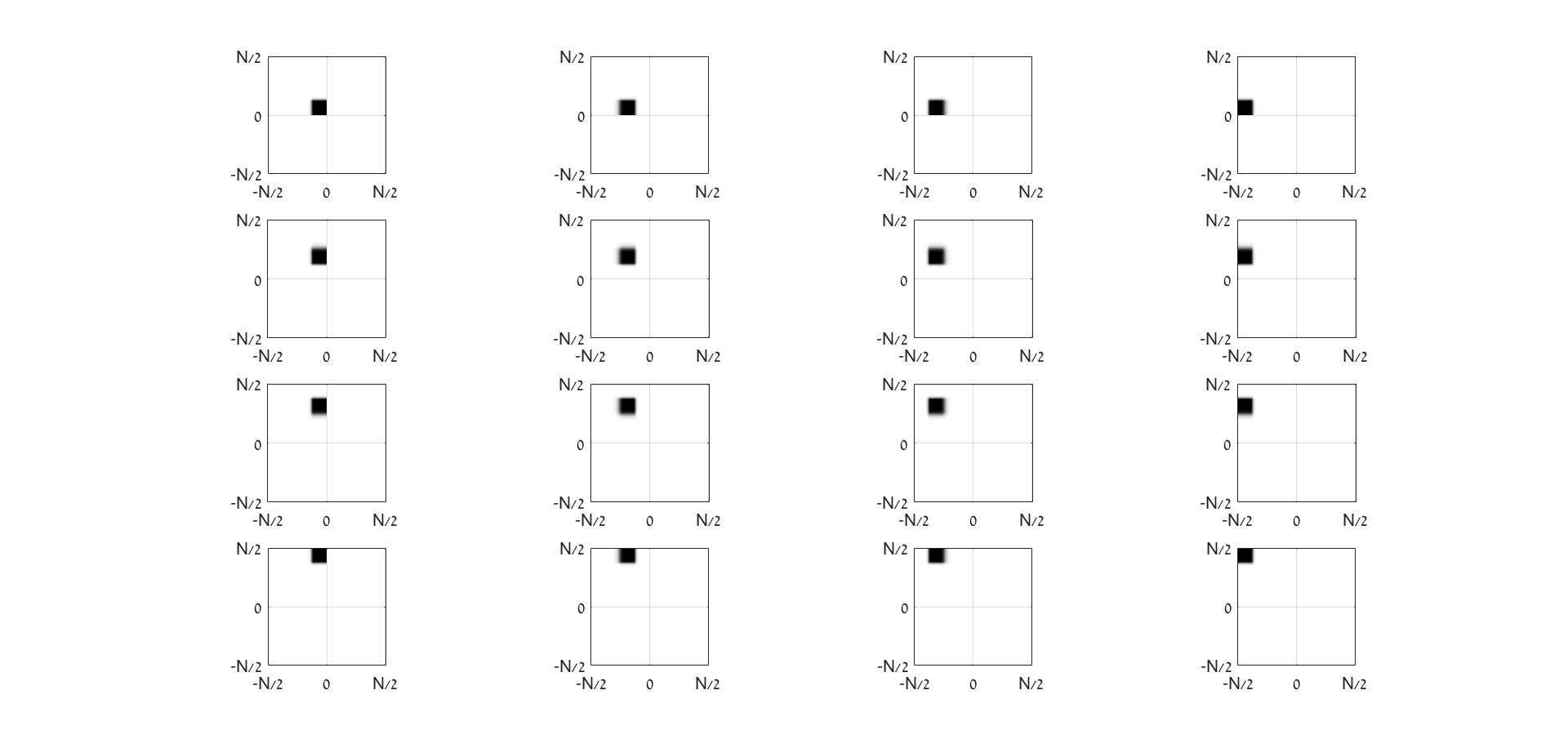}%
\caption{Magnitude \sp a of 2D qWPs  $\Psi_{+-[2],j ,l}^{10}$  from the second \d\ level}
\label{fmp_2}
\end{figure}

The complex 2D qWPs are defined  as follows:
\begin{eqnarray*} \label{qwp_2d}
\Psi_{++[m],j , l}^{2r}[k,n] &\srr& \Psi_{+[m],j}^{2r}[k]\,\Psi_{+[m], l}^{2r}[n], \\\nn
  \Psi_{+-[m],j ,l}^{2r}[k,n] &\srr& \Psi_{+[m],j}^{2r}[k]\,\Psi_{-[m], l}^{2r}[n],
\end{eqnarray*}
where $  m=1,...,M,\;j ,l=0,...,2^{m}-1,$ and $k ,n=-N/2,...,N/2-1$.
The real and imaginary parts of these 2D qWPs are
\begin{equation}
\label{vt_pm}
\begin{array}{lll}
 \vt_{+[m],j ,l}^{2r}[k,n] &\srr& \mathfrak{Re}(\Psi_{++[m],j ,l}^{2r}[k,n]) =  \psi_{[m],j ,l}^{2r}[k,n]-\f_{[m],j ,l}^{2r}[k,n], \\
\vt_{-[m],j ,l}^{2r}[k,n] &\srr&  \mathfrak{Re}(\Psi_{+-[m],j ,l}^{2r}[k,n]) =  \psi_{[m],j ,l}^{2r}[k,n]+\f_{[m],j ,l}^{2r}[k,n],\\%\label{th_pm}
\end{array}
\end{equation}
\begin{equation}
\label{th_pm}
\begin{array}{lll}
 \th_{+[m],j ,l}^{2r}[k,n] &\srr& \mathfrak{Im}(\Psi_{++[m],j ,l}^{2r}[k,n]) =  \psi_{[m],j}^{2r}[k]\,\f_{[m], l}^{2r}[n]+\f_{[m],j}^{2r}[k]\,\psi_{[m], l}^{2r}[n], \\
\th_{-[m],j ,l}^{2r}[k,n] &\srr&  \mathfrak{Im}(\Psi_{+-[m],j ,l}^{2r}[k,n]) =   \f_{[m],j}^{2r}[k]\,\psi_{[m], l}^{2r}[n]-\psi_{[m],j}^{2r}[k]\,\f_{[m], l}^{2r}[n].
\end{array}
\end{equation}
The DFT \sp a of the 2D qWPs $\Psi_{++[m],j ,l}^{2r},\;j ,l=0,...,2^{m}-1,$ are the tensor products of the one-sided \sp a of the qWPs:
\begin{equation*}\label{fPsi_pp}
\hat{ \Psi}_{++[m],j ,l}^{2r}[p,q] =\hat{ \Psi}_{+[m],j}^{2r}[p]\,\hat{\Psi}_{+[m], l}^{2r}[q]
\end{equation*}
and, as such,  they fill the  quadrant $k ,n=0,...,N/2-1$ of the \ff\ domain, while the \sp a of $\Psi_{+-[m],j ,l}^{2r},\;j ,l=0,...,2^{m}-1,$ fill the  quadrant $k =0,...,N/2-1,\;n=-N/2,...,-1$. Figures \ref{fpp_2} and \ref{fmp_2} display the magnitude \sp a of the tenth-order 2D qWPs $\Psi_{++[2],j ,l}^{10}$ and $\Psi_{+-[2],j ,l}^{10}$ from the second \d\ level, respectively.

\br\label{rem:dif_lev}The 2D qWPs $\Psi_{+\pm[m],j , l}^{2r}$ are the tensor products of 1D qWPs from the \d\ level $m$. However, there is no problems to design the 2D qWPs as a tensor products of 1D qWPs from \df t \d\ levels such as
\(\Psi_{+\pm[m,s],j , l}^{2r}[k,n] \srr \Psi_{+[m],j}^{2r}[k]\,\Psi_{\pm[s], l}^{2r}[n].\)
\er
\subsubsection{Directionality of real-valued 2D WPs}\label{sec:sss512}
It is seen in Fig. \ref{fpp_2} that the DFT \sp a of the  qWPs $\Psi_{+\pm[m],j ,l}^{10}$ effectively occupy relatively small squares in the \ff\ domain. For deeper \d\ levels,  sizes of the corresponding  squares decrease on geometric progression. Such configurations of the \sp a lead to the directionality of the real-valued 2D WPs $ \vt_{\pm[m],j ,l}^{2r}$ and $ \th_{\pm[m],j ,l}^{2r}$.

\begin{SCfigure}
\centering
\caption{Magnitude \sp a of 2D qWP $\Psi$ (left) and  $\mathfrak{Re}(\Psi)=\vt$  (right)}
\includegraphics[width=3.2in]{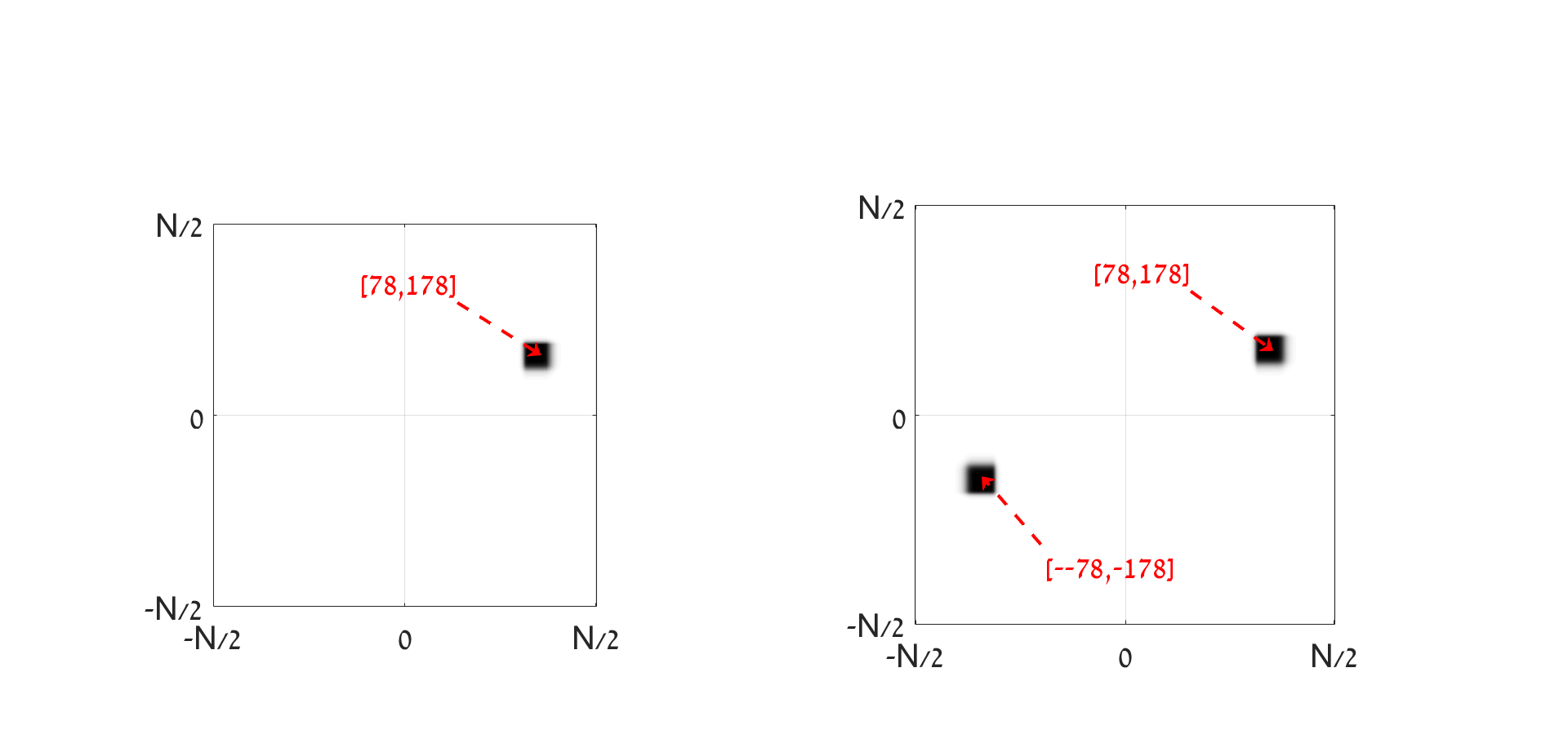}
%\hfil
%\includegraphics[width=3.2in]{png/fpm_2.png}%
\label{78_178}
\end{SCfigure}

Assume, for example, that $N=512, \;m=3,\; j=2,\;  l=5$ and denote $\Psi[k,n]\srr \Psi_{++[3],2 ,5}^{2r}[k,n]$ and $ \vt[k,n]\srr\mathfrak{Re}(\Psi[k,n])$. Its magnitude \sp um $\left|\hat{\Psi}[j,l]\right|$, displayed in Fig, \ref{78_178}, effectively occupies the square of size $40\times 40$ \emph{pixels} centered around the point $\mathbf{C}=[\k_{0},\n_{0}]$, where $\k_{0}=78, \;\n_{0}=178$. Thus, the WP $\Psi$ is \ry ed by
\begin{eqnarray*}
\label{psi78_178}
  \Psi[k,n] &=& \frac{1}{N^{2}}\sum_{\k,\n=0}^{N/2-1}\w^{k\k+n\n}\, \hat{\Psi}[\k,\n]\approx{\w^{\k_{0}k+\n_{0}n}}\,\underline{\Psi}[k,n] \\\nn
 \underline{\Psi}[k,n]  &\srr& \frac{1}{N^{2}}\sum_{\k,\n=-20}^{19}\w^{k\k+n\n} \, \hat{\Psi}[\k+\k_{0},\n+\n_{0}].
\end{eqnarray*}
Consequently, the real-valued  WP $\vt $ is \ry ed as follows:
\begin{eqnarray*}
\label{th78_178}
  \vt[k,n]  \approx{\cos\frac{2\pi(\k_{0}k+\n_{0}n)}{N}}\,\underline{\vt}[k,n] ,\quad \underline{\vt}[k,n] \srr\mathfrak{Re}(\underline{\Psi}[k,n]).
\end{eqnarray*}
The \sp um of the 2D \ss\ $\underline{\vt}$ comprises only f low frequencies in both directions and the \ss\ $\underline{\vt}$ does not have a directionality. But the 2D \ss\ $\cos\frac{2\pi(\k_{0}k+\n_{0}n)}{N}$ is oscillating in the direction of the vector $\vec{V}_{++[2],2 ,5}=178\vec{i}+78\vec{j}$. The 2D WP $\vt[k,n]$ is well localized in the spatial domain as is seen from \eh{vt_pm} and the same is true for the low-\ff\ \ss\ $\underline{\vt}$. Therefore, WP $\vt[k,n]$  can be regarded as the directional cosine modulated  by the localized low-\ff\ \ss\ $\underline{\vt}$.

\begin{figure}[H]
%\resizebox{12cm}{4cm}{
\centering
\includegraphics[width=3.2in]{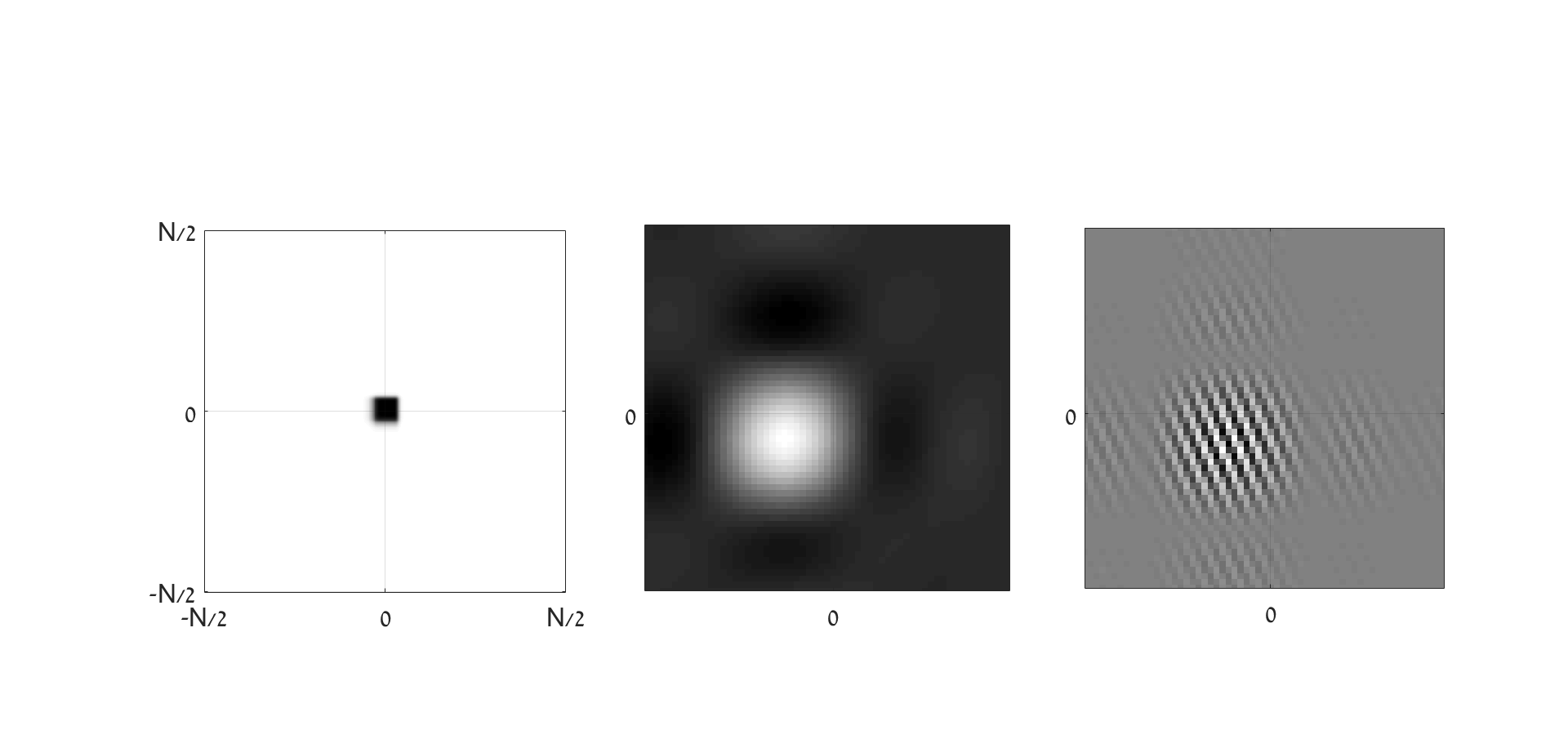}
%\includegraphics{png/AB3_2_5.png}
%}
\caption{Center: Low-\ff\ \ss\ $\underline{\vt}$. Left: Its magnitude \sp um. Right:  2D WP $\vt[k,n]$}
     \label{AB_3_2_5}
\end{figure}

The same arguments, which to some extent are similar to the discussion in Section 6.2 of \cite{bhan_zhao}, are applicable to all four real-valued 2D WPs defined in  Eqs. \rf{vt_pm} -- \rf{th_pm}. Figure \ref{AB_3_2_5} displays the low-\ff\ \ss\ $\underline{\vt}$, its magnitude \sp um and the  2D WP $\vt[k,n]$.

\begin{figure}[H]
\centering
\includegraphics[width=6.2in]{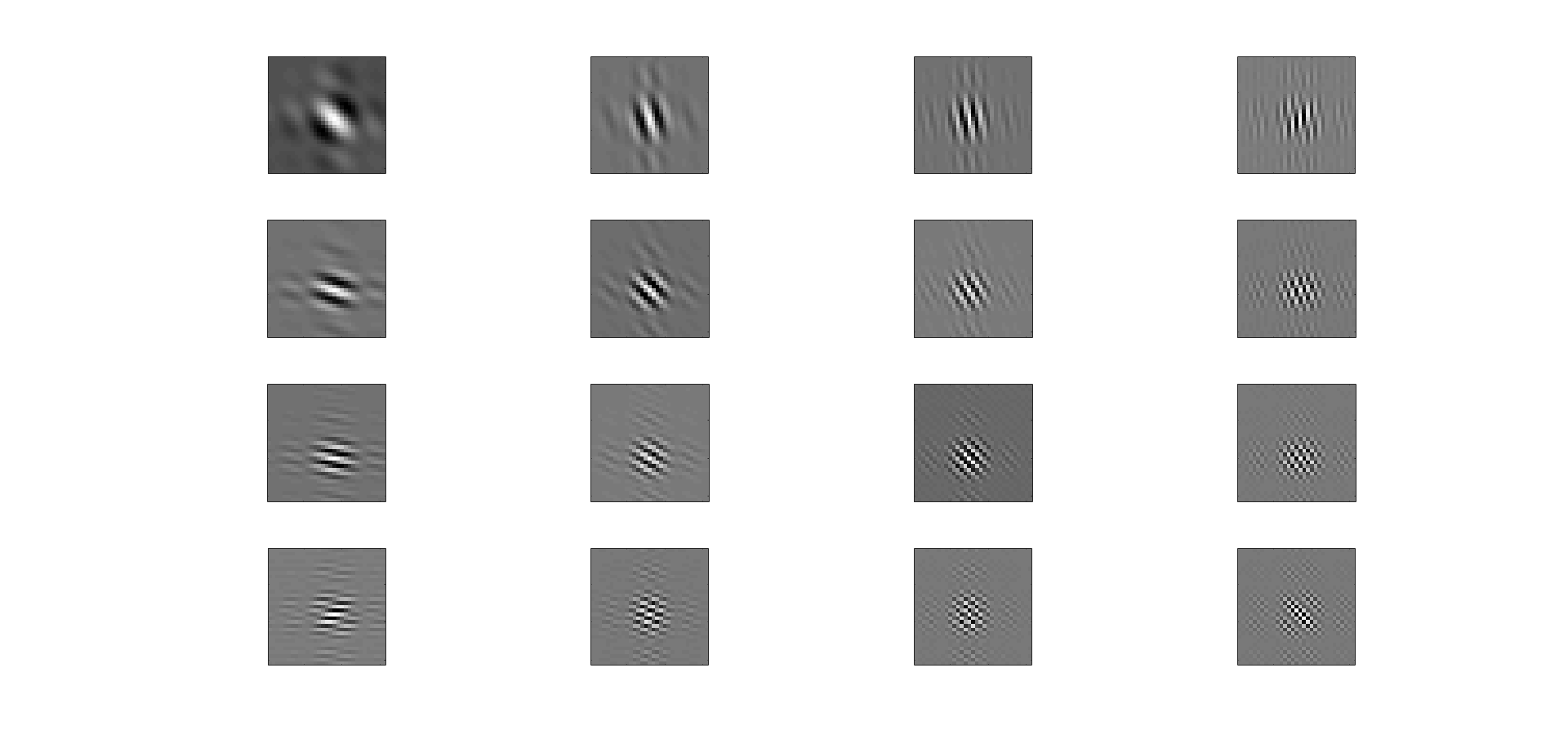}
%\hfil
%\includegraphics[width=3.2in]{png/fpp_2_2d.png}%
\caption{WPs $\vt_{+[2],j ,l}^{10}$ from the second \d\ level}
\label{pp_2_2d}
\end{figure}

\begin{figure}[H]
\centering
\includegraphics[width=5.2in]{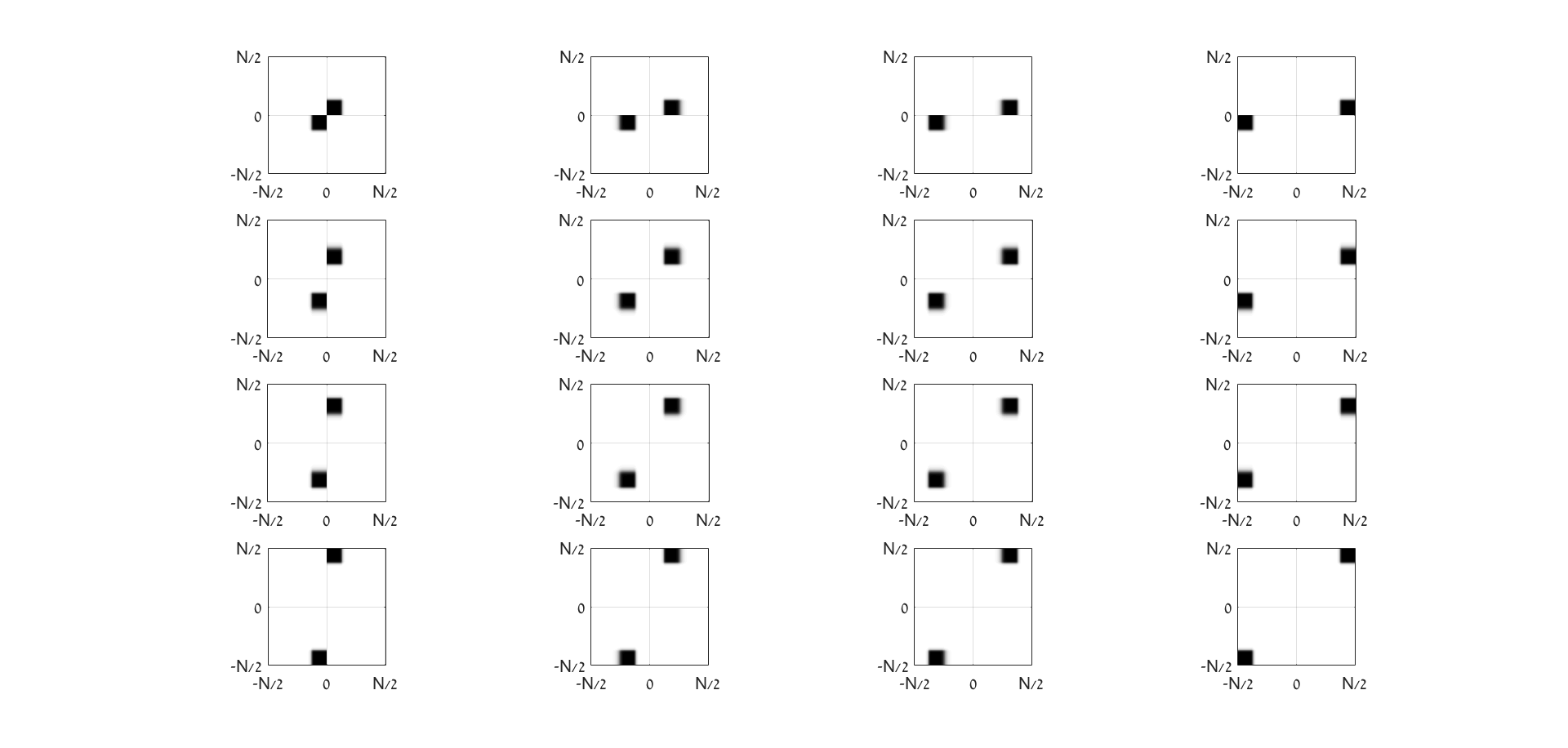}
%\hfil
%\includegraphics[width=3.2in]{png/fpp_2_2d.png}%
\caption{Magnitude \sp a  of WPs $\vt_{+[2],j ,l}^{10}$ from the second \d\ level }
\label{fpp_2_2d}
\end{figure}

\begin{figure}[H]
\centering
\includegraphics[width=6.2in]{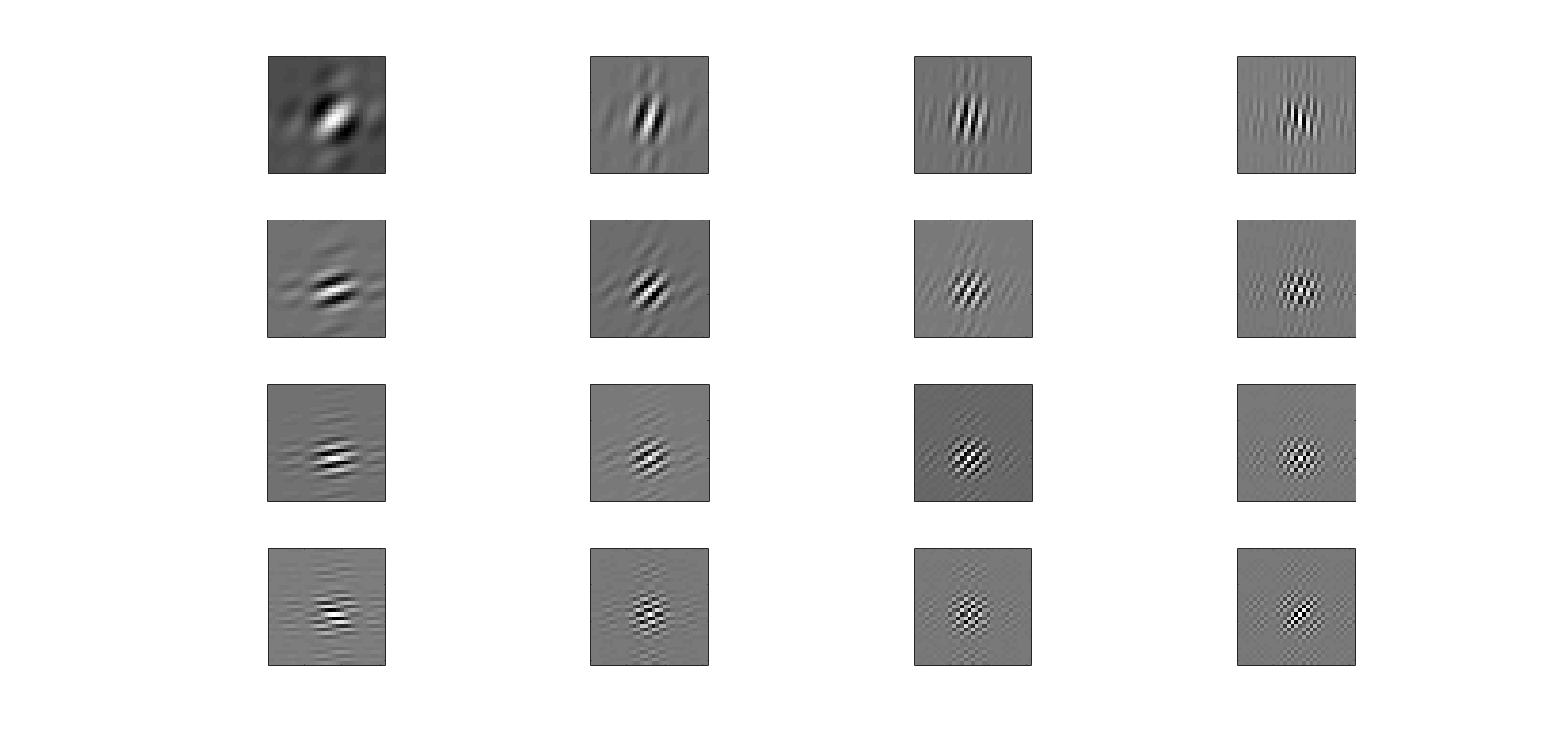}
%\hfil
%\includegraphics[width=3.2in]{png/fpp_2_2d.png}%
\caption{WPs $\vt_{-[2],j ,l}^{10}$ from the second \d\ level}
\label{pm_2_2d}
\end{figure}

\begin{figure}[H]
\centering
\includegraphics[width=5.2in]{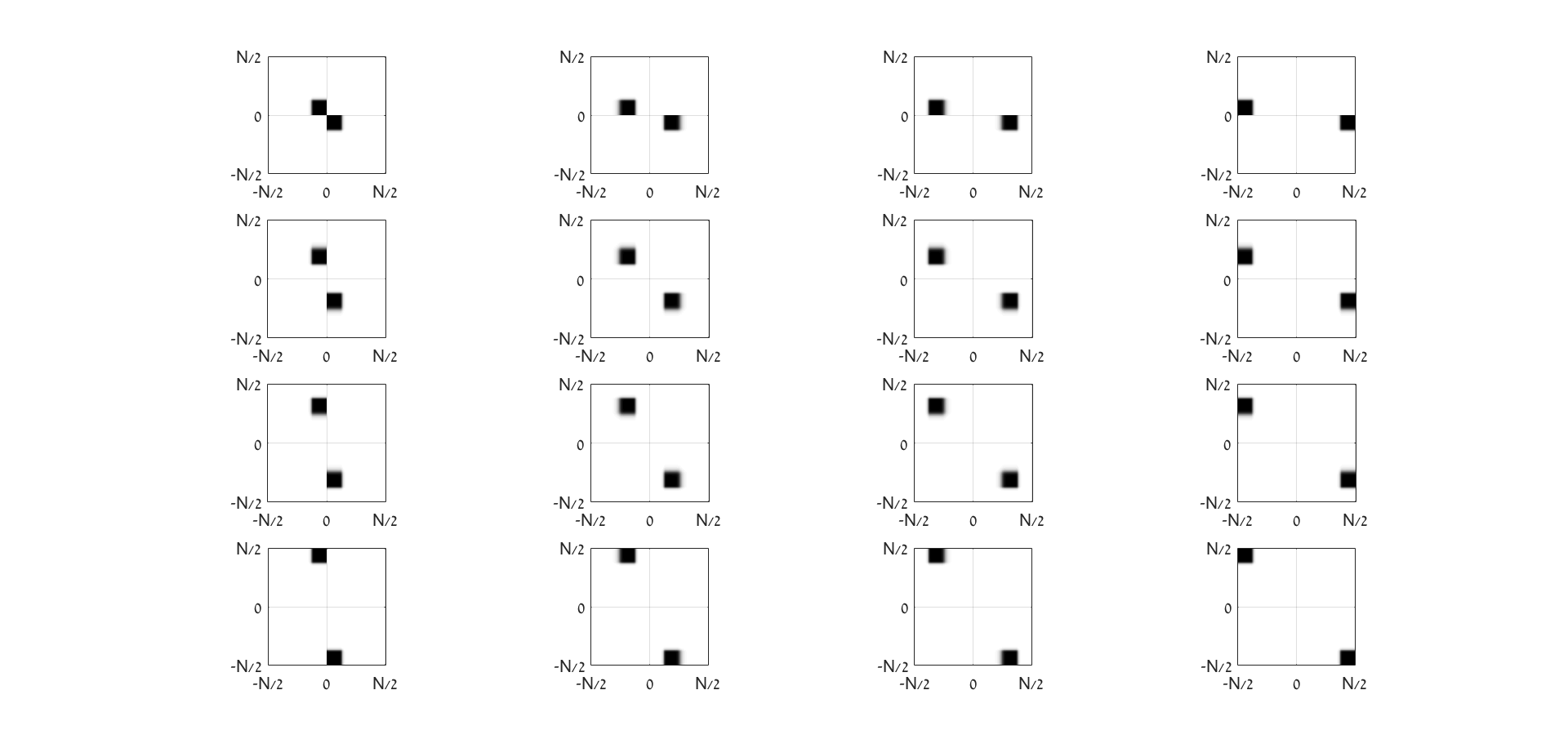}
%\hfil
%\includegraphics[width=3.2in]{png/fpp_2_2d.png}%
\caption{Magnitude \sp a  of WPs $\vt_{-[2],j ,l}^{10}$ from the second \d\ level }
\label{fpm_2_2d}
\end{figure}

Figures \ref{pp_2_2d} and  \ref{fpp_2_2d} display WPs  $\vt_{+[2],j ,l}^{10},\;j,l=0,1,2,3,$ from the second \d\ level and their magnitude \sp a, respectively.
Figures \ref{pm_2_2d} and  \ref{fpm_2_2d} display WPs  $\vt_{-[2],j ,l}^{10},\;j,l=0,1,2,3,$ from the second \d\ level and their magnitude \sp a, respectively.

\begin{figure}[H]
\centering
\includegraphics[width=6.2in]{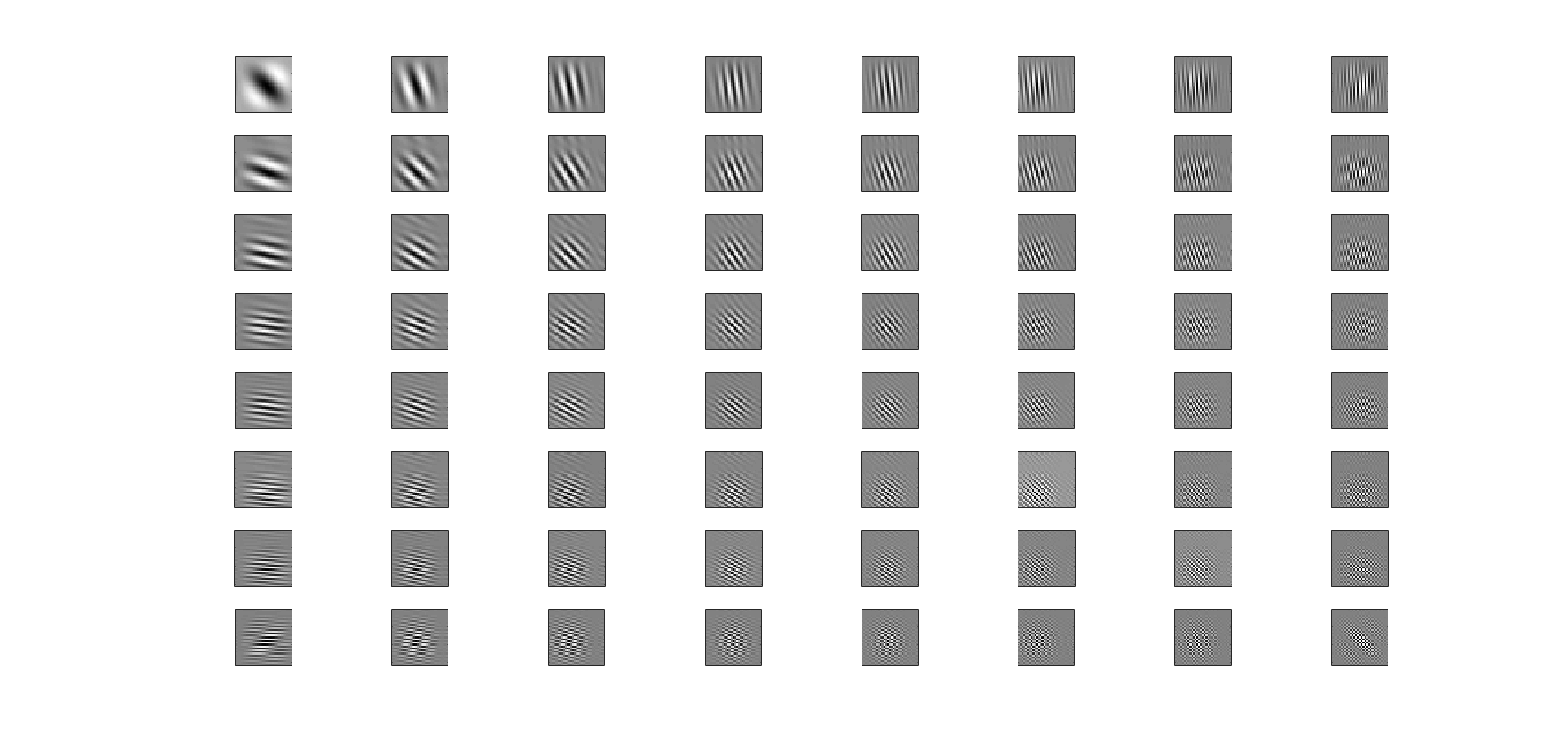}
%\hfil
%\includegraphics[width=3.2in]{png/pm_3_2d.png}%
\caption{WPs $\vt_{+[3],j ,l}^{10}$   from the third \d\ level }
\label{pp_3_2d}
\end{figure}

\begin{figure}[H]
\centering
%\includegraphics[width=3.2in]{png/pp_3_2d.png}
%\hfil
\includegraphics[width=6.2in]{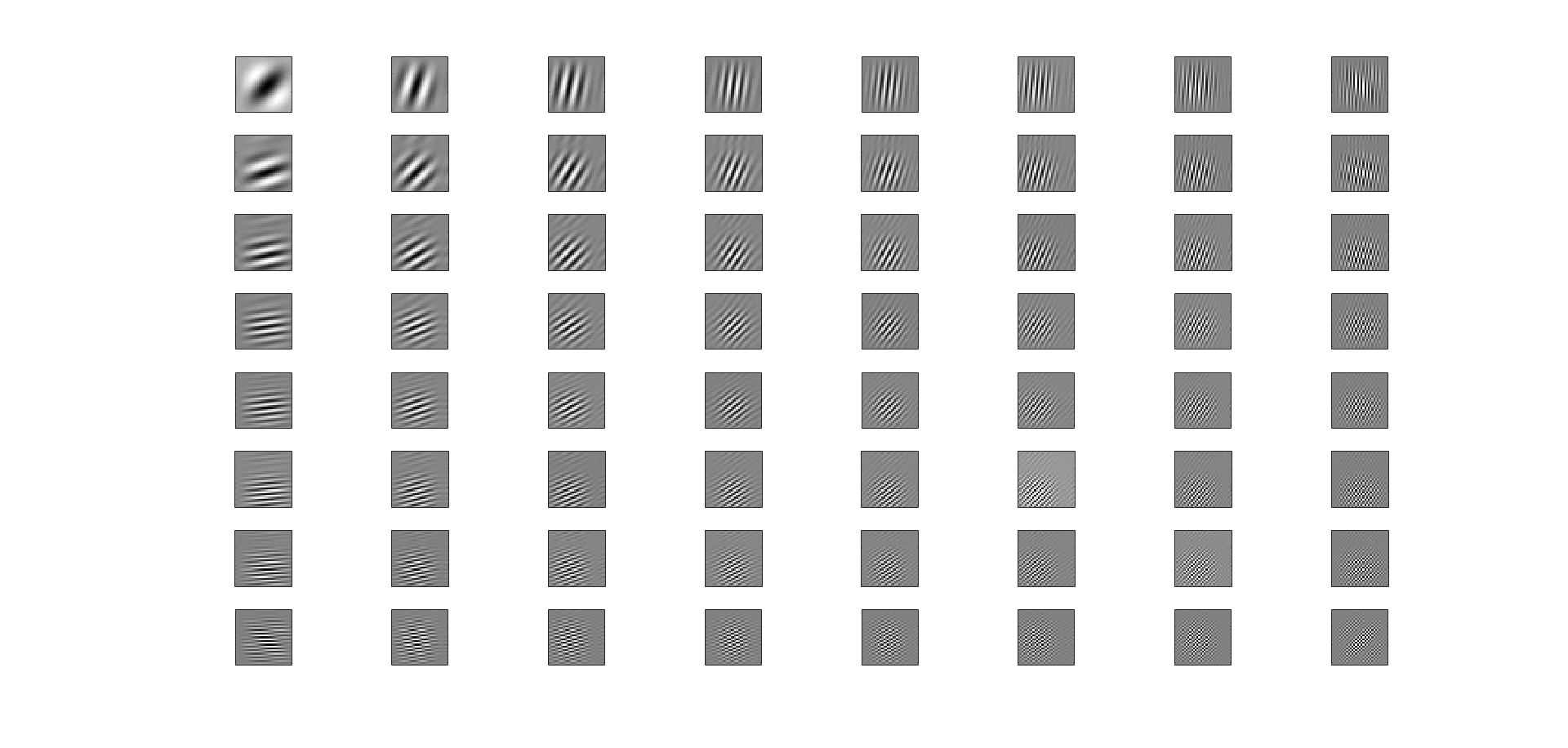}%
\caption{WPs $\vt_{-[3],j ,l}^{10}$ (right) from the third \d\ level }
\label{pm_3_2d}
\end{figure}

%\begin{figure}[H]
%\centering
%\includegraphics[width=3.2in]{png/pm_2_2d.png}
%\hfil
%\includegraphics[width=3.2in]{png/fpm_2_2d.png}%
%\caption{WPs $\vt_{-[2],j ,l}^{10}$ from the second \d\ level (left) and their magnitude \sp a (right)}
%\label{fpm_2_2d}
%\end{figure}
\br\label{direc_rem}Note that orientations of the \v s $\vec{V}_{++[m],j, l}$ and $\vec{V}_{++[m],j+1, l+1}$ are \as ely the same. These \v s   determine the orientations of the WPs $\vt_{+[m],j ,l}^{2r}$ and $\vt_{+[m],j+1 , l+1}^{2r}$, respectively. Thus, these WPs have \as ely the same orientation. Consequently, the WPs from the $m$-th \d\ level are oriented in $2^{m+1}-1$ \df t directions. The same is true for the WPs $\vt_{-[m],j ,l}^{2r}$. Thus, altogether, at the level $m$ we have WPs oriented in $2(2^{m+1}-1)$ \df t directions.  It is seen in Figs. \ref{pp_2_2d}, \ref{pm_2_2d} and in Figs. \ref{pp_3_2d}, \ref{pm_3_2d} that display the WPs $\vt_{\pm[3],j ,l}^{2r}$.\er
\section{Implementation of 2D qWP \t s}\label{sec:s6}
\begin{figure}[H]
\centering
\includegraphics[width=3.0in]{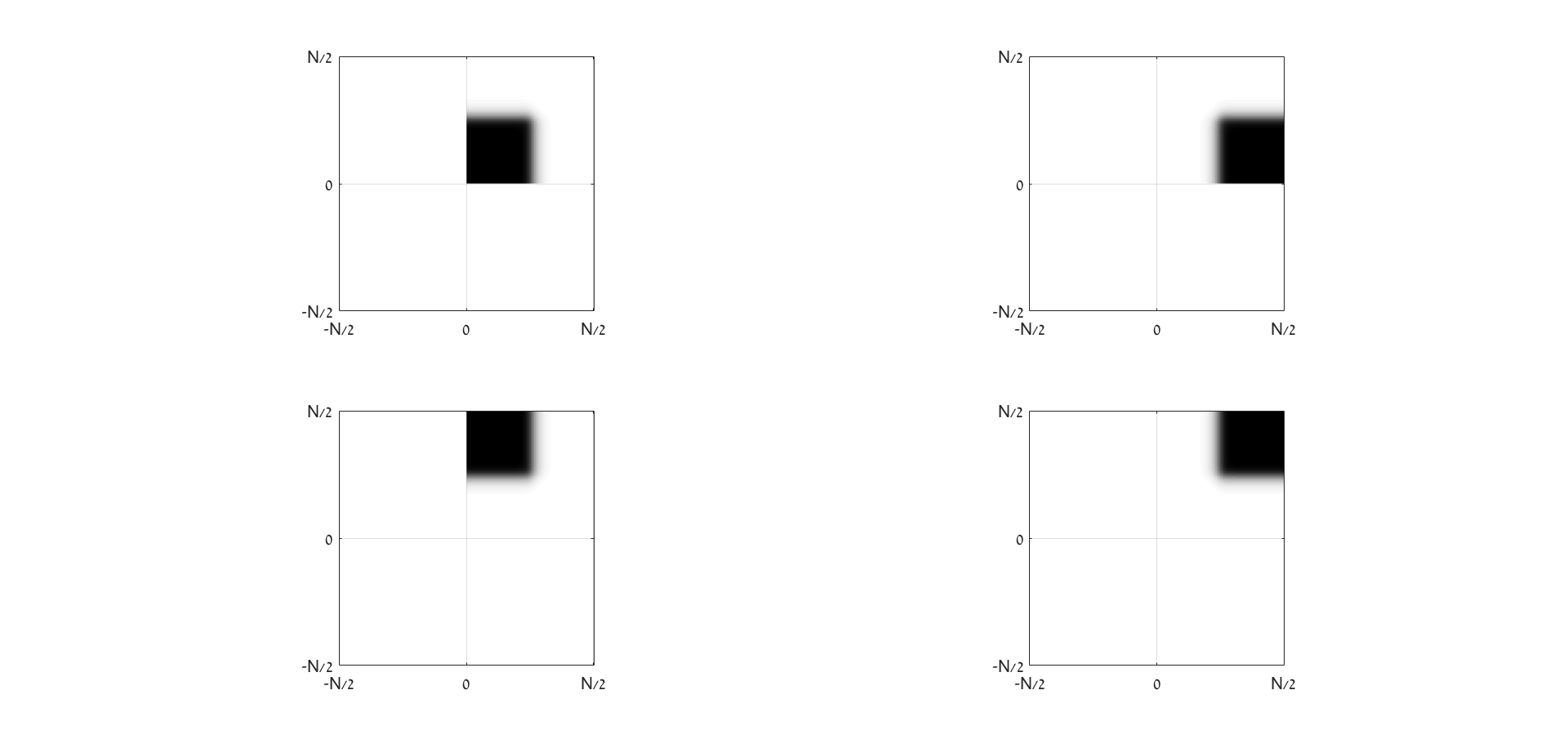}
\hfil
\includegraphics[width=3.0in]{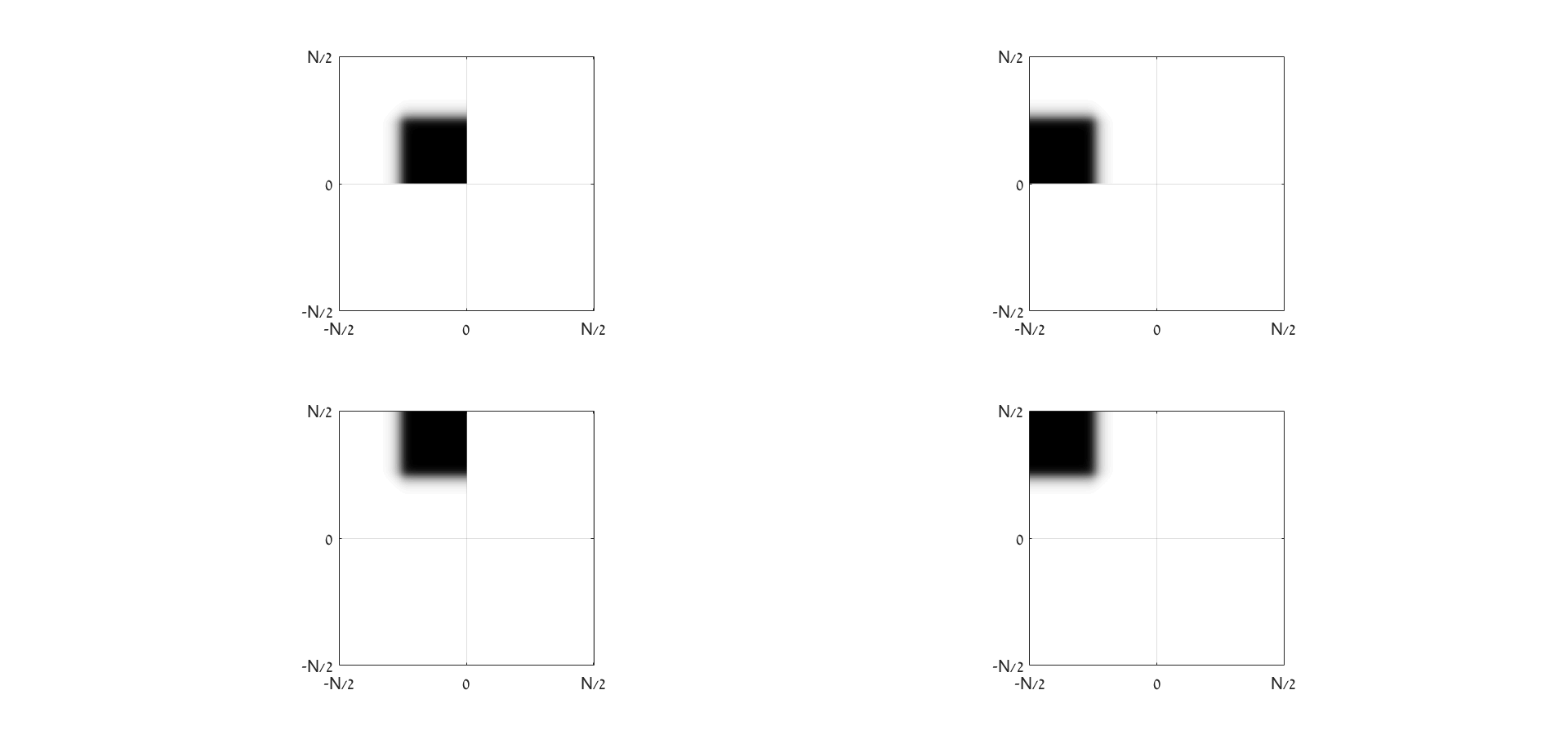}%
\caption{Magnitude \sp a of qWPs $\Psi_{++[1],j ,l}^{10}$  (left)  and $\Psi_{+-[1],j ,l}^{10}$ (right) from the first \d\ level }
\label{fpm_1_2d}
\end{figure}

The \sp a of 1D qWPs $\left\{\Psi_{+[m],j }^{2r}\right\},\;j =0,...,2^{m}-1$, fill the non-negative half-band $[0,N/2]$, and vice versa for the qWPs $\left\{\Psi_{-[m],j }^{2r}\right\},\;j =0,...,2^{m}-1$ . Therefore,  the \sp a of 2D qWPs $\left\{\Psi_{++[m],j ,l}^{2r}\right\},\;j ,l=0,...,2^{m}-1$ fill the quadrant  $[0,N/2-1]\times[0,N/2-1]$  of the \ff\ domain, while  the \sp a of 2D qWPs $\left\{\Psi_{+-[m],j ,l}^{2r}\right\}$ fill the quadrant  $[0,N/2-1]\times[-N/2,-1]$. It is clearly seen in Fig. \ref{fpm_1_2d}.
%\begin{SCfigure}%[H]
%\centering
%\caption{Magnitude \sp a of WPs $\Psi_{++[1],j ,l}^{10}$  (left)  and $\Psi_{+-[1],j ,l}^{10}$ (right) from the first \d\ level }
%\includegraphics[width=2.0in]{png/fp_1_2d.png}
%\hfil
%\includegraphics[width=2.0in]{png/fpm_1_2d.png}%
%\label{fpm_1_2d}
%\end{SCfigure}

Consequently, the \sp a of the real-valued 2D WPs $\left\{\vt_{+[m],j ,l}^{2r}\right\},\;j ,l=0,...,2^{m}-1$, and  $\left\{\vt_{-[m],j ,l}^{2r}\right\}$ fill the pairs  of quadrant  $\mathbf{Q}_{+}\srr[0,N/2-1]\times[0,N/2-1]\bigcup[-N/2,-1]\times[-N/2,-1]$ and $\mathbf{Q}_{-}\srr[0,N/2-1]\times[-N/2,-1]\bigcup[-N/2,-1]\times[0,N/2-1]$, respectively (Figs. \ref{fpp_2_2d} and \ref{fpm_2_2d}).

By this reason, none linear combination of the WPs $\left\{\vt_{+[m],j ,l}^{2r}\right\}$  and their shifts can serve as a basis in the \ss\ space $\Pi[N,N]$. The same is   true for WPs $\left\{\vt_{-[m],j ,l}^{2r}\right\}$. However,  combinations of  the WPs $\left\{\vt_{\pm[m],j ,l}^{2r}\right\}$ provide frames of the  space $\Pi[N,N]$.
\subsection{One-level 2D \t s}\label{sec:ss61}
The one-level 2D qWP \t s of a \ss\ $\mathbf{X}=\left\{X[k,n] \right\}\in\Pi[N,N]$ are implemented by a  tensor-product scheme mentioned in Section \ref{sec:ss25}.

\subsubsection{Direct transforms with qWPs $\Psi^{2r}_{+\pm[1]}$}\label{sec:sss611}
  Denote by $\tilde{\mathbf{T}}_{\pm}^{h}$ the 1D \t s of row \ss s from $\Pi[N]$ with the \aa\ \mv ces $\tilde{\mathbf{M}}_{\pm}^{q}$ which are defined in \eh{aa_modma10p}. Application of these \t s  to rows of a \ss\  \textbf{X} produces the \c\ arrays

\begin{eqnarray*}\label{tTh+x}
        \tilde{\mathbf{T}}_{+}^{h}\mathbf{\cdot}\mathbf{X} &=&  \left(\za_{+}^{0},\za_{+}^{1}\right),\quad \za_{+}^{j}[k,n]=\eta^{j}[k,n]-i\,\xi^{j}[k,n],
                                                             \\\nn \tilde{\mathbf{T}}_{-}^{h}\mathbf{\cdot}\mathbf{X} &=&  \left(\za_{-}^{0},\za_{-}^{1}\right),\quad \za_{-}^{j}[k,n]=\eta^{j}[k,n]+i\,\xi^{j}[k,n]=(\za_{+}^{j}[k,n])^{*},
                                                             \\\nn
        {\eta}^{j}[k,n] &=&\left\langle \mathbf{X}[k,\cdot],{\psi}^{2r}_{[1],j}[\cdot -2n]\right\rangle,\quad  {\xi}^{j}[k,n]=\left\langle \mathbf{X}[k,\cdot],{\f}^{2r}_{[1],j}[\cdot -2n]\right\rangle,\;j=0,1.
      \end{eqnarray*}

 Denote by ${\mathbf{T}}_{\pm}^{h}$ the 1D inverse \t s with the \sa\ \mv ces ${\mathbf{M}}_{\pm}^{q}$. Due to Proposition \ref{pro:Mq_z}, application of these \t s  to rows of the \c\ arrays  $\za_{\pm}=\left(\za_{\pm}^{0},\za_{\pm}^{1}\right)$, respectively, produces the  2D  \az\ \ss s:
 \begin{equation}\label{Th_za}
  {\mathbf{T}}_{\pm}^{h}\mathbf{\cdot}(\za_{\pm}^{0},\za_{\pm}^{1})=2\bar{\mathbf{X}}_{\pm}=2(\mathbf{X}\pm i\,H(\mathbf{X})),
 \end{equation}
 where H(\textbf{X}) is the 2D \ss\ consisting of the HTs of rows of the \ss\ \textbf{X}.

Denote by $\tilde{\mathbf{T}}_{+}^{v}$     the direct  1D \t\ determined by the \mv x $\tilde{\mathbf{M}}_{+}^{q}$  applicable to columns of the corresponding \ss s. The next step of  the tensor product \t\ consists of  the application of the  1D \t\ $\tilde{\mathbf{T}}_{+}^{v}$ to  columns of the arrays  ${\za}^{j},\;j=0,1.$  As a result, we get   four \t\ \c s arrays:

\begin{eqnarray*}\label{tTv_za0}
        \tilde{\mathbf{T}}_{+}^{v}\mathbf{\cdot}\za_{+}^{l} &=& \tilde{\mathbf{T}}_{+}^{v}\mathbf{\cdot}\left(\eta^{l}-i\,\xi^{l} \right)=\left\{
                                                             \begin{array}{l}
                                                             \left( \a^{0,l}-i\,\b^{0,l} \right)-i\left( \g^{0,l}-i\,\da^{0,l} \right) \\
                                                              \left( \a^{1,l}-i\,\b^{1,l} \right)-i\left( \g^{1,l}-i\,\da^{1,l} \right)
                                                             \end{array}
                                                           \right.
         \\\nn&=&\left\{
                   \begin{array}{l}
                      \mathbf{Z}_{+[1]}^{0,l}= \mathbf{Y}_{+[1]}^{0,l}-i\, \mathbf{C}_{+[1]}^{0,l},\quad \mathbf{Y}_{+[1]}^{0,l}=\a^{0,l}-\da^{0,l},\quad \mathbf{C}_{+[1]}^{0,l}=\b^{0,l}+\g^{0,l} \\
                    \mathbf{Z}_{+[1]}^{1,l}= \mathbf{Y}_{+[1]}^{1,l}-i\, \mathbf{C}_{+[1]}^{1,l},\quad \mathbf{Y}_{+[1]}^{1,l}=\a^{1,l}-\da^{1,l},\quad \mathbf{C}_{+[1]}^{1,l}=\b^{1,l}+\g^{1,l}
                   \end{array}
                 \right.
        ,\\\nn
        \a^{j,l}[k,n]&=&\sum_{\la,\mu=0}^{N-1}X[\la,\mu],{\psi}^{2r}_{[1],j}[\la -2k]{\psi}^{2r}_{[1],l}[\mu -2n],\\\nn
        \da^{j,l}[k,n]&=&\sum_{\la,\mu=0}^{N-1}X[\la,\mu],{\f}^{2r}_{[1],j}[\la -2k]{\f}^{2r}_{[1],l}[\mu -2n],\\\nn
        \b^{j,l}[k,n]&=&\sum_{\la,\mu=0}^{N-1}X[\la,\mu],{\psi}^{2r}_{[1],j}[\la -2k]{\f}^{2r}_{[1],l}[\mu -2n],\\\nn
        \g^{j,l}[k,n]&=&\sum_{\la,\mu=0}^{N-1}X[\la,\mu],{\f}^{2r}_{[1],j}[\la -2k]{\psi}^{2r}_{[1],l}[\mu -2n],\quad j,l=0,1.
      \end{eqnarray*}
Hence, it follows that
\begin{equation}
\label{Zjl+}
\begin{array}{lll}
  Y_{+[1]}^{j,l}[k,n]&=&\sum_{\la,\mu=0}^{N-1}X[\la,\mu]\,{\vt}^{2r}_{+[1],j,l}[\la -2k,\mu -2n],\\ C^{j,l}_{+[1]}[k,n]&=&\sum_{\la,\mu=0}^{N-1}X[\la,\mu]\,{\th}^{2r}_{+[1],j,l}[\la -2k,\mu -2n],\\
  Z_{+[1]}^{j,l}[k,n]&=&\sum_{\la,\mu=0}^{N-1}X[\la,\mu]\,{\Psi}^{2r}_{++[1],j,l}[\la -2k,\mu -2n],\quad j,l=0,1. %C^{j,l}_{+[1]}[k,n]&=&\sum_{\la,\mu=0}^{N-1}X[\la,\mu],{\th}^{2r}_{+[1],j,l}[\la -2k,\mu -2n].
\end{array}
\end{equation}
\br\label{rem:+rec}Recall that the DFT \sp a of WPs ${\vt}^{2r}_{+[1],j,l}$ and ${\th}^{2r}_{+[1],j,l},\;j,l=0,1,$ which are the real and imaginary parts of the qWP ${\Psi}^{2r}_{++[1],j,l}$,  are confined within the area $\mathbf{Q}_{+}$ of the \ff\ domain. It is seen from Eq. \rf{Zjl+}  that  if at least a part of the \sp um of a \ss\ $\mathbf{X}\in\Pi[N,N]$ is located in the area $\mathbf{Q}_{-}$, then the \ss\ $\mathbf{X}\in\Pi[N,N]$ cannot be fully restored  from the \t\ \c s $Z_{+[1]}^{j,l}[k,n]$,  although their number   is the same as the number of samples in the \ss\ $\mathbf{X}$. To achieve a perfect \r, the \c s from the arrays $\mathbf{Z}_{-[1]}^{j,l}$ should be incorporated.\er
The \c\ arrays $\mathbf{Z}_{-[1]}^{j,l}$ are derived in the same way as the arrays $\mathbf{Z}_{+[1]}^{j,l}$. The only \df ce is that, for the 1D \t\ $\tilde{\mathbf{T}}^{h}_{-}$  the  \mv x
 $\tilde{\mathbf{M}}_{-}^{q}$ is used  instead of $\tilde{\mathbf{M}}_{+}^{q}$. For the  \t\ $\tilde{\mathbf{T}}^{v}_{-}$,   the  \mv x
 $\tilde{\mathbf{M}}_{+}^{q}$ is used. Consequently, to derive the \c\ arrays $\mathbf{Z}_{-[1]}^{j,l}$,  the \t\ $\tilde{\mathbf{T}}_{+}^{v}$ should be applied to columns of the arrays $\za^{l}_{-}=(\za^{l}_{+})^{*}$. As a result, we get
 \begin{eqnarray*}
 \label{Zjl-}
   \mathbf{Z}_{-[1]}^{j,l}= \mathbf{Y}_{-[1]}^{j,l}+i\,\mathbf{C}_{-[1]}^{j,l} =\sum_{\la,\mu=0}^{N-1}X[\la,\mu],{\Psi}^{2r}_{+-[1],j,l}[\la -2k,\mu -2n],\; j,l=0,1.
 \end{eqnarray*}

\subsubsection{Inverse transforms with qWPs $\Psi^{2r}_{+\pm[1]}$}\label{sec:sss612}
Denote by $\mathbf{T}_{+}^{v}$ the 1D inverse \t\ with the \sa\ \mv x ${\mathbf{M}}_{+}^{q}$ applicable to columns of the \c\ arrays. Denote by $H(\za_{\pm}^{l}),\;l=0,1,$
the HTs of the arrays consisting of columns of the \c\ arrays  $\za_{\pm}^{j}$.
  Proposition \ref{pro:Mq_z} implies that
\begin{eqnarray*}
  \mathbf{T}_{+}^{v}\mathbf{\cdot}  \left(
                                         \begin{array}{c}
                                           \mathbf{Z}_{+[1]}^{0,l} \\
                                           \mathbf{Z}_{+[1]}^{1,l} \\
                                         \end{array}
                                       \right)=2\bar{\za}_{+}^{l},\quad \mathbf{T}_{+}^{v}\mathbf{\cdot}  \left(
                                         \begin{array}{c}
                                           \mathbf{Z}_{-[1]}^{0,l} \\
                                           \mathbf{Z}_{-[1]}^{1,l} \\
                                         \end{array}
                                       \right)=2\bar{\za}_{-}^{l},%=(\za_{+}^{l})^{*}+i\,(H(\za_{+}^{l}))^{*},
\end{eqnarray*}
where $\;l=0,1$ and $\bar{\za}_{\pm}^{l}=\za_{\pm}^{l}+i\,H(\za_{\pm}^{l})$ are \az\ \c\ arrays. Denote by $\mathbf{G}$ a \ss\ from $\Pi[N,N]$ such that
\begin{equation}\label{sig_G}
 \tilde{\mathbf{T}}^{h}_{\pm}\mathbf{\cdot}\mathbf{G}=\left(H(\za_{\pm}^{0}), H(\za_{\pm}^{1})\right)\Longrightarrow
{\mathbf{T}}_{\pm}^{h}\mathbf{\cdot}\left(H(\za_{\pm}^{0}), H(\za_{\pm}^{1})\right)=4(\mathbf{G}\pm i\,H(\mathbf{G}).
\end{equation}

Equations \rf{Th_za} and \rf{sig_G} imply that the applications of the \t s $ {\mathbf{T}}_{\pm}^{h}$ to rows of the respective \c\ arrays  results in the following relations:
 \begin{eqnarray}\label{Th_zahza}
  \mathbf{X}_{+}&\srr&{\mathbf{T}}_{+}^{h}\mathbf{\cdot}\left(\bar{\za}_{+}^{0}),\bar{\za}_{+}^{1}\right)=
4\left(\mathbf{X}+i\,H(\mathbf{X})+i\mathbf{G}-H(\mathbf{G}))\right),\\\nn
\mathbf{X}_{-}&\srr&{\mathbf{T}}_{-}^{h}\mathbf{\cdot}\left(\bar{\za}_{-}^{0}),\bar{\za}_{-}^{1}\right)
=4\left(\mathbf{X}-i\,H(\mathbf{X})+i\mathbf{G}+H(\mathbf{G}))\right).
 \end{eqnarray}

Finally, we have the \ss\ \textbf{X} restored by $\mathbf{X}=\mathfrak{Re}(\mathbf{X}_{+}+\mathbf{X}_{-})/8$.
%%%%%%%%%%%%%%%%%%%%%%%

\begin{figure}[H]
\centering
\includegraphics[width=5.2in]{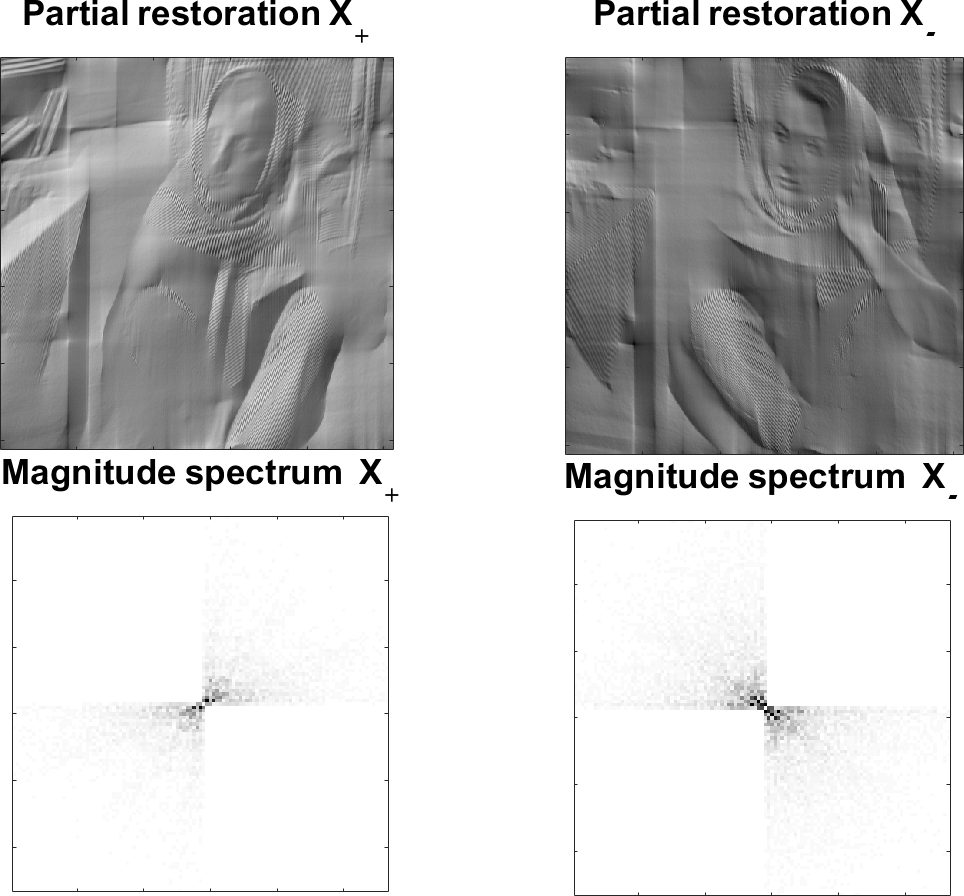}
\caption{Top left: Partially restored ``Barbara" image by $\mathfrak{Re}(\mathbf{X}_{+})$. Top right: Partially restored  image by $\mathfrak{Re}(\mathbf{X}_{-})$.
Bottom left: Magnitude DFT \sp um of $\mathbf{X}_{+}$. Bottom right: Magnitude DFT \sp um of $\mathbf{X}_{-}$}
\label{xp_xm_x5}
\end{figure}

\begin{figure}[H]
\centering
\includegraphics[width=4.0in]{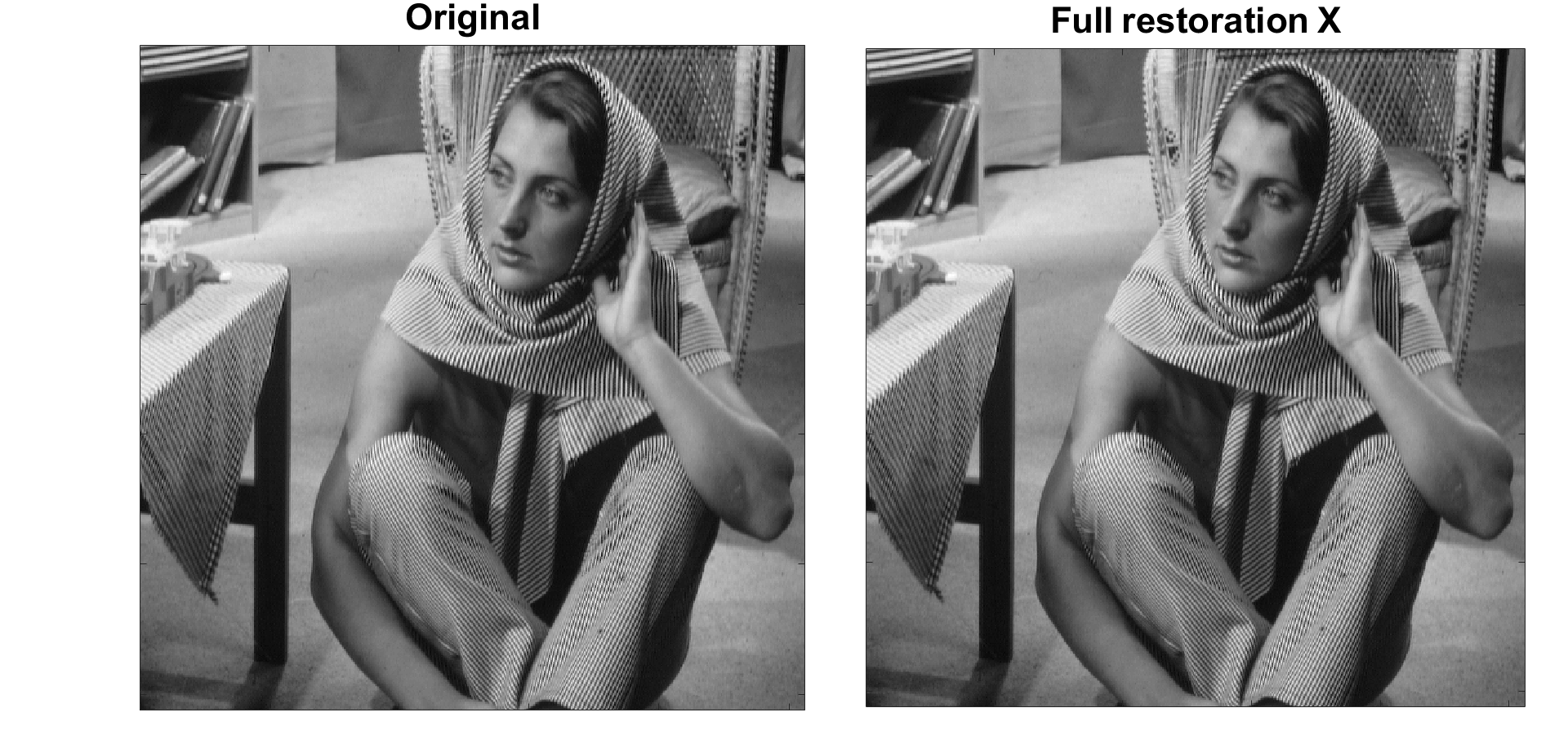}
\caption{Left: Original ``Barbara" image. Right: Fully restored  image by $\mathfrak{Re}(\mathbf{X}_{+}+\mathbf{X}_{-})/8$}
\label{xp_xm_x}
\end{figure}

Figures \ref{xp_xm_x5} and \ref{xp_xm_x} illustrate the  image ``Barbara" restoration by the 2D \ss s $\mathfrak{Re}(\mathbf{X}_{\pm})$ and $\mathbf{X}=\mathfrak{Re}(\mathbf{X}_{+}+\mathbf{X}_{-})/8$. The \ss\ $\mathfrak{Re}(\mathbf{X}_{+})$ captures edges oriented to \emph{north-east}, while $\mathfrak{Re}(\mathbf{X}_{-})$ captures edges oriented to \emph{north-west}. The \ss\ $\mathbf{X}$ perfectly restores the image achieving  PSNR=313.8596 dB.

\subsection{Multi-level 2D \t s}\label{sec:ss62}
It was established in Section \ref{sec:ss42} that the 1D qWP \t s of a \ss\ $\mathbf{x}\in\Pi[N]$  to the second and further  \d\ levels are implemented by the  iterated  application of
the \fb s,  that  are determined by their \aa\ \mv ces $\tilde{\mathbf{M}}[2^{m}n],\;m=1,...,M-1,$ to the \c s arrays $\mathbf{z}_{\pm[m]}^{\la}$. The \t s  applied to the arrays $\mathbf{z}_{\pm[m]}^{\la}$  produce the arrays $\mathbf{z}_{\pm[m+1]}^{\rr}$, respectively. The inverse \t\ consists of  the iterated  application of the
\fb s that are determined by their \sa\ \mv ces ${\mathbf{M}}[2^{m}n],\;m=1,...,M-1,$ to the \c s arrays $\mathbf{z}_{\pm[m+1]}^{\rr}$.  In that way the first-level \c\ arrays $\mathbf{z}_{\pm[1]}^{\la},\;\la=0,1$ are restored.

The tensor-product of the 2D \t s of a \ss\ $\mathbf{X}\in\Pi[N,N]$ consists of the subsequent application of the 1D \t s to columns and rows of the \ss\ and \c s arrays. By application of \fb s, which are determined by the \aa\ \mv x $\tilde{\mathbf{M}}[2n]$ to columns and rows of  \c s array $\mathbf{Z}_{\pm[1]}^{j,l}$, we derive four second-level arrays
$\mathbf{Z}_{\pm[2]}^{\rr,\tau},\;\rr=2j,2j+1;\;\tau=2l,2l+1$. The arrays $\mathbf{Z}_{\pm[1]}^{j,l}$ are restored by the application of the
\fb s that are determined by their \sa\ \mv ces ${\mathbf{M}}[2n] $ to rows and columns of the \c s arrays $\mathbf{Z}_{\pm[2]}^{\rr,\tau},\;\rr=2j,2j+1;\;\tau=2l,2l+1$. The transition from the second to further levels and back are executed similarly using the \mv ces $\tilde{\mathbf{M}}[2^{m}n]$ and ${\mathbf{M}}[2^{m}n]$, respectively. The inverse \t s produce the \c s arrays $\mathbf{Z}_{\pm[1]}^{j,l},\;j,l=0,1,$ from which the \ss\ $\mathbf{X}\in\Pi[N,N]$ is restored using the \sa\ \mv ces ${\mathbf{M}}_{\pm}^{q}[n] $ as it is explained in Section \ref{sec:sss612}.

All the computations are implemented in the \ff\ domain using the FFT.  For example, the Matlab execution of the  2D qWP \t\ of a $512\times 512$ image down to the sixth \d\ level  takes 1.34 seconds. The four-level \t\ takes 0.28 second.
%\br\label{rem:dob_tre}
\paragraph{Summary}
The 2D qWP \pr\ of a \ss\ $\mathbf{X}\in\Pi[N,N]$ is implemented by a dual-tree scheme. The first step produces two sets of the \c s arrays: $\mathbf{Z}_{+[1]}=\left\{\mathbf{Z}_{+[1]}^{j,l}\right\}, \;j,l,=0,1,$ which are derived using the \aa\ \mv x $\tilde{\mathbf{M}}_{+}^{q}[n]$, and $\mathbf{Z}_{-[1]}=\left\{\mathbf{Z}_{-[1]}^{j,l}\right\}, \;j,l,=0,1,$ which are derived using the \aa\ \mv x $\tilde{\mathbf{M}}_{-}^{q}[n]$.  Further \d\ steps are implemented in parallel on the sets $\mathbf{Z}_{+[1]}$ and  $\mathbf{Z}_{-[1]}$  using the same \aa\ \mv x $\tilde{\mathbf{M}}[2^{m}n]$, thus producing two multi-level sets of the \c s arrays $\left\{\mathbf{Z}_{+[m]}^{j,l}\right\}$ and $\left\{\mathbf{Z}_{-[m]}^{j,l}\right\},\;m=2,...,M,\;j,l=0,2^{m}-1$.

By parallel implementation of the inverse \t s on the \c s from the   sets  $\left\{\mathbf{Z}_{+[m]}^{j,l}\right\}$ and $\left\{\mathbf{Z}_{-[m]}^{j,l}\right\}$   using the same \sa\ \mv x ${\mathbf{M}}[2^{m}n]$, the sets  $\mathbf{Z}_{+[1]}$ and  $\mathbf{Z}_{-[1]}$ are restored, which, in turn, provide the \ss s   $\mathbf{X}_{+}$ and  $\mathbf{X}_{-}$, using the \sa\ \mv ces  ${\mathbf{M}}_{+}^{q}[n]$ and ${\mathbf{M}}_{-}^{q}[n]$, respectively. Typical  \ss s   $\mathbf{X}_{\pm}$ and their DFT \sp a are displayed in Fig. \ref{xp_xm_x5}.

Prior to the \r, some structures, possibly \df t, are defined in the sets  $\left\{\mathbf{Z}_{+[m]}^{j,l}\right\}$ and $\left\{\mathbf{Z}_{-[m]}^{j,l}\right\},\;m=1,...M,$ (2D \ww\  or Best Basis structures,  for example) and some manipulations on the \c s, (thresholding, $l_1$ minimization, for example) are executed.
\section{Numerical examples}\label{sec:s7}
In this section, we present two examples of application of the 2D qWPs to image restoration. These examples illustrate the ability of  the qWPs to restore edges and texture details even from severely damaged images.  Certainly, this ability stems from the fact that the designed 2D qWP \t s provide a variety of 2D \we s oriented in multiple directions, from perfect \ff\ resolution of these \we s and, last but not least, from oscillatory structure of many \we s.
\subsection{Denoising examples}\label{sec:ss71}An image $\mathbf{I}$ is \ry ed by the 2D \ss\ $\mathbf{X}\in\Pi[N,N]$. The image, which is corrupted by additive Gaussian noise with STD=$\o$,  is \ry ed by the 2D \ss\ $\mathbf{X}_{\o}$. We apply the following image denoising scheme:
\begin{itemize}
  \item 2D  \t\ of the  \ss\ $\mathbf{X}$ with directional qWPs $\Psi^{2r}_{++[m]}$ and $\Psi^{2r}_{+-[m]}$ down to level $M$ is implemented to generate two sets of the \c s arrays  $\left\{\mathbf{Z}_{+[m]}^{j,l}\right\}$ and $\left\{\mathbf{Z}_{-[m]}^{j,l}\right\},\;m=1,...M,\;j,l=0,...,2^{m}-1$.
 \item For comparison, the 2D WP \t\ of the  \ss\ $\mathbf{X}$ with non-directional WP $\psi^{2r}_{[m]}$ down to level $M$ is implemented thus generating the set  $\left\{\mathbf{y}_{[m]}^{j,l}\right\},\;m=1,...M,\;j,l=0,...,2^{m}-1,$ of the \c s arrays.
  \item In each of the sets $\left\{\mathbf{Z}_{\pm[m]}^{j,l}\right\}$ and $\left\{\mathbf{y}_{[m]}^{j,l}\right\}$ the ``Best Basis" is selected by a standard procedure of comparison the cost function (entropy or $l_1$ norm) of the ``parent" \c s block with the cost functions of its ``offsprings" (see \cite{coiw1}). The selected bases are designated by $\mathbf{B}_{\pm[M]}$ and $\mathbf{b}_{[M]}$, respectively.
  \item The number of the \t\ \c s $\mathbf{Z}_{\pm[B]}$ and  $\mathbf{y}_{[B]}$  associated with each basis is the same as the number $N^{2}$ of pixels in the image.
\item Denoising of the \ss\ $\mathbf{X}_{\o}$ is implemented by hard thresholding of the \c s $\mathbf{Z}_{\pm[B]}$ and  $\mathbf{y}_{[B]}$.
\item The thresholds for the \c s arrays are defined by the following naive scheme:
\begin{enumerate}
  \item The absolute values of the coefficients from each set  $\mathbf{Z}_{\pm[B]}$ and  $\mathbf{y}_{[B]}$ are arranged in an ascending order thus  forming the  \sq s $\mathbf{A}_{\pm}$ and  $\mathbf{A}$, respectively.
  \item The values of the $L$-th term $T_{\pm}={A}_{\pm}[L]$ and $T={A}[L]$ are selected as the thresholds for the \c s arrays $\mathbf{Z}_{\pm[B]}$ and  $\mathbf{y}_{[B]}$, respectively. The number $L$ is chosen  depending on the noise intensity.
\end{enumerate}
\item The sets of the \t\ \c s $\mathbf{Z}_{\pm[B]}$ and  $\mathbf{y}_{[B]}$ are subjected to thresholding  with the selected values $T_{\pm}$ and $T$, respectively, and, after that the inverse \t s are applied to produce the 2D \ss s $\mathbf{X}_{d}=\mathfrak{Re}(\mathbf{X}_{+}+\mathbf{X}_{-})/8$ from the directional qWPT and $\mathbf{X}_{nd}$ from the non-directional tensor product WPT.
\end{itemize}
\subsubsection{Example I: ``Pentagon" image}\label{sec:sss711}The ``Pentagon" image of size $1024\times 1024$ (1048576 pixels) was corrupted by additive Gaussian noise with STD=30 dB. As a result, the PSNR of  the corrupted image was 18.59 dB. The corrupted image was decomposed by the directional qWPs $\Psi^{4}_{++[m]}$ and $\Psi^{4}_{+-[m]}$ originating from the fourth-order \ds s down to a fourth \d\ level. In this way, two sets $\left\{\mathbf{Z}_{+[m]}^{j,l}\right\}$ and $\left\{\mathbf{Z}_{-[m]}^{j,l}\right\},\;m=1,...4,\;j,l=0,...,2^{m}-1$, of the \t\ \c s were produced. For comparison, the corrupted image was decomposed by the non-directional WPs $\psi^{4}_{[m]}$, which resulted in the set  $\left\{\mathbf{y}_{[m]}^{j,l}\right\},\;m=1,...4,\;j,l=0,...,2^{m}-1,$ of the \c s arrays.

The ``Best Bases" $\mathbf{B}_{\pm[4]}$ and $\mathbf{b}_{[4]}$ were designed for the \c\ arrays. The thresholds $T_{\pm}={A}_{\pm}[L]$ and $T={A}[L]$ were selected for each set of the \c\ arrays. $L=1010000$ was chosen. Thus, we had $T_{+}=130.42$, $T_{-}=126.8$ and $T=73.30$. 
\br\label{rem:2thre}The threshold $T=73.30$ is \as ely two times less than the threshold $T_{\pm}$. It happens because the \t\ \c s  $\left\{\mathbf{Z}_{+[m]}^{j,l}\right\}$ are complex-valued and the threshold is operating with absolute values of the \c s.\er
\br\label{rem_penta}qWPs oriented in 62 \df t directions were involved. The Matlab implementation of all the above procedures takes 4.2 seconds.\er

\begin{figure}[H]
\resizebox{17cm}{9cm}{
\includegraphics{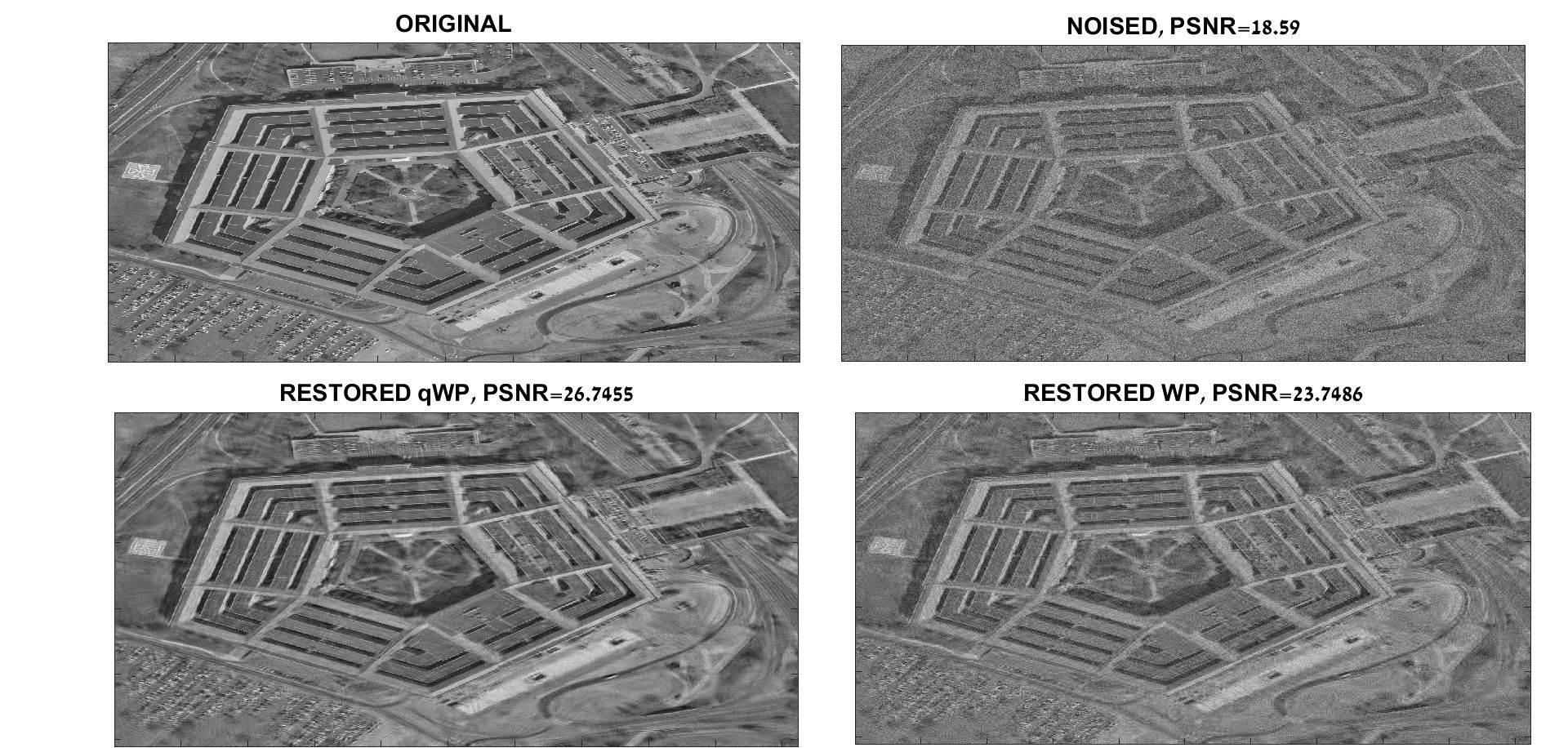}
}
\caption{Top left: Original ``Pentagon" image. Top right: Image corrupted by noise with STD=30 dB. Bottom  left: The  qWPT-based restored  image  $\mathbf{X}_{d}$ . Bottom right: WPT-based restored image  $\mathbf{X}_{nd}$}
\label{penta30}
\end{figure}

Figure \ref{penta30} is  the outputs of the image \r\ by the  directional qWPT  and  the non-directional tensor product WPT from the thresholded \c s arrays. The qWPT-restored image  image  $\mathbf{X}_{d}$ has PSNR=26.75 dB versus PSNR=23.75 dB for the WPT-restored image  $\mathbf{X}_{nd}$. Visually, image  $\mathbf{X}_{d}$ is cleaner in comparison to $\mathbf{X}_{nd}$ and more fine details are restored.

\subsubsection{Example II: ``Barbara" image}\label{sec:sss712}
We present two cases with the ``Barbara" image. In one case, the image was corrupted by an additive Gaussian noise with STD=30 dB and in the other, the noise was more intensive with STD=50 dB.
\begin{description}
  \item[Noise  with STD=30 dB:] In this case, the PSNR of the corrupted image was 18.59 dB.
In order to avoid boundary effects, the image of size $512\times 512$ was \sy ally extended to the size $1024\times 1024$. After \pr, the results were shrunk to the original size.
The corrupted image was decomposed by the directional qWPs $\Psi^{4}_{++[m]}$ and $\Psi^{4}_{+-[m]}$ originating from the fourth-order \ds s down to third \d\ level. In this way, the two sets $\left\{\mathbf{Z}_{+[m]}^{j,l}\right\}$ and $\left\{\mathbf{Z}_{-[m]}^{j,l}\right\},\;m=1,...3,\;j,l=0,...,2^{m}-1$, of the \t\ \c s were produced. For comparison, the corrupted image was decomposed by the non-directional WPs $\psi^{4}_{[m]}$, which resulted in the set  $\left\{\mathbf{y}_{[m]}^{j,l}\right\},\;m=1,...3,\;j,l=0,...,2^{m}-1,$ of the \c s arrays.

The ``Best Bases" $\mathbf{B}_{\pm[3]}$ and $\mathbf{b}_{[3]}$ were designed for the \c\ arrays. The thresholds $T_{\pm}={A}_{\pm}[L]$ and $T={A}[L]$ were selected for each set of the \c\ arrays. $L=1011560$ was chosen. Thus, we had $T_{+}=144.51$, $T_{-}=144.66$ and $T=82.83$.

Figures \ref{barb30} and \ref{barb30F} are the outputs of the image \r\ by using the  directional qWPT  and  the non-directional tensor product WPT from the thresholded \c s arrays. The qWPT-restored image  $\mathbf{X}_{d}$ has PSNR=28.52 dB versus PSNR=25.20 dB for the WPT-restored image  $\mathbf{X}_{nd}$. Visually, the image  $\mathbf{X}_{d}$ is much cleaner in comparison to $\mathbf{X}_{nd}$ and almost all edges and the texture structure are restored.

\br\label{rem_barb50}qWPs oriented in 31 \df t directions are involved. In this case, the image under \pr\ has four times more pixels than the original image.  The Matlab implementation of all the  procedures including 3-level qWP and WP \t s, design of ``Best Bases", thresholding of the \c s arrays and the inverse \t s  takes 3.5 seconds. Processing the image without the extension takes 0,88 seconds and produces the restores image with PSNR=28.32 dB.\er
\begin{figure}[H]
\begin{center}
\resizebox{14cm}{9cm}{
\includegraphics{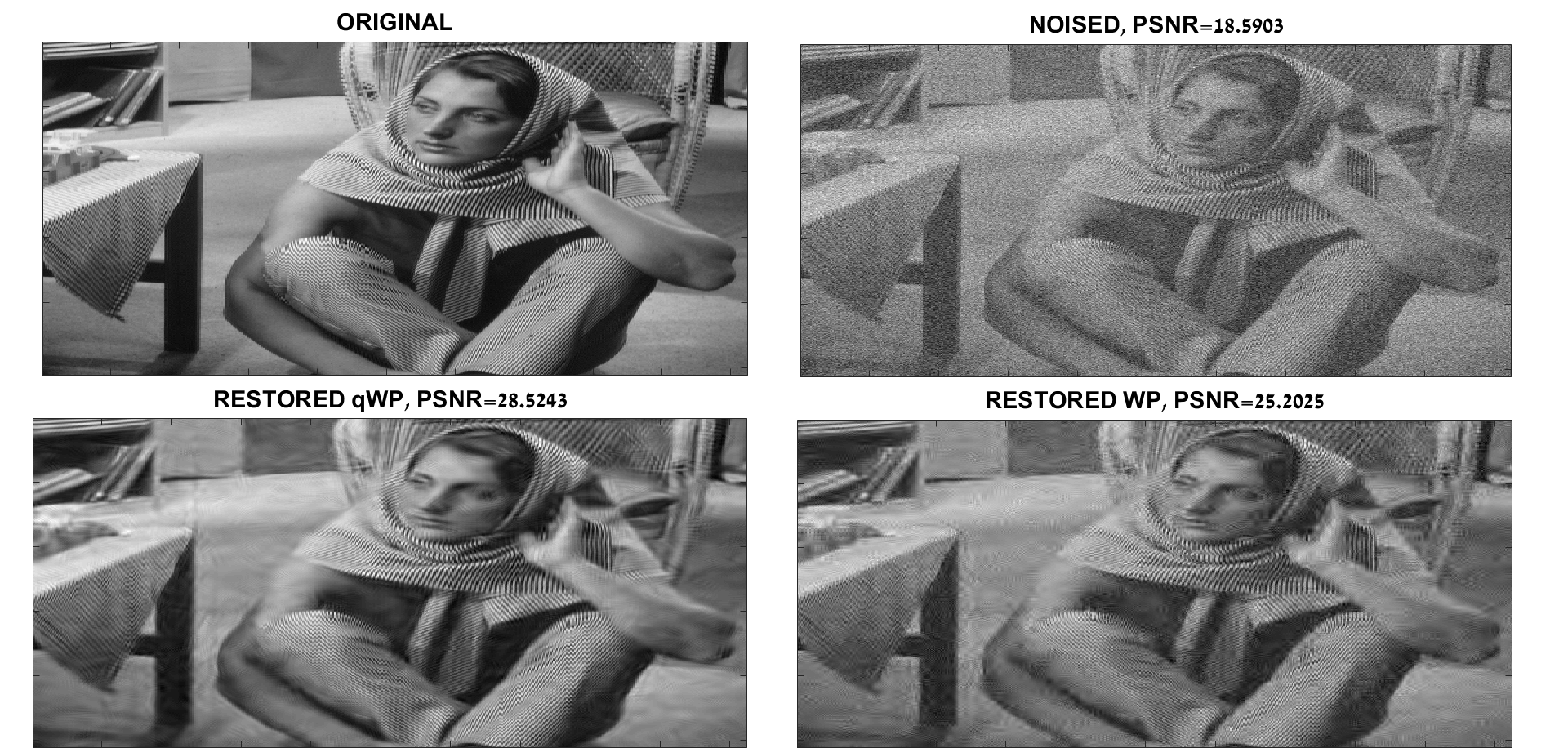}
}
\end{center}
%\centering
%\includegraphics[width=7.0in]{png/barb30.png}
\caption{Top left: Original ``Barbara" image. Top right: Image corrupted by noise with STD=30 dB. Bottom  left: The  qWPT-based restored  image  $\mathbf{X}_{d}$ . Bottom right: WPT-based restored image  $\mathbf{X}_{nd}$}
\label{barb30}
\end{figure}

\begin{figure}[H]
\begin{center}
\resizebox{13cm}{9cm}{
\includegraphics{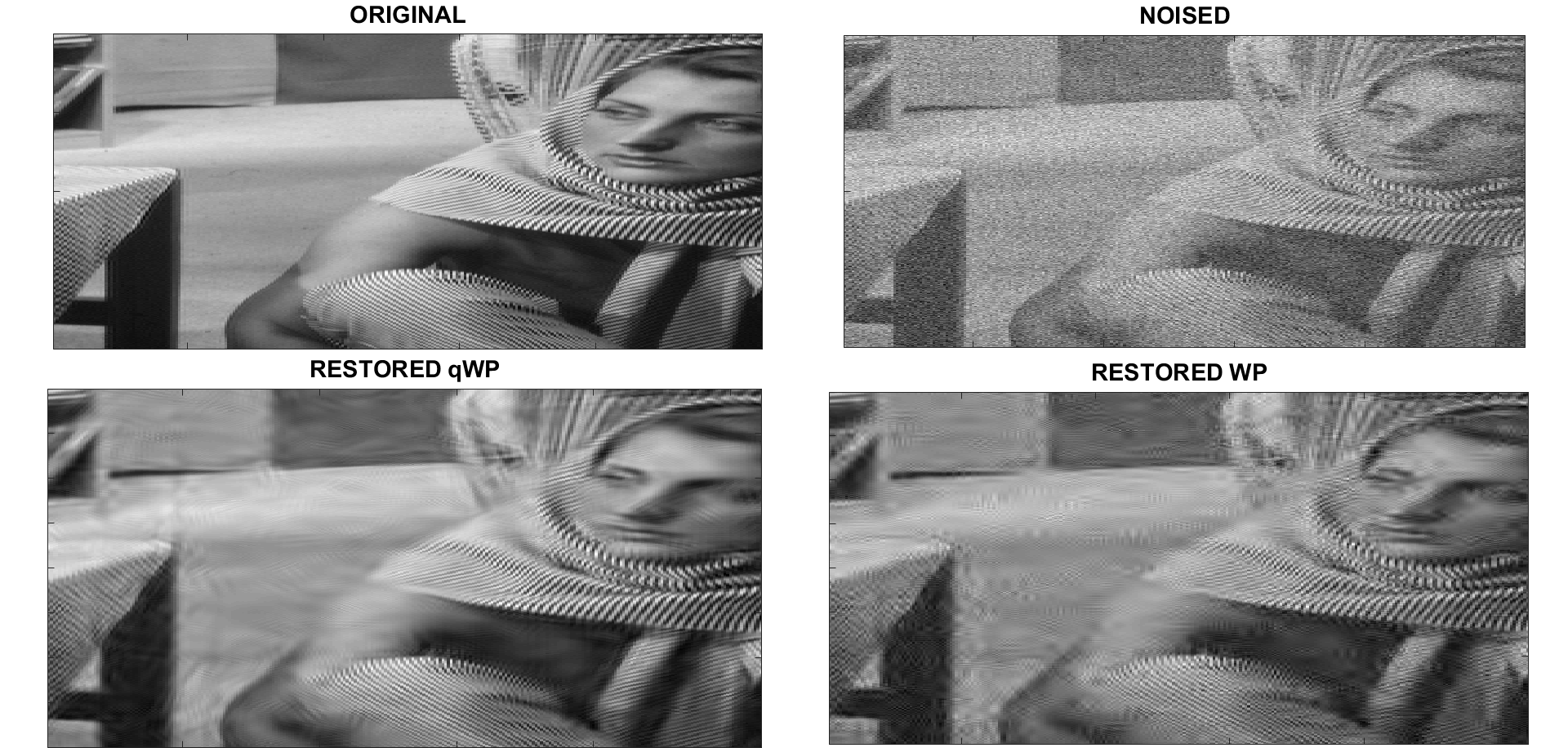}
}
\end{center}
%\centering
%\includegraphics[width=7.0in]{png/barb30F.png}
\caption{Fragments of the images shown in Fig. \ref{barb30}}
\label{barb30F}
\end{figure}

  \item[Noise  with STD=50 dB:] In this case, the PSNR of the corrupted image was 14.14 dB.  The same operations as in the previous case were applied to the corrupted image. $L=1019600$ was chosen. Thus, we had $T_{+}=234.18$, $T_{-}=234.16$ and $T=137.29$.

Figures \ref{barb50} and \ref{barb50F} are the outputs of the image \r\ by using the directional qWPT  and  the non-directional tensor product WPT from the thresholded \c s arrays.  The qWPT-based restored image   $\mathbf{X}_{d}$ has PSNR=25.66 dB versus PSNR=22.20 dB for the WPT-based restored image  $\mathbf{X}_{nd}$. Visually, the image  $\mathbf{X}_{d}$ is cleaner in comparison to $\mathbf{X}_{nd}$, which, in addition comprises many artifacts.  Many edges and the texture structure are retained in  the image   $\mathbf{X}_{d}$.
\begin{figure}[H]
\begin{center}
\resizebox{14cm}{9cm}{
\includegraphics{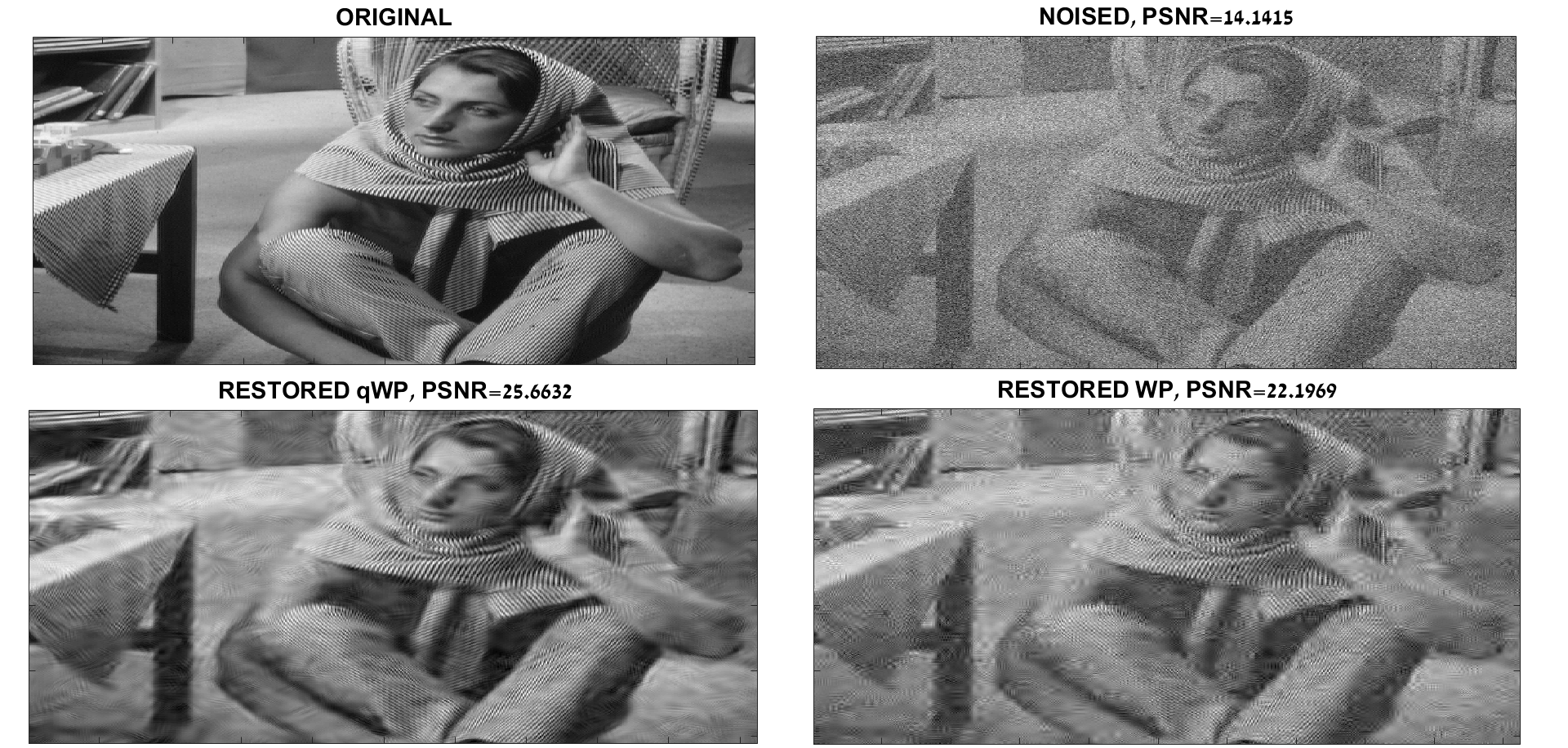}
}
\end{center}
%\centering
%\includegraphics[width=7.0in]{png/barb50.png}
\caption{Top left: Original ``Barbara" image. Top right: Image corrupted by noise with STD=50 dB. Bottom  left: The  qWPT-based restored  image  $\mathbf{X}_{d}$ . Bottom right: WPT-based restored image  $\mathbf{X}_{nd}$}
\label{barb50}
\end{figure}
\begin{figure}[H]
\begin{center}
\resizebox{13cm}{9cm}{
\includegraphics{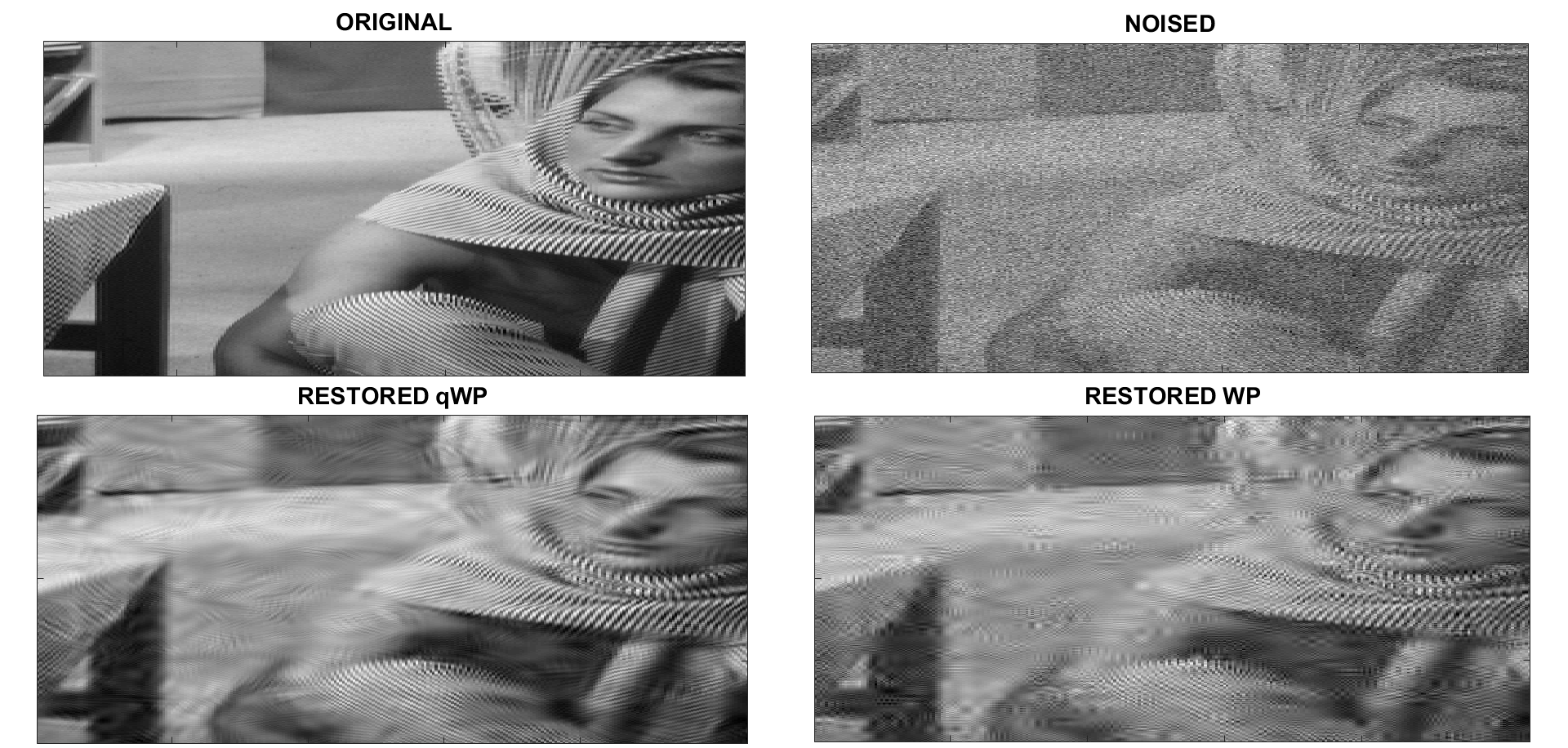}
}
\end{center}
%\centering
%\includegraphics[width=7.0in]{png/barb50F.png}
\caption{Fragments of the images shown in Fig. \ref{barb50}}
\label{barb50F}
\end{figure}

\paragraph{Comment}In the above examples, we did not apply sophisticated adaptive denoising schemes such as Gaussian scale mixture model (\cite{porti_strela}) or bivariate shrinkage (\cite{sen_seles}). Our goal here was to compare performance of the directional qWP \t s with the performance of the standard tensor-product WP \t s.

\end{description}
\subsection{Image restoration  examples}\label{sec:ss72}In this section we present a few cases  of   image restoration using directional qWPs. Images to be restored  were degraded by blurring, aggravated by  random noise and
random loss of significant number of pixels. In our previous work (\cite{Inv_frame} and Chapter 18 in \cite{ANZ_book1}) we developed the  image restoration scheme utilizing  2D  \ww\ frames designed in  Chapter 18 of \cite{ANZ_book1}. In the examples presented below we use, generally, the same scheme as   in \cite{ANZ_book1} with the \df ce that the directional qWPs designed in Section \ref{sec:s6} are used instead of \ww\ frames.
\subsubsection{Brief outline of the restoration scheme}\label{sec:sss721}

Images are restored by the application of the \emph{split
Bregman iteration} (SBI) scheme \cite{gold_os} that uses the
so-called \emph{\aa-based} approach (see for example \cite{ji_shen_xu}).

Denote by $\mathbf{u}=\left\{u[\k,\n]\right\}$ the original image array to be restored from the degraded array
\(
  \mathbf{f}=\mathbf{K}\,\mathbf{u}+\e,
\) where $\mathbf{K}$ denotes the operator of  2D \dd\ \cc\ of the array
$\mathbf{u}$ with a kernel $\mathbf{k}=\left\{k[\k,\n]\right\}$,  and
$\e=\left\{e_{k,n}\right\}$ is the random error array.
$\mathbf{K^*}$ denotes the  conjugate operator of $\mathbf{K}$,
which implements the \dd\ \cc\  with the transposed kernel
$\mathbf{k}^{T}$. If some number of pixels are missing then the image
$\mathbf{u}$  should be restored from the available data
 \begin{equation}\label{av_dat}
  \mathbf{P}_{\Lambda}\,\mathbf{f}=\mathbf{P}_{\Lambda}\,(\mathbf{K}\,\mathbf{u}+\e),
\end{equation}
where $\mathbf{P}_{\La}$ denotes the projection on the
remaining set of pixels.

The solution scheme is based on the assumption that the original
image $\mathbf{u}$  can be
sparsely represented in the qWP domain. Denote by
$\mathbf{\tilde{F}}$ the  operator of qWP  expansion of the image
$\mathbf{u}$. To be specific, the 2D  \t\ of the  \ss\ $\mathbf{X}$ with directional qWP $\Psi^{2r}_{++[m]}$ and $\Psi^{2r}_{+-[m]}$ down to level $M$ is implemented to generate two sets of the \c s arrays  $\left\{\mathbf{Z}_{+[m]}^{j,l}\right\}$ and $\left\{\mathbf{Z}_{-[m]}^{j,l}\right\},\;m=1,...M,\;j,l=0,...,2^{m}-1$.
In each of the sets $\left\{\mathbf{Z}_{\pm[m]}^{j,l}\right\}$  either   ``Best Basis" or ``basis", which  consist of shifts of  all the WPs from the \d\ level $M$, are selected. The bases are designated by $\mathbf{B}_{\pm[M]}$.
 The number of the \t\ \c s $\mathbf{Z}_{\pm[B]}$  associated with each basis is the same as the number $N^{2}$ of pixels in the image. Thus,  $\mathbf{C}\srr\mathbf{\tilde{F}}\,\mathbf{u}=\mathbf{Z}_{+[B]}\bigcup\mathbf{Z}_{-[B]}$  is the set of the  \t\  \c s.

Denote by $\mathbf{F}$ the \r\ operator  of  the image $\mathbf{u}$
from the set of the \t\ \c s. We get
$\mathbf{F}\,\mathbf{C}=\mathbf{u}=\mathfrak{Re}(\mathbf{u}_{+}+\mathbf{u}_{-})/8$,
$\mathbf{F}\,\mathbf{\tilde{F}}=\mathbf{I}$, where $\mathbf{I}$ is
the identity operator.

An  approximate solution to  \eh{av_dat} is derived
via minimization of the \fd al
\begin{equation}\label{av_dat_min}
  \min_{u}\frac{1}{2}\left\|\mathbf{P}_{\La}\,(\mathbf{K}\,\mathbf{u}-f)\right\|_{2}^{2}+\la\,\left\|\mathbf{\tilde{F}}\,\mathbf{u}\right\|_{1},
\end{equation}
where $\left\|\cdot\right\|_{1}$ and $\left\|\cdot\right\|_{2}$ are
the  $l_{1}$ and the $l_{2}$ norms of the \sq s, respectively.  If
$\mathbf{x}=\left\{x[\k,\n]\right\},\;{\k}=0,...,k,\;{\n}=0,...,n$,
then
\[\left\|\mathbf{x}\right\|_{1}\srr\sum_{\k=0}^{k-1}\sum_{\n=0}^{n-1}|x[\k,\n]|,\quad \left\|\mathbf{x}\right\|_{2}\srr\sqrt{\sum_{\k=0}^{k-1}\sum_{\n=0}^{n-1}|x[\k,\n]|^{2}}.\]

Denote by $\mathbf{T}_{\vt}$ the operator of soft thresholding:
\[\mathbf{T}_{\vt}\,\mathbf{x}=\left\{x_{\vt}[\k,\n]\right\},\quad x_{\vt}[\k,\n]\srr \mathrm{sgn}(x[\k,\n])\,\max \left\{0, |x[\k,\n]|-\vt\right\}.\]
Following \cite{ji_shen_xu}, we solve the minimization problem in
\eh{av_dat_min} by an iterative  SBI algorithm. We begin with the
initialization $\mathbf{u}^{0}=0,\;\mathbf{d}^{0}=\mathbf{b}^{0}=0$.
Then,
\begin{equation}\label{breg_iter}
  \begin{array}{l}
    \mathbf{u}^{k+1}:=(\mathbf{K^*}\,\mathbf{P}_{\La}
    \,\mathbf{K} +\mu\,\mathbf{I})\,\mathbf{u}=\mathbf{K^*}\,\mathbf{P}_{\La}\,\mathbf{f}
+\mu\,\mathbf{F}\,(\mathbf{d}^{k}-\mathbf{b}^{k}),\\
    \mathbf{d}^{k+1} =\mathbf{T}_{\la/\mu}(\mathbf{\tilde{F}}\,\mathbf{u}^{k+1}+\mathbf{b}^{k}),\\
    \mathbf{b}^{k+1}=\mathbf{b}^{k}+(\mathbf{\tilde{F}}\,\mathbf{u}^{k+1}- \mathbf{d}^{k+1}).
  \end{array}
\end{equation}
The linear system in the first line of \eh{breg_iter} is solved by
the application of the \emph{conjugate gradient} algorithm. The
operations in the second and third lines are straightforward. The
choice of the parameters $\la$ and $\mu$ depends on experimental
conditions.
\subsubsection{Examples}\label{sec:sss722}
\begin{description}
  \item[Example I: ``Barbara" blurred, missing 50\% of pixels:] The ``Barbara" image was restored after it was blurred by a \cc\ with the Gaussian kernel (MATLAB \fd\ \\
\texttt{fspecial}('gaussian',[5 5])) and its PSNR became 23.32 dB.
Then,  50\% of its pixels were randomly removed. This reduced the
PSNR to 7.56 dB. Random noise was not added.
 The image was restored by 50 SBI using the parameters $\la=0.0015,\;\mu=0.00014$ in \eh{breg_iter}. The conjugate gradient solver used 150 iterations.
 qWPs originating from \ds s of sixth order were used. For ``bases", 8-samples shifts of  all the WPs from the third \d\ level were selected.  Matlab implementation of the restoration procedures took 59.6 seconds.

Figure \ref{barb50B} displays the restoration result. The image is  deblurred and the fine  texture is restored  completely with PSNR=32.09 dB. Note that the best result in an identical experiment reported in  \cite{ANZ_book1} achieved PSNR=30.32 dB.
  \begin{figure}[H]%[ht!]
 \begin{center}
 \resizebox{14cm}{9cm}{
 \includegraphics{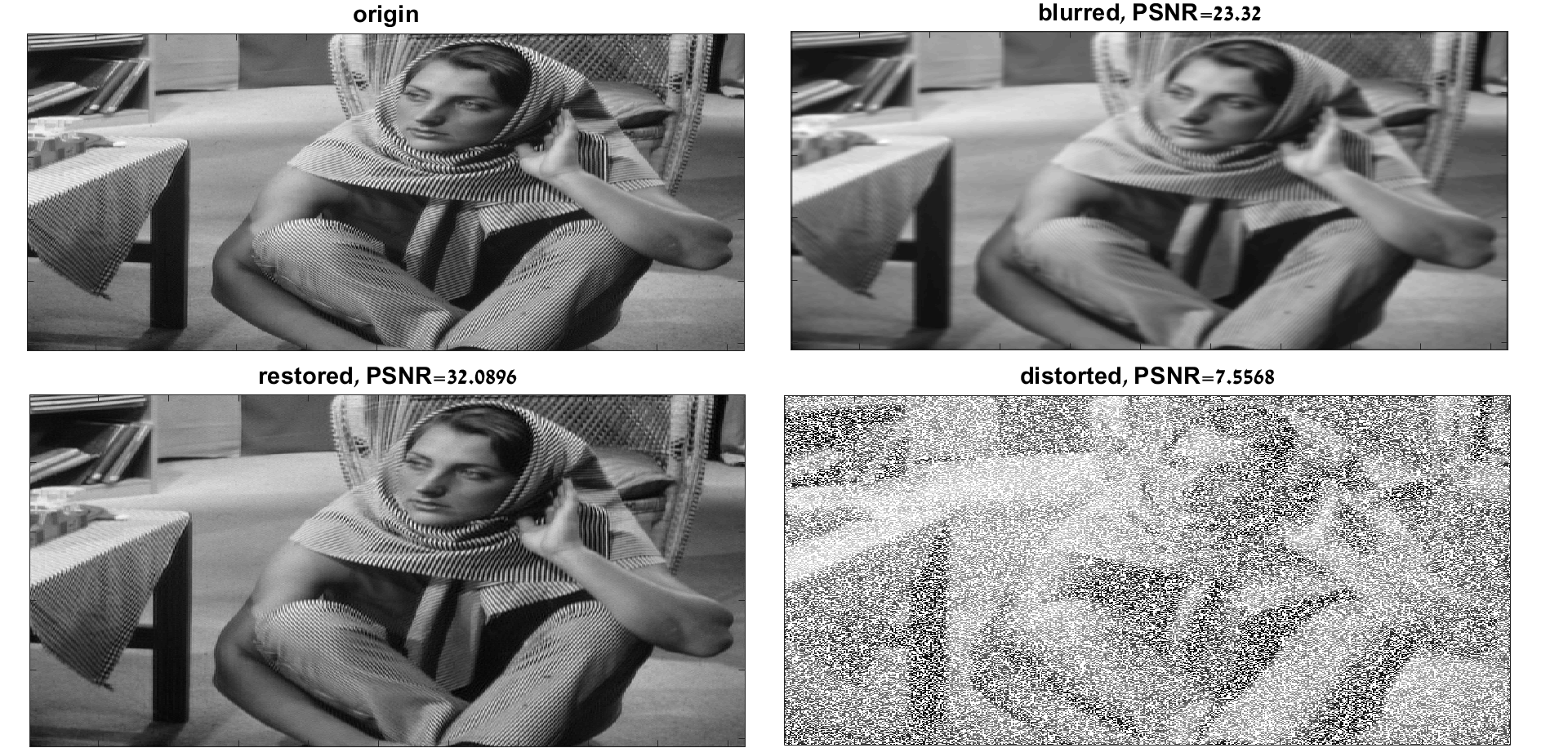}
 }
    \caption {Top left: Source input - ``Barbara" image. Top right: Blurred, PSNR=23.32 dB. Bottom right: After random removal of 50\% of its pixels. PSNR=7.56 dB.
    Bottom left: The  image restored by the directional qWPT. PSNR=32.09 dB}\label{barb50B}
\end{center}
\end{figure}
\item[Example II: ``Barbara" blurred, added noise, missing 50\% of pixels:] The ``Barbara" image was restored after it was blurred by a \cc\ with the Gaussian kernel (MATLAB \fd\
\texttt{fspecial}('gaussian',[5 5])). Random Gaussian noise with STD=10 dB was  added  and the image  PSNR became 22.08 dB.
Then,  50\% of its pixels were randomly removed. This reduced the PSNR to 7.53 dB. The image was restored by 70 SBI using the parameters $\la=3,\;\mu=0.025$ in \eh{breg_iter}. The conjugate gradient solver used 15 iterations.
 qWPs originating from \ds s of fourth order were used. For the ``bases", 16-samples shifts of  all the WPs from the fourth \d\ level were selected.  Matlab implementation of the restoration procedures took 51.9 seconds.

Figure \ref{barb50BN} displays the restoration result. The image is  deblurred, noise is removed  and the fine  texture is partially restored producing PSNR=24.31 dB. Note that the best result in an identical experiment reported in  \cite{ANZ_book1} achieved PSNR=24.19 dB.
 \begin{figure}[H]%[ht!]
 \begin{center}
 \resizebox{14cm}{9cm}{
 \includegraphics{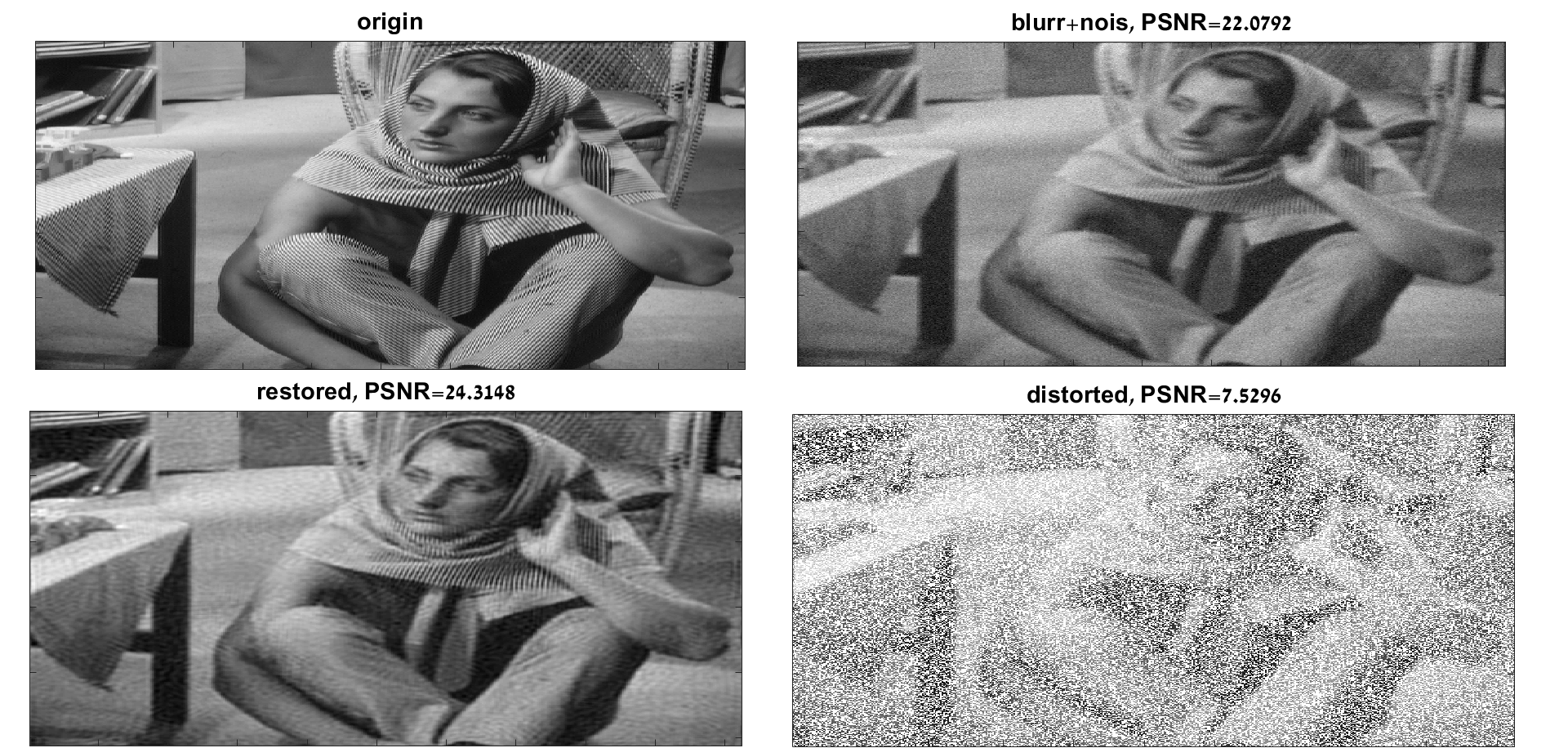}
 }
    \caption {Top left: Source input - ``Barbara" image. Top right: Blurred and noised, PSNR=22.08 dB. Bottom right: After random removal of 50\% of its pixels. PSNR=7.53 dB.
    Bottom left: The  restored image by the directional qWPT. PSNR=24.31 dB}\label{barb50BN}
\end{center}
\end{figure}
\item[Example III: ``Barbara" blurred, missing 90\% of pixels:] The ``Barbara" image was restored after it was blurred by a \cc\ with the Gaussian kernel (MATLAB \fd\ \\
\texttt{fspecial}('gaussian',[5 5])) and  90\% of its pixels were randomly removed. This reduced the PSNR to 5.05 dB. The image was restored by 150 SBI using the parameters $\la=0.0025,\;\mu=0.000025$ in \eh{breg_iter}. The conjugate gradient solver used 150 iterations.
 qWPs originating from \ds s of eighth order were used. For the ``bases", 16-samples shifts of  all the WPs from the fourth \d\ level were selected.  Matlab implementation of the restoration procedures took 226.8 seconds.

Figure \ref{barb90} displays the restoration result. The image is  deblurred and the fine  texture is partially restored. The output has PSNR=25.24 dB.
 \begin{figure}[H]%[ht!]
 \begin{center}
 \resizebox{14cm}{9cm}{
 \includegraphics{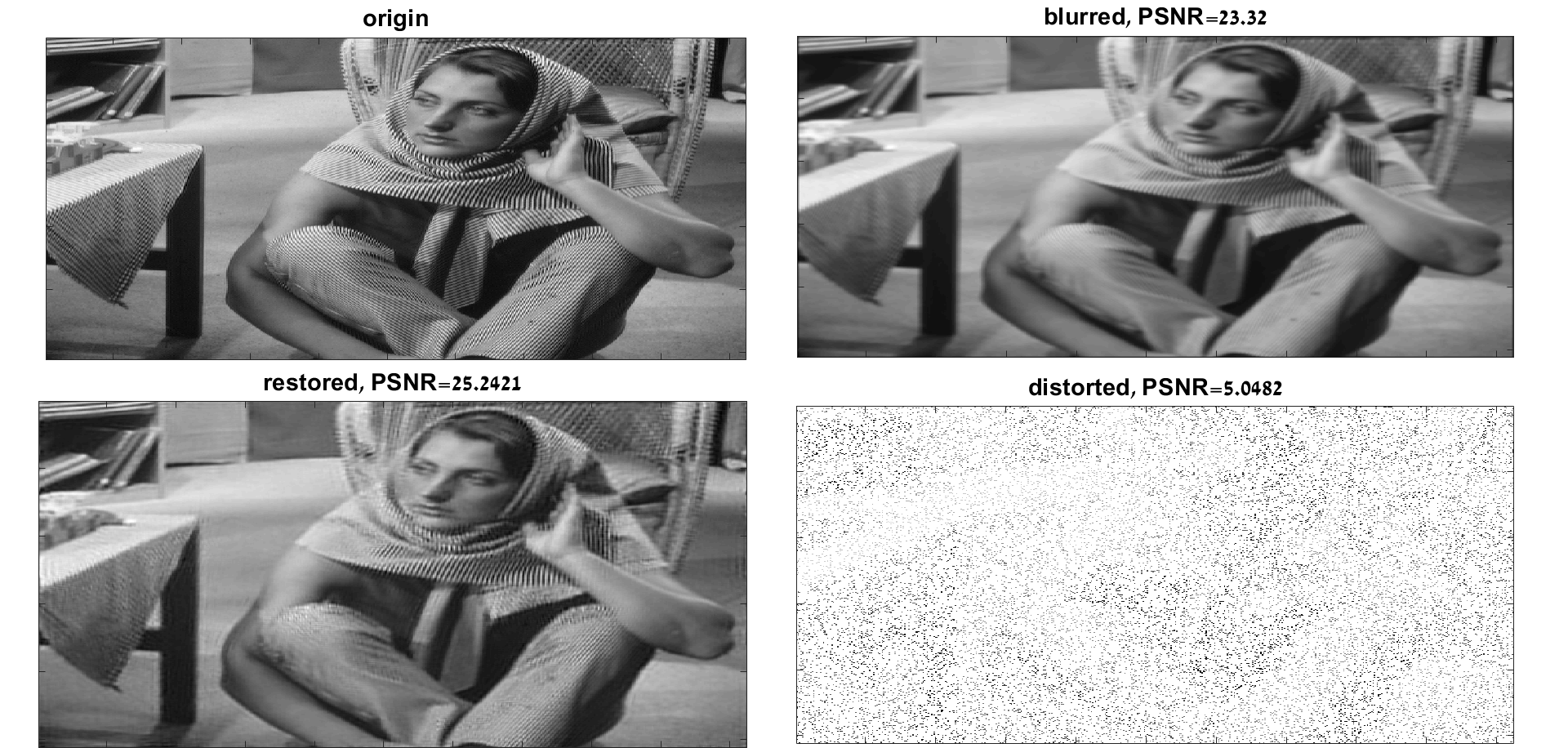}
 }
    \caption {Top left: Source input - ``Barbara" image. Top right: Blurred, PSNR=23.32 dB. Bottom right: After random removal of 90\% of its pixels. PSNR=5.05 dB.
    Bottom left: The  image restored by the directional qWPT. PSNR=25.24 dB}\label{barb90}
\end{center}
\end{figure}

\end{description}
\section{Discussion}\label{sec:s8}We presented a library of complex \dt\ \wq s operating in one-  or two-dimensional spaces of \p\ \ss s. Seemingly, the requirement of \p ity imposes some limitations on the scope of \ss s available for \pr, but actually these limitations are easily circumvented. Any limited \ss\ can be regarded as one period of a \p\ \ss. In order to prevent boundary effects, the \ss s can be \sy ally extended beyond the boundaries before \pr\ and shrunk to the original size after that.  We used such a trick in the ``Barbara" denoising examples.

On the other hand, the \p\ setting provides a lot of substantial opportunities for the design and implementation of   WP \t s such as
\begin{itemize}
  \item A unified computational scheme based on 1D and 2D FFT.
  \item Opportunity to use \fr s with infinite \ir s, which enables us to design a variety of \on\ WP systems where WPs can have any number of local \vm s.
  \item The number of local \vm s does not affect the computational cost of the \t s implementation.
\item A simple explicit scheme of expansion of real  WPs to \az\ and quasi-\az\ WPs with perfect \ff\ separation.
\end{itemize}
The library of qWP \t s described in the paper has a number of free parameters enabling to adapt the \t s to the problem under consideration:
\begin{itemize}
  \item Order of the generating \s, which determines the  number of local \vm s.
  \item Depth of \d, which in 2D case determines the    \de ality of qWPs. For example, fourth-level qWPs are oriented in 62 \df t \de s.
  \item Selection of an optimal structure, such as, for example,  separate Best Bases in   the real and imaginary parts of 1D qWP \t s, separate ``Best Bases" in positive and negative branches of 2D dual-tree qWP \t s, a \ww-basis structure or the set of all \wq s from a single level.
\item Controllable redundancy rate of the \ss\ \ry ation. The minimal rate is 2 when one of options listed in a previous item is utilized. However, several basis-type structures can be involved, for example, all  \wq s from several  levels can be used for the \ss\ \r\ and results can be averaged.
\end{itemize}
The goal of the paper is to design qWPs with an efficient computational scheme for the corresponding \t s. A few experimental results   highlight exceptional properties of these WPs. The directional qWPs are tested for image restoration examples. In the denoising examples, the goal was not to achieve the best output in RSNR values but rather to compare the performance of \de al versus standard WPs that are tensor-product-based. We did not use sophisticated adaptive denoising schemes such as for example  Gaussian scale mixture model (\cite{porti_strela}) or bivariate shrinkage (\cite{sen_seles}). Instead, after \d\ of an image down to level $M$, we selected ``Best Bases" in the positive and negative branches of the qWP \t\ and in the WP \t. Then, we  discarded all the \c s in these bases except for the $L$ largest \c s.  Even with such naive scheme, the \de al  qWPs significantly outperform  the standard WPs in both PSNR values and in the visual perception. The edges, lines and oscillating texture structures in the ``Barbara" image, which was corrupted by Gaussian noise with STD=30 dB, were almost perfectly restored. When  noise STD was 50 dB, an essential part of these structures was restored as well. We emphasise that the Matlab implementation of all the procedures including 3-level qWP and WP \t s, the design of ``Best Bases", thresholding of the \c s arrays and inverse \t s  took 0.88 seconds.  To eliminate  boundary effects, the image was extended from $512\times512$ to $1024\times1024$ and the \pr\ took 3.5 seconds.

The second group of experimental results  dealt with the restoration of the ``Barbara" image which was blurred my the convolving the image with a  Gaussian kernel and degraded  by removing randomly either 50\% or 90\% of the pixels. The image was restored by using a constrained $l_{1}$ minimization of the qWP \t\ \c s from a certain \d\ level and implemented via the split Bregman Iterations procedure.
In  Example I with missing 50\%  of the pixels, the image was almost perfectly restored with PSNR-32.1 dB and practically all the fine structure reconstructed although it was blurred even before the removal of the pixels.
Addition of the Gaussian noise with STD 10 dB to the blurred image in Example II depleted the \r\ result. Although  the noise became suppressed and the image was deblurred, most of the fine structure was lost. Restoration results were better in  Example III where, instead of adding noise, the number of  pixels missing from the blurred image was raised to 90\%. The image was deblurred and an essential part of fine structure was restored.

Summarizing, we can state that, having such a versatile and flexible tool at hand, we are prepared to address multiple data \pr\ problems such as \ss\ and image deblurring and denoising, target detection, segmentation, inpainting,  superresolution, to name a few. In one of the applications, whose results are to be  published soon, \de al qWPs are used with Compressed Sensing methodology  for the   conversion of a regular digital photo camera to an hyper\sp al imager. Preliminary results appear in \cite{Hyper_camera}. Special efforts are on the way for the  design of a denoising scheme which fully explores the  \de ality properties of qWPs.

\pa{Acknowledgment}
This research was partially supported by the Israel Science Foundation (ISF, 1556/17),
Blavatnik Computer Science Research Fund
Israel Ministry of Science and Technology 3-13601 and by Academy of Finland (grant 311514).

\include{ANA_AP}
\bibliographystyle{plain}
   % mathematics and physical sciences
\bibliography{BookBib_SA}
\end{document}